\documentclass[20pt,reqno]{article}
\usepackage{float}
\usepackage[T1]{fontenc}
\usepackage{caption}
\usepackage{placeins}
\usepackage{amsfonts,mathrsfs}
\usepackage{fancyhdr,amscd}
\usepackage{anysize}
\usepackage{setspace}
\usepackage{graphicx,epsfig}
\usepackage{amsthm,latexsym,amsmath,amssymb}
\usepackage[perpage,symbol*]{footmisc}
\usepackage[colorlinks=true]{hyperref}
\usepackage{xcolor}
\usepackage{amsfonts,amsmath,amssymb,latexsym}
\usepackage{enumerate}
\usepackage{mathrsfs}
\usepackage{CJKfntef}
\usepackage{epsfig,CJK,fancyhdr}
\usepackage{CJK,fancyhdr,amscd}
\usepackage{indentfirst}
\newtheorem{thm}{Theorem}[section]
\newtheorem{defi}[thm]{Definition}
\newtheorem{prop}[thm]{Proposition}
\newtheorem{lem}[thm]{Lemma}
\newtheorem{cor}[thm]{Corollary}

\newtheorem{rem}[thm]{Remark}
\numberwithin{equation}{section}

\def\be{\begin{equation}}
\def\ee{\end{equation}}
\def\bea{\begin{eqnarray}}
\def\eea{\end{eqnarray}}

\def\y{\begin{eqnarray*}}
\def\ey{\end{eqnarray*}}

\allowdisplaybreaks[4]
\title{Explicit resonance and bifurcation of symmetric periodic solutions for differential systems with multiple distributed delays}

\author{Huiting He$^{a}$, Chungen Liu$^{a,}$\footnote{Corresponding author. \newline
		The second author of this paper is supported by the National Natural Science Foundation of China (Grant No.11790271, 12171108).	\newline
			Email addresses: \newline
			Heht21@163.com(H.He),\ liucg@nankai.edu.cn(C.Liu).} 
			\\
	\small  $^a$ School of Mathematics and Information Science, Guangzhou University, Guangzhou, 510006, PR China\\
	} 
\date{}
\begin{document}
\maketitle
\begin{abstract} 
	We investigate the bifurcation of symmetric periodic solutions for differential systems with multiple distributed delays. For general multiple-delay systems with a quadratic potential, we establish a Hamiltonian formulation with a \(P_{n,m}\)-boundary condition and derive explicit Maslov index and nullity formulas. For the special case of two distributed delays, we extend the analysis to general nonlinear potentials. These results yield explicit resonance criteria, bifurcation conditions, and Rabinowitz type	and Fadell-Rabinowitz type alternatives for $\mathbb{R} $-distinct solutions.
\vspace{0.6cm}

\end{abstract}
\maketitle
\section{Introduction and main results}\label{sec1}
Consider the following autonomous delay differential  system involving multiple distributed delays with a parameter $\lambda $ in a  real interval $\Lambda $ and a fixed delay $\tau > 0$
\begin{equation}\label{e1}
\begin{cases}
	\dot{x}(t)=-\nabla F\left ( \lambda , \int_{0}^{\tau } \nabla G\left ( \lambda ,x\left ( t-s \right ) \right )   ds \right )- \cdots-\nabla F\left ( \lambda , \int_{\left ( m-1 \right ) \tau }^{m\tau } \nabla G\left (\lambda , x\left ( t-s \right ) \right )  ds \right ),\\
	x\left ( t+m\tau  \right ) =-x\left ( t \right ),\quad  t\in \mathbb{R}, \ x(t)\in \mathbb{R}^{n},
\end{cases}
\end{equation}
where  $F, G :\Lambda \times \mathbb{R}^{n}\to \mathbb{R}$ are  continuous functions satisfying the following condition:
\begin{enumerate}[$({\bf FG})$]
	\item  For each fixed $\lambda \in \Lambda$, the functions $F(\lambda ,\cdot )$ and $G(\lambda ,\cdot )$ are of class $C^{2}$ and are even, namely, 
	\begin{equation}
		F(\lambda, -x )=F(\lambda,x), \quad G(\lambda, -x )=G(\lambda,x), \quad \forall  x\in \mathbb{R}^{n}.
	\end{equation}
	Moreover, the gradients $\nabla F(\lambda ,x)$ and $\nabla G(\lambda ,x)$ with respect to $x$, as well as the corresponding Hessian matrices $\nabla^{2} F(\lambda ,x)$ and $\nabla^{2} G(\lambda ,x)$,  depend continuously on $(\lambda ,x)\in \Lambda \times \mathbb{R} ^{n}$.
\end{enumerate}
A solution $x(t)$ of system \eqref{e1} satisfying the above anti-periodic boundary condition is called a \textit{symmetric $2m\tau$-periodic solution}.

The study of symmetric periodic solutions of delay differential equations dates back to the pioneering work of Kaplan and Yorke \cite{K}, who reduced the search for periodic solutions of certain scalar delay equations to the study of associated ordinary differential systems. Their approach also motivated the investigation of equations with multiple discrete delays, and Nussbaum  \cite{Nussbaum1978} subsequently established a generalization of the Kaplan-Yorke conjecture. Building on these ideas, Hamiltonian and variational methods for symmetric periodic solutions were developed in \cite{Li3,Li1,Li2}. Liu \cite{Liu3} introduced the Maslov $P$-index theory for Hamiltonian systems with $P$-boundary conditions and applied it to discrete delay systems. Further developments and applications of this theory can be found  in \cite{Liu10,Liu2,Wang3,Zhou}.  Alternative approaches based on critical point theory and equivariant index methods were  also explored, see for instance \cite{Fei1,Fei2,Guo1,Guo2,Zheng}.

Compared with discrete delay systems, distributed delay systems present additional difficulties in the construction of suitable Hamiltonian structures. Xiao and Guo \cite{Xiao1} developed a Hamiltonian reduction for a class of distributed delay systems, which was subsequently extended to higher-dimensional systems in \cite{Xiao2}.
Further progress was made by Zhong, Wang, and Liu in \cite{Liu8,Zhong,Zhong2}, who established existence and multiplicity results for symmetric periodic solutions of systems with a single distributed delay. More recently, the authors in this paper \cite{He2,He} investigated symmetric periodic solutions for systems involving multiple distributed delay intervals. These developments provide an important foundation for further investigations of the qualitative behavior of symmetric periodic solutions for distributed delay systems.
 
Beyond existence and multiplicity, the bifurcation of symmetric periodic solutions is another important problem. Classical approaches to bifurcation in functional differential systems include global bifurcation theory \cite{Nussbaum1975}, the Fuller index methods \cite{Chow}, equivariant degree methods for symmetric Hopf bifurcation \cite{Erbe,Krawcewicz,Wu}, and techniques based on periodic normal forms and Lyapunov-Schmidt reduction \cite{Guo,Lentjes}. From a Hamiltonian system perspective, Lu \cite{Lu1}  developed a bifurcation theory based on the Maslov-type index for  systems with generalized boundary conditions, yielding bifurcation criteria and, under additional assumptions, Rabinowitz type and Fadell-Rabinowitz type bifurcation alternatives. Related developments for potential-type delay equations can be found in \cite{Lu,Lu2}.

Motivated by these developments, we investigate the bifurcation of symmetric $2m\tau$-periodic solutions for system \eqref{e1}. Our approach combines Hamiltonian reduction with Maslov $P$-index theory. More precisely, we consider two complementary settings, distinguished by the structural assumptions imposed on the potential $G$.

First, we consider systems with multiple distributed delay intervals under the assumption that $G(\lambda,\cdot)$ is quadratic. Specifically, we impose the following assumption:
\begin{enumerate}[$({\bf G})$]
	\item There exists a continuous mapping $\Theta :\Lambda \to \mathcal{L} _{s}(\mathbb{R}^{n} )\cap \mathrm{GL}(n,\mathbb{R} )$ such that 
	\begin{equation}
		\nabla G(\lambda,x)= \Theta(\lambda)x, \quad \forall (\lambda ,x)\in \Lambda \times \mathbb{R} ^{n}.
	\end{equation}
	Here, $\mathcal{L} _{s}(\mathbb{R}^{n} )$ denotes the space of all $n\times n$ symmetric matrices, while $\mathrm{GL}(n,\mathbb{R} )$ denotes the group of all $n\times n$ invertible matrices.
\end{enumerate}
Under assumption $({\bf G})$, system \eqref{e1} reduces to 
\begin{equation}\label{e2}
	\begin{cases}
		\dot{x}(t)=-\nabla F\left ( \lambda , \int_{0}^{\tau } \Theta(\lambda)x\left ( t-s \right )   ds \right )- \cdots-\nabla F\left ( \lambda , \int_{\left ( m-1 \right ) \tau }^{m\tau } \Theta(\lambda) x\left ( t-s \right )   ds \right ),\\
		x\left ( t+m\tau  \right ) =-x\left ( t \right ),\quad  t\in \mathbb{R}, \ x(t)\in \mathbb{R}^{n}.
	\end{cases}
\end{equation}
Two solutions $x$ and $y$ of system \eqref{e2} with a fixed $\lambda$ are called $\mathbb{R} $-distinct if and only if there is no $s\in \mathbb{R} $ such that $x=s\cdot y$, where the $\mathbb{R} $-action  is given by $s\cdot y(t)=y(t+s)$ for all $t\in \mathbb{R}$. Since $F(\lambda ,\cdot )$ and $G(\lambda ,\cdot )$ are even, it follows that 
\begin{equation}
	\nabla F(\lambda ,0)=\nabla G(\lambda ,0)=0, \quad \forall \lambda\in \Lambda.
\end{equation}
Consequently, $x\equiv0$ is a solution
of \eqref{e2} for every $\lambda\in\Lambda$,
giving rise to the trivial branch $\left \{ \left ( \lambda ,0 \right )|\lambda \in \Lambda   \right \} $. A point $(\mu,0)$ is called
a bifurcation point of \eqref{e2}
if there exists a sequence
$\left \{ \lambda _{k} \right \} _{k=1}^{\infty }\subset \Lambda$ with
$\lambda_k\to\mu$ and nonzero
solutions $x_k$ of \eqref{e2}
corresponding to $\lambda=\lambda_k$
such that
\[
\|x_k\|_{C^1([0,\tau],\mathbb{R}^n)}
\to 0
\quad\text{as }k\to\infty.
\]

Under assumption $({\bf G})$, we transform the original delay system into a first-order Hamiltonian system with a $P_{n,m}$-boundary condition. Explicit formulas for the associated  Maslov $P_{n,m}$-index and nullity are established in Lemma \ref{abv} for scalar systems and in Lemma \ref{abv2} for higher-dimensional systems with commuting symmetric matrices. These formulas yield resonance conditions, bifurcation criteria, and under suitable crossing assumptions, Fadell-Rabinowitz type and Rabinowitz type alternatives. The corresponding results are stated in Theorems \ref{n1} and \ref{thm2}.

Second, we consider the special case of two distributed delay intervals
\begin{equation}\label{e3'}
	\left\{\begin{matrix}
		\dot{x}(t)=-\nabla F\left ( \lambda ,\int_{0}^{\tau }\nabla G\left ( \lambda ,x\left ( t-s \right ) \right ) ds \right )-\nabla F\left ( \lambda , \int_{\tau }^{2\tau }\nabla G\left ( \lambda ,x\left ( t-s \right ) \right )ds \right ), \\
		x\left ( t+2\tau  \right ) =-x\left ( t \right ),\quad  t\in \mathbb{R}, \ x(t)\in \mathbb{R}^{n},\quad \quad \quad \quad \quad \quad \quad \quad \quad \quad \quad \quad \quad \quad \quad \ \ 
	\end{matrix}\right.
\end{equation}
without requiring $G$ to be quadratic. By exploiting the anti-periodic symmetry, we construct an equivalent first-order Hamiltonian system with a $P_{n,2}$-boundary condition for general even nonlinear potentials $F$ and $G$. Applying the explicit $P_{n,2}$-index and nullity formulas established in Lemmas \ref{abv} and \ref{abv2} with $m=2$, we derive resonance conditions and bifurcation criteria in both scalar and higher-dimensional settings. We also establish the corresponding Fadell-Rabinowitz type and Rabinowitz type bifurcation alternatives, see Theorems \ref{th4} and \ref{thm5}.

 Our main results are stated below.

\begin{figure}[!htbp]
	\centering
	\includegraphics[width=\linewidth]{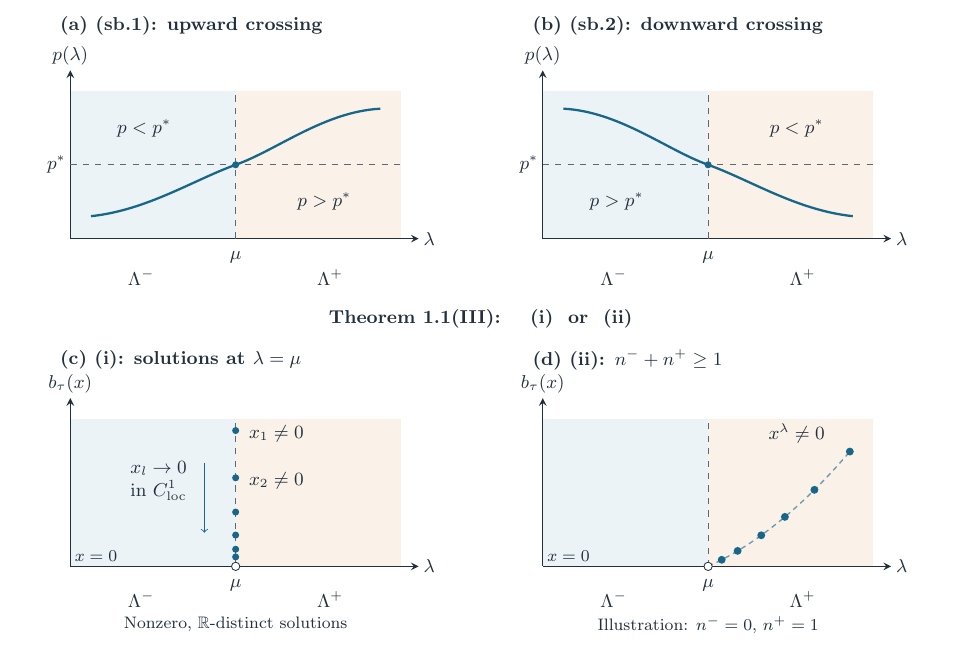}
	\caption{Fadell--Rabinowitz type alternatives.
		The diagrams are schematic: the crossing conditions do not require
		$p$ to be monotone, and no continuous dependence of the solution
		families on $\lambda$ is asserted.}
	\label{fig:HL-thm11}
\end{figure}

\begin{thm}\label{n1}
	Under assumption $({\bf FG})$ and $({\bf G})$, with $n=1$,  delay differential system \eqref{e2} with multiple distributed delays becomes
	\begin{equation}\label{fx}
		\begin{cases}
			\dot{x}(t)=-F'  \left ( \lambda , \theta (\lambda)\int_{0}^{\tau }    x\left ( t-s \right )  ds \right )- \cdots-F' \left ( \lambda ,  \theta (\lambda)\int_{\left ( m-1 \right ) \tau }^{m\tau } x\left ( t-s \right )  ds \right ),\\
			x\left ( t+m\tau  \right ) =-x\left ( t \right ),\quad  t\in \mathbb{R}, \ x(t)\in \mathbb{R}, 
		\end{cases}
	\end{equation}
	where $\theta :\Lambda\to  \mathbb{R}\setminus \left \{ 0 \right \} $ is continuous. Here, $F'(\lambda ,\cdot )$ and $F''(\lambda ,\cdot )$ denote, respectively,  the first and second derivatives of $F(\lambda ,\cdot )$ with respect to the second variable. Let $f(\lambda )=F'' (\lambda ,0)$, then 
	\begin{enumerate}[{\rm (I)}]
		\item {\rm (Necessary condition):} If $\left ( \mu ,0 \right )$ with $\mu \in \Lambda $ is a bifurcation point of system \eqref{fx}, then $m\sqrt{2f(\mu  )\theta (\mu  )}\tau \in  \left ( 2\mathbb{N}-1  \right ) \pi $.
		\item {\rm (Sufficient condition):} Let $\mu$ be an interior point of $\Lambda$ satisfying the following conditions: 
		\begin{enumerate}[{\rm (a)}]
			\item $f(\mu  )\theta (\mu  )>0$ and $m\sqrt{2f(\mu  )\theta (\mu  )}\tau \in  \left ( 2\mathbb{N}-1  \right ) \pi $.
			\item There exist sequences $\left ( \lambda _{q}^{\pm } \right ) _{q}\subset \Lambda $ converging to $\mu$ such that $f(\lambda _{q}^{-} )\theta (\lambda _{q}^{-} )< f(\mu  )\theta (\mu  )< f(\lambda _{q}^{+} )\theta (\lambda _{q}^{+} )$ for all $q\in \mathbb{N} $.
		\end{enumerate}
		Then, $\left ( \mu ,0 \right )$ is a bifurcation point of system \eqref{fx}.
		\item {\rm (Alternative bifurcations of Fadell-Rabinowitz type):} Let the assumptions of {\rm (II)} hold, with condition {\rm (b)} replaced by the stronger condition:
		\begin{enumerate}[{\rm (sb)}]
			\item For some $\varepsilon >0$, one of the following two conditions holds:
			\begin{enumerate}[{\rm (sb.1)}]
				\item $f(\lambda )\theta (\lambda )< f(\mu  )\theta (\mu  )< f(\lambda'  )\theta (\lambda' )$ for all $\left ( \lambda ,\lambda ' \right )\in \left ( \mu-\varepsilon ,\mu  \right )\times \left ( \mu ,\mu +\varepsilon \right ) $.
				\item $f(\lambda )\theta (\lambda )> f(\mu  )\theta (\mu  )> f(\lambda'  )\theta (\lambda' )$ for all $\left ( \lambda ,\lambda ' \right )\in \left ( \mu-\varepsilon ,\mu  \right )\times \left ( \mu ,\mu +\varepsilon \right ) $.
			\end{enumerate}
		\end{enumerate}
		Then, for system \eqref{fx}, there exist the following alternatives:
		 \begin{enumerate}[{\rm (i)}]
		 	\item The system \eqref{fx} with $\lambda=\mu$ admits a sequence of nonzero $\mathbb{R} $-distinct solutions $\left \{ x_{l} \right \} _{l}^{\infty }$ such that $x_{l}\to 0$ in $C_{{\rm loc}}^{1}\left ( \mathbb{R},\mathbb{R}  \right ) $ as $l\to \infty $.
		 	\item There exist left and right neighborhoods $\Lambda^{-}$ and $\Lambda^{+}$ of $\mu$ in $\Lambda$ and nonnegative integers $n^{+}$ and $n^{-}$ satisfying $n^{+}+n^{-}\ge 1$, such that for $\lambda \in \Lambda ^{-}\setminus \left \{ \mu  \right \} (resp. \ \lambda \in \Lambda ^{+}\setminus \left \{ \mu  \right \})$, system \eqref{fx} with parameter $\lambda$ has at least $n^{-}(resp.\ n^{+})$ $\mathbb{R} $-distinct nonzero solutions $x^{\lambda}_{i}(i=1,\cdots,n^{-}(resp.\ n^{+}))$ which converge to zero in $C_{{\rm loc}}^{1}\left ( \mathbb{R},\mathbb{R}  \right ) $ as $\lambda\to \mu$.
		 \end{enumerate}
	\end{enumerate}
\end{thm}
\paragraph{Geometric interpretation.}
Figure~\ref{fig:HL-thm11} illustrates Theorem \ref{n1}(III). Set $p(\lambda)=f(\lambda)\theta(\lambda)$ and $p^*=p(\mu)=(2j-1)^2\pi^2/(2m^2\tau^2)$.
Under the theorem's hypotheses, either strict crossing in $\mathbf{(a)}$  or $\mathbf{(b)}$
yields $\mathbf{(c)}$ or $\mathbf{(d)}$. In the solution diagrams, the vertical coordinate is
$b_\tau(x)=\|x|_{[0,\tau]}\|_{C^1}$. Whenever a sequence or family of solutions converges to zero in $C^1_{\mathrm{loc}}$, the corresponding values  $b_\tau$ tend to zero.

\FloatBarrier
\clearpage
\begin{thm}\label{thm2}
	Suppose that $({\bf FG})$ and $({\bf G})$ hold with $n>1$, then we have 
	\begin{enumerate}[{\rm (I)}]
		\item {\rm (Necessary condition):} If $\left ( \mu ,0 \right )$ with $\mu \in \Lambda $ is a bifurcation point of system  \eqref{e2}, and if the $n\times n$ real symmetric matrices $\nabla^{2}F(\mu ,0)$ and $\Theta (\mu )$ commute, that is, there exists a real orthogonal $n\times n$ matrix $T$ such that 
		\begin{equation}
			T^{-1}\nabla ^{2}F(\mu ,0)T=\mathrm{diag}(f_{1}(\mu ),\cdots ,f_{n}(\mu ))\ \mathrm{and} \ T^{-1}\Theta (\mu )T=\mathrm{diag}(\theta _{1}(\mu),\cdots ,\theta _{n}(\mu )),
		\end{equation}
		then $m\sqrt{2f_{j}(\mu  )\theta_{j} (\mu  )}\tau \in  \left ( 2\mathbb{N}-1  \right ) \pi $ for some $j\in \left \{ 1,\cdots ,n \right \} $.
		\item {\rm (Sufficient condition):} For each $\lambda\in \Lambda$, suppose $\nabla ^{2}F(\lambda ,0)$ and $\Theta (\lambda )$ commute. Let $\mu$ be an interior point of $\Lambda$, assume further that the matrix $\nabla ^{2}F(\mu ,0)$ has $n$ distinct eigenvalues $f_{j}(\mu)$, $j=1,\cdots,n$. 
		Suppose there exists a partition $\left \{ j_{1},\cdots ,j_{k} \right \} \cup \left \{ j_{k+1},\cdots ,j_{n} \right \} $ of $\left \{ 1,\cdots ,n \right \} $ with $1\le k\le n$, such that the following conditions hold,
		\begin{enumerate}[{\rm (A)}]
			\item For each $l=k+1,\cdots,n$, either $2f_{j_{l}}(\mu)\theta _{j_{l}}(\mu)< \left ( \pi /m\tau  \right ) ^{2}$ or 
			\begin{equation}\label{2r}
				\left ( 2\rho _{l}(\mu)-1 \right ) \pi <m\sqrt{2f_{j_{l}}(\mu)\theta _{j_{l}}(\mu)} \tau < \left ( 2\rho _{l}(\mu)+1 \right ) \pi ,
			\end{equation}
			for some $\rho _{l}(\mu)\in\mathbb{N} $.
			\item For each $l=1,\cdots,k$, $f_{j_{l}}(\mu)\theta _{j_{l}}(\mu)> 0$ and $m\sqrt{2f_{j_{l}}(\mu)\theta _{j_{l}}(\mu)}\tau =\left ( 2\rho _{l}(\mu)-1 \right ) \pi  $, $\rho _{l}(\mu)\in\mathbb{N} $.
			\item There exist sequence $\left ( \lambda _{q}^{\pm } \right ) _{q}\subset \Lambda $  converging to $\mu$ such that $f_{j_{l}}(\lambda_{q}^{-} )\theta _{j_{l}}(\lambda_{q}^{-} )<f_{j_{l}}(\mu)\theta _{j_{l}}(\mu)<f_{j_{l}}(\lambda_{q}^{+} )\theta _{j_{l}}(\lambda_{q}^{+} )$ if $ f_{j_{l}}(\mu  )>0$, and $f_{j_{l}}(\lambda_{q}^{-} )\theta _{j_{l}}(\lambda_{q}^{-} )>f_{j_{l}}(\mu)\theta _{j_{l}}(\mu)>f_{j_{l}}(\lambda_{q}^{+} )\theta _{j_{l}}(\lambda_{q}^{+} )$ if $ f_{j_{l}}(\mu  )<0$, for all $q\in \mathbb{N} $ and $l=1,\cdots,k$.
		\end{enumerate}
		then $\left ( \mu ,0 \right )$ is a bifurcation point of system  \eqref{e2}.
		\item {\rm (Alternative bifurcations of Fadell-Rabinowitz type and of Rabinowitz type):} Let the assumptions of {\rm (II)} hold, with {\rm (C)} strengthened as follows,
		\begin{enumerate}[{\rm (SC)}]
			\item	For some $\varepsilon >0$, one of the following two conditions holds,
			\begin{enumerate}[{\rm (SC.1)}]
				\item $f_{j_{l}}(\lambda )\theta _{j_{l}}(\lambda )< f_{j_{l}}(\mu  )\theta_{j_{l}} (\mu  )< f_{j_{l}}(\lambda'  )\theta_{j_{l}} (\lambda' )$ for all $\left ( \lambda ,\lambda ' \right )\in \left ( \mu-\varepsilon ,\mu  \right )\times \left ( \mu ,\mu +\varepsilon \right ) $ if $ f_{j_{l}}(\mu  )>0$, and $f_{j_{l}}(\lambda )\theta _{j_{l}}(\lambda )> f_{j_{l}}(\mu  )\theta_{j_{l}} (\mu  )> f_{j_{l}}(\lambda'  )\theta_{j_{l}} (\lambda' )$ for all $\left ( \lambda ,\lambda ' \right )\in \left ( \mu-\varepsilon ,\mu  \right )\times \left ( \mu ,\mu +\varepsilon \right ) $ if $ f_{j_{l}}(\mu  )<0$, $l=1,\cdots,k$.
				\item $f_{j_{l}}(\lambda )\theta _{j_{l}}(\lambda )> f_{j_{l}}(\mu  )\theta_{j_{l}} (\mu  )> f_{j_{l}}(\lambda'  )\theta_{j_{l}} (\lambda' )$ for all $\left ( \lambda ,\lambda ' \right )\in \left ( \mu-\varepsilon ,\mu  \right )\times \left ( \mu ,\mu +\varepsilon \right ) $ if $ f_{j_{l}}(\mu  )>0$, and $f_{j_{l}}(\lambda )\theta _{j_{l}}(\lambda )< f_{j_{l}}(\mu  )\theta_{j_{l}} (\mu  )< f_{j_{l}}(\lambda'  )\theta_{j_{l}} (\lambda' )$ for all $\left ( \lambda ,\lambda ' \right )\in \left ( \mu-\varepsilon ,\mu  \right )\times \left ( \mu ,\mu +\varepsilon \right ) $ if $ f_{j_{l}}(\mu  )<0$, $l=1,\cdots,k$.
			\end{enumerate}
		\end{enumerate}
		Then for system \eqref{e2}, there exist two possible alternatives:
		\begin{enumerate}[{\rm (i)}]
			\item The system \eqref{e2} with $\lambda=\mu$ admits a sequence of nonzero $\mathbb{R} $-distinct solutions $\left \{ x_{l} \right \} _{l}^{\infty }$ such that $x_{l}\to 0$ in $C_{{\rm loc}}^{1}\left ( \mathbb{R},\mathbb{R}^{n}  \right ) $ as $l\to \infty $.
			\item There exist left and right neighborhoods $\Lambda^{-}$ and $\Lambda^{+}$ of $\mu$ in $\Lambda$ and nonnegative integers $n^{+}$ and $n^{-}$ satisfying $n^{+}+n^{-}\ge k$, such that for $\lambda \in \Lambda ^{-}\setminus \left \{ \mu  \right \} (resp. \ \lambda \in \Lambda ^{+}\setminus \left \{ \mu  \right \})$, system \eqref{e2} with parameter $\lambda$ has at least $n^{-}(resp.\ n^{+})$ $\mathbb{R} $-distinct nonzero solutions $x_{i}^{\lambda}(i=1,\cdots,n^{-}(resp.\ n^{+}))$ which converge to zero in $C_{{\rm loc}}^{1}\left ( \mathbb{R},\mathbb{R}^{n}  \right ) $ as $\lambda\to \mu$.
			\item[]
			\hspace{-\labelwidth}\hspace{-\labelsep}
			Moreover, if $k>1$, then at least one of {\rm (i)} and the following {\rm (iii)} and {\rm (iv)} holds,
			\item For every $\lambda \in \Lambda \setminus \left \{ \mu  \right \} $ near $\mu$, there is a nonzero solution $x^{\lambda}$ of system \eqref{e2} with parameter $\lambda$, such that   $x^{\lambda}\to 0$ in $C_{{\rm loc}}^{1}\left ( \mathbb{R},\mathbb{R}^{n}  \right ) $ as $\lambda\to \mu$.
			\item  For a given $\varepsilon >0$, there exists a one-sided neighborhood $\Lambda^{0}$ of $\mu$ in $\Lambda$ such that, for any $\lambda\in \Lambda ^{0}\setminus \left \{ \mu \right \} $, system \eqref{e2} with parameter $\lambda$ has either infinitely many nonzero $\mathbb{R} $-distinct solutions $\left \{ x_{l }^{\lambda} \right \}_{l=1}^{\infty } $ such that $\left \| x_{l}^{\lambda} \mid _{\left [ 0,\tau \right ] }\right \| _{C^{1}}<\varepsilon $ for all $l\in \mathbb{N} $, or at least two nonzero $\mathbb{R} $-distinct solutions $\bar{x} ^{\lambda } $ and $\hat{x} ^{\lambda } $ such that $\left \| \bar{x}^{\lambda } \mid _{\left [ 0,\tau \right ] }\right \| _{C^{1}}<\varepsilon $, $\left \| \hat{x}^{\lambda } \mid _{\left [ 0,\tau \right ] }\right \| _{C^{1}}<\varepsilon $, and 
			\begin{equation}\label{c}
				\begin{split}
					&-\frac{1}{4} \int_{0}^{m\tau } \left ( \Theta (\lambda )\left ( \bar{x}^{\lambda } (t)-\bar{x}^{\lambda }  (t-\tau ) \right ) ,\bar{x}^{\lambda }  (t)-\bar{x}^{\lambda }  (t-\tau )  \right )dt\\
					&+\int_{0}^{m\tau }F\left ( \lambda ,\int_{0}^{\tau }\Theta (\lambda )\bar{x}^{\lambda }  (t-s)ds  \right ) dt   \\
					\ne& -\frac{1}{4} \int_{0}^{m\tau } \left ( \Theta (\lambda )\left ( \hat{x}^{\lambda } (t)-\hat{x}^{\lambda }  (t-\tau ) \right ) ,\hat{x}^{\lambda }  (t)-\hat{x}^{\lambda }  (t-\tau )  \right )dt\\
					&+\int_{0}^{m\tau }F\left ( \lambda ,\int_{0}^{\tau }\Theta (\lambda )\hat{x}^{\lambda }  (t-s)ds  \right ) dt.   
				\end{split}
			\end{equation}						
		\end{enumerate}
	\end{enumerate}
\end{thm}
\begin{rem} Since the continuous, real symmetric $n\times n$ matrix-valued functions $\nabla ^{2}F(\lambda ,0)$ and $\Theta (\lambda )$ commute for all $\lambda\in \Lambda	$, 
	 Lemma 6.2 in \cite{Lu} assures that in a neighborhood $\Lambda_{0}$ of $\mu$ in $\Lambda$ there exist continuous real functions $f_{j}, \theta_{j}$ for $j=1,\cdots,n$, and a continuous real orthogonal $n\times n$ matrix-valued function $T$ such that for each $\lambda\in \Lambda_{0}$,
	\begin{equation}
		T(\lambda )^{-1}\nabla ^{2}F(\lambda ,0)T(\lambda )=\mathrm{diag}(f_{1}(\lambda ),\cdots ,f_{n}(\lambda ))\ \mathrm{and} \ T(\lambda )^{-1}\Theta (\lambda )T(\lambda )=\mathrm{diag}(\theta _{1}(\lambda ),\cdots ,\theta _{n}(\lambda )).
	\end{equation}
\end{rem}
\paragraph{Geometric interpretation.} Figure~\ref{fig:HL-thm12-crossings} illustrates Theorem \ref{thm2}(III). For any $l\in \left \{ 1,\cdots,k \right \} $, denote
$p_l(\lambda)=f_{j_l}(\lambda)\theta_{j_l}(\lambda)$ and $p_l^*=p_l(\mu)=(2\rho_l(\mu)-1)^2\pi^2/(2m^2\tau^2)$.
\begin{figure}[!htbp]
	\centering
	\includegraphics[width=\linewidth]{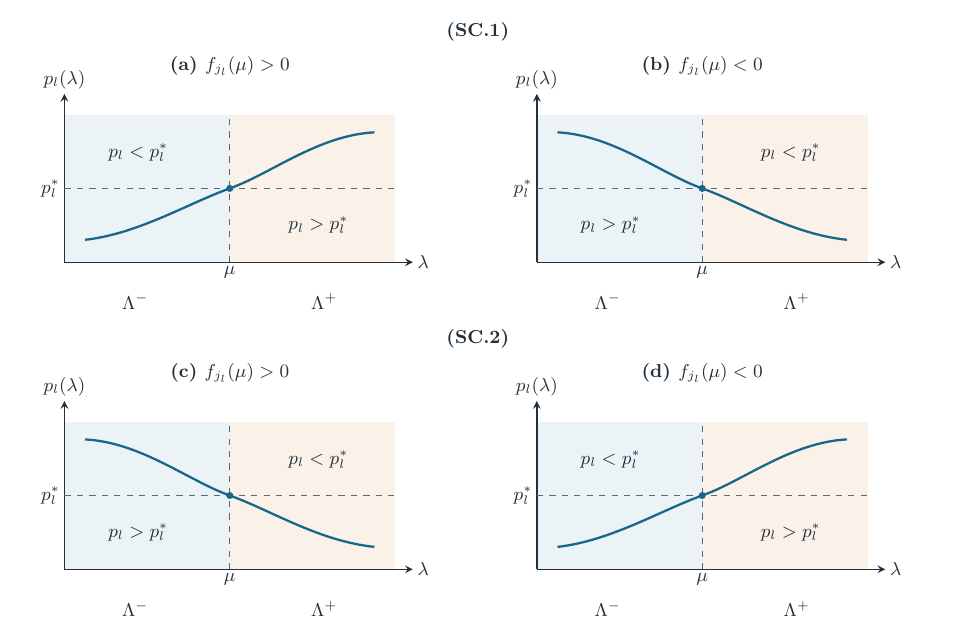}
	\caption{Admissible crossings in Theorem \ref{thm2}(III).}
	\label{fig:HL-thm12-crossings}
\end{figure}
\FloatBarrier

\begin{figure}[!htbp]
	\ContinuedFloat
	\centering
	\includegraphics[width=\linewidth]{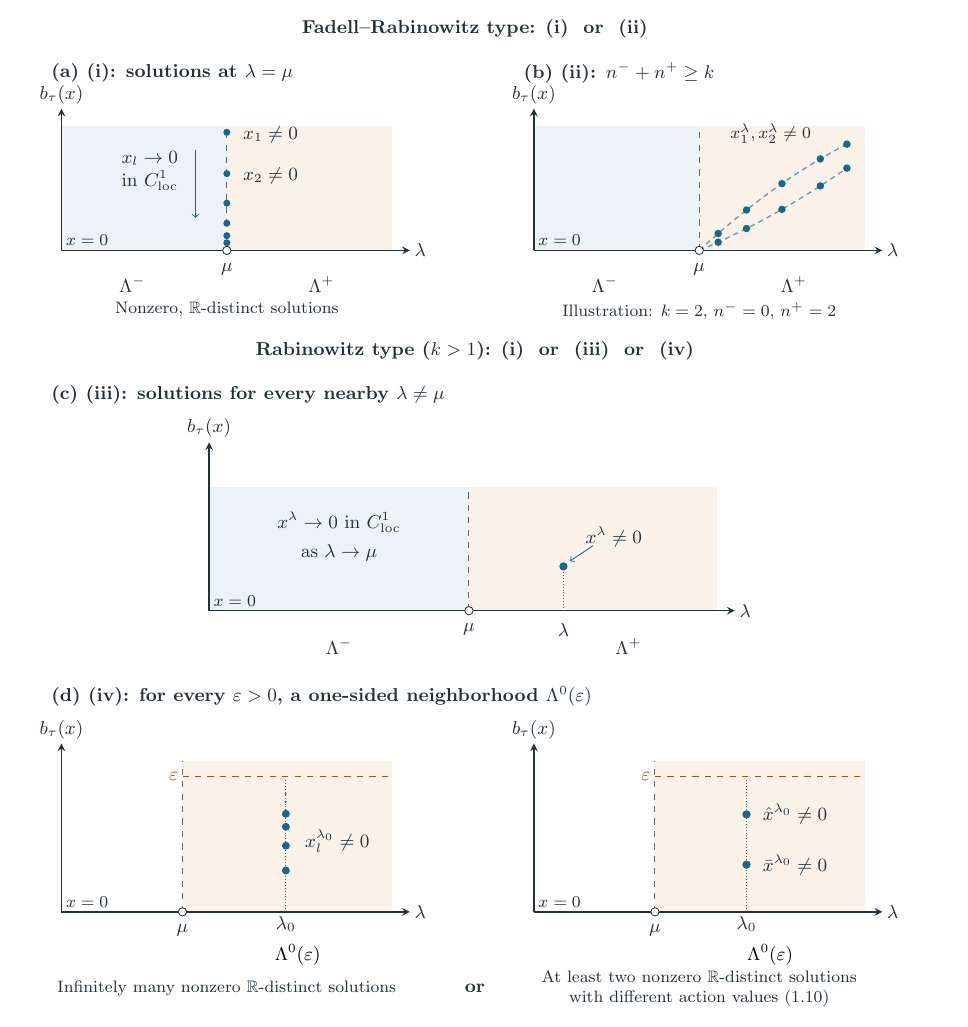}
	\caption{(Continued.) Schematic solution alternatives in Theorem \ref{thm2}(III).
		}
	\label{fig:HL-thm12-alternatives}
\end{figure}
\FloatBarrier

The following corollary is a consequence of Theorem \ref{thm2}.
\begin{cor}\label{co1}
	Let $\mathbf{F} :\mathbb{R}^{n}\to \mathbb{R} $ be an even $C^{2}$ function, and let $\mathbf{\Theta }  \in \mathcal{L} _{s}(\mathbb{R}^{n} )\cap \mathrm{GL}(n,\mathbb{R} )$ be a symmetric nonsingular matrix. For a nonzero real $\lambda$, consider the distributed delay differential system 
	\begin{equation}\label{e3}
		\begin{cases}
			\dot{x}(t)=-\nabla\mathbf{F} \left (  \int_{0}^{\lambda } \mathbf{\Theta }  x\left ( t-s \right )   ds \right )- \cdots-\nabla \mathbf{F} \left (\int_{\left ( m-1 \right ) \lambda }^{m\lambda } \mathbf{\Theta } x\left ( t-s \right )   ds \right ),\\
			x\left ( t+m\lambda  \right ) =-x\left ( t \right ),\quad  t\in \mathbb{R}, \ x(t)\in \mathbb{R}^{n}.
		\end{cases}
	\end{equation} 
	Assume that matrices $\nabla ^{2}\mathbf{F} (0)$ and $\mathbf{\Theta}  $ commute, and thus there exists an orthogonal matrix $T$ of order $n$  such that 
	\begin{equation}\label{T2}
		T^{-1}\nabla ^{2}\mathbf{F}(0)T={\rm diag}\left ( f_{1},\cdots ,f_{n} \right )\ {\rm and}\ T^{-1}\mathbf{\Theta} T={\rm diag}\left ( \theta _{1},\cdots ,\theta _{n} \right ).
	\end{equation}
	If $\left ( \mu,0 \right ) $ with $\mu \in \mathbb{R}  \setminus \left \{ 0 \right \} $ is a bifurcation point of system \eqref{e3}, then $m\left | \mu  \right | \sqrt{2f_{j}\theta _{j}} \in \left ( 2\mathbb{N}-1  \right )\pi $ for some $j$. Conversely, suppose $\mu \ne 0$ and there exists a partition $\left \{ j_{1},\cdots ,j_{k} \right \} \cup \left \{ j_{k+1},\cdots ,j_{n} \right \} $ of $\left \{ 1,\cdots ,n \right \} $ with $1\le k\le n$ such that the following two conditions hold,
	\begin{enumerate}
		\item[{\rm (D)}] For each $l=k+1,\cdots,n$, either $2\left | \mu  \right | ^{2}f_{j_{l}}\theta _{j_{l}}< \left ( \pi /m \right ) ^{2}$ or 
		\begin{equation}\label{2r12}
			\left ( 2\rho _{l}(\mu)-1 \right ) \pi <m\left | \mu  \right | \sqrt{2f_{j_{l}}\theta _{j_{l}}} < \left ( 2\rho _{l}(\mu)+1 \right ) \pi ,
		\end{equation}
		for some $\rho _{l}(\mu)\in\mathbb{N} $.
		\item[{\rm (E)}] For each $l=1,\cdots,k$, $f_{j_{l}}\theta _{j_{l}}> 0$,  and $m\left | \mu  \right | \sqrt{2f_{j_{l}}\theta _{j_{l}}} =\left ( 2\rho _{l}(\mu)-1 \right ) \pi  $ for some $\rho _{l}(\mu)\in\mathbb{N} $. Moreover, $f _{j_{1}},\cdots,f _{j_{k}}$ have the same sign.
	\end{enumerate}
	Then, for system \eqref{e3}, there exist two possible alternatives,
	\begin{enumerate}[{\rm (i)}]
		\item The system \eqref{e3} with $\lambda=\mu$ admits a sequence of nonzero $\mathbb{R} $-distinct solutions $\left \{ x_{l} \right \} _{l}^{\infty }$ such that $x_{l}\to 0$ in $C_{{\rm loc}}^{1}\left ( \mathbb{R},\mathbb{R}^{n}  \right ) $ as $l\to \infty $.
		\item There exist left and right neighborhoods $\Lambda^{-}$ and $\Lambda^{+}$ of $\mu$ in $\mathbb{R} $, respectively, containing no $0$,   and nonnegative integers $n^{+}$ and $n^{-}$ satisfying $n^{+}+n^{-}\ge k$, such that for $\lambda \in \Lambda ^{-}\setminus \left \{ \mu  \right \} (resp. \ \lambda \in \Lambda ^{+}\setminus \left \{ \mu  \right \})$, system \eqref{e3} with parameter $\lambda$ has at least $n^{-}(resp.\ n^{+})$ $\mathbb{R} $-distinct nonzero solutions $x_{i}^{\lambda}(i=1,\cdots,n^{-}(resp.\ n^{+}))$ which converge to zero in $C_{{\rm loc}}^{1}\left ( \mathbb{R},\mathbb{R}^{n}  \right ) $ as $\lambda\to \mu$.
		\item[]
		\hspace*{-\leftmargin}%
		\parbox{\dimexpr\linewidth+\leftmargin\relax}{%
			Moreover, if $k>1$, then at least one of {\rm (i)} and the following {\rm (iii)} and {\rm (iv)} holds,}
		\item For every nonzero $\lambda \in \Lambda \setminus \left \{ \mu  \right \} $ near $\mu$, there is a nonzero solution $x^{\lambda}$ of system \eqref{e3} with parameter $\lambda$, such that   $x^{\lambda}\to 0$ in $C_{{\rm loc}}^{1}\left ( \mathbb{R},\mathbb{R}^{n}  \right ) $ as $\lambda\to \mu$.
		\item  For a given $\varepsilon >0$, there exists a one-sided neighborhood $\Lambda^{0}$ of $\mu$ in $\Lambda$ such that, for any nonzero $\lambda\in \Lambda ^{0}\setminus \left \{ \mu \right \} $, system \eqref{e3} with parameter $\lambda$ has either infinitely many nonzero $\mathbb{R} $-distinct solutions $\left \{ x_{l }^{\lambda} \right \}_{l=1}^{\infty } $ such that $\left \| x_{l}^{\lambda} \mid _{\left [ 0,\left | \lambda  \right |  \right ] }\right \| _{C^{1}}<\varepsilon $ for all $l\in \mathbb{N} $, or at least two nonzero $\mathbb{R} $-distinct solutions $\bar{x} ^{\lambda } $ and $\hat{x} ^{\lambda } $ such that $\left \| \bar{x}^{\lambda } \mid _{\left [ 0,\left | \lambda  \right |  \right ] }\right \| _{C^{1}}<\varepsilon $, $\left \| \hat{x}^{\lambda } \mid _{\left [ 0,\left | \lambda  \right |  \right ] }\right \| _{C^{1}}<\varepsilon $, and 
		\begin{equation}\label{c1}
			\begin{split}
				&-\frac{1}{4} \int_{0}^{m\lambda } \left ( \mathbf{\Theta }  \left ( \bar{x}^{\lambda } (t)-\bar{x}^{\lambda }  (t-\lambda  ) \right ) ,\bar{x}^{\lambda }  (t)-\bar{x}^{\lambda }  (t-\lambda )  \right )dt+\int_{0}^{m\lambda  }\mathbf{F} \left (\int_{0}^{\lambda }\mathbf{\Theta }  \bar{x}^{\lambda }  (t-s)ds  \right ) dt   \\
				\ne& -\frac{1}{4} \int_{0}^{m\lambda } \left ( \mathbf{\Theta }  \left ( \hat{x}^{\lambda } (t)-\hat{x}^{\lambda }  (t-\lambda ) \right ) ,\hat{x}^{\lambda }  (t)-\hat{x}^{\lambda }  (t-\lambda  )  \right )dt+\int_{0}^{m\lambda }\mathbf{F} \left (\int_{0}^{\lambda }\mathbf{\Theta } \hat{x}^{\lambda }  (t-s)ds  \right ) dt.   
			\end{split}
		\end{equation}						
	\end{enumerate}
\end{cor}
 \begin{thm}\label{th4}
 	Under assumption $({\bf FG})$, with $n=1$,  delay differential system \eqref{e3'} with two distributed delays becomes
 	\begin{equation}\label{e4'}
 		\begin{cases}
 			\dot{x}(t)=-F' \left ( \lambda , \int_{0}^{\tau }   G'(\lambda, x\left ( t-s \right ))   ds \right )-F' \left ( \lambda ,  \int_{ \tau }^{2\tau }G'(\lambda, x\left ( t-s \right ) )  ds \right ),\\
 			x\left ( t+2\tau  \right ) =-x\left ( t \right ),\quad  t\in \mathbb{R}, \ x(t)\in \mathbb{R}.
 		\end{cases}
 	\end{equation}
  Here, $F'(\lambda ,\cdot )$, $F''(\lambda ,\cdot )$ and $G'(\lambda ,\cdot )$, $G''(\lambda ,\cdot )$ denote, respectively,  the first and second derivatives of $F(\lambda ,\cdot )$ and $G(\lambda ,\cdot )$ with respect to the second variable. Let $f(\lambda )=F'' (\lambda ,0)$ and $g(\lambda )=G'' (\lambda ,0)$, then 
 	\begin{enumerate}[{\rm (I)}]
 		\item {\rm (Necessary condition):} If $\left ( \mu ,0 \right )$ with $\mu \in \Lambda $ is a bifurcation point of system \eqref{e4'}, then $2\sqrt{2g(\mu)f(\mu  )}\tau \in  \left ( 2\mathbb{N}-1  \right ) \pi $.
 		\item {\rm (Sufficient condition):} Let $\mu$ be an interior point of $\Lambda$ satisfying the following conditions: 
 		\begin{enumerate}
 			\item[{\rm (c)}] $g(\mu)f(\mu  )>0$ and $2\sqrt{2g(\mu)f(\mu  )}\tau \in  \left ( 2\mathbb{N}-1  \right ) \pi $.
 			\item[{\rm (d)}] There exist sequences $\left ( \lambda _{q}^{\pm } \right ) _{q}\subset \Lambda $ converging to $\mu$ such that $g(\lambda _{q}^{-} )f(\lambda _{q}^{-} )< g(\mu  )f(\mu  )< g(\lambda _{q}^{+} )f(\lambda _{q}^{+} )$ for all $q\in \mathbb{N} $.
 		\end{enumerate}
 		Then, $\left ( \mu ,0 \right )$ is a bifurcation point of system \eqref{e4'}.
 		\item {\rm (Alternative bifurcations of Fadell-Rabinowitz type):} Let the assumptions of {\rm (II)} hold, with condition {\rm (d)} replaced by the stronger condition:
 		\begin{enumerate}[{\rm (sd)}]
 			\item For some $\varepsilon >0$, one of the following two conditions holds:
 			\begin{enumerate}[{\rm (sd.1)}]
 				\item $g(\lambda )f(\lambda )< g(\mu  )f(\mu  )< g(\lambda'  )f(\lambda'  )$ for all $\left ( \lambda ,\lambda ' \right )\in \left ( \mu-\varepsilon ,\mu  \right )\times \left ( \mu ,\mu +\varepsilon \right ) $.
 				\item $g(\lambda )f(\lambda )> g(\mu  )f(\mu  )> g(\lambda'  )f(\lambda'  )$ for all $\left ( \lambda ,\lambda ' \right )\in \left ( \mu-\varepsilon ,\mu  \right )\times \left ( \mu ,\mu +\varepsilon \right ) $.
 			\end{enumerate}
 		\end{enumerate}
 		Then, for system \eqref{e4'}, there exist the following alternatives:
 		\begin{enumerate}[{\rm (i)}]
 			\item The system \eqref{e4'} with $\lambda=\mu$ admits a sequence of nonzero $\mathbb{R} $-distinct solutions $\left \{ x_{l} \right \} _{l}^{\infty }$ such that $x_{l}\to 0$ in $C_{{\rm loc}}^{1}\left ( \mathbb{R},\mathbb{R}  \right ) $ as $l\to \infty $.
 			\item There exist left and right neighborhoods $\Lambda^{-}$ and $\Lambda^{+}$ of $\mu$ in $\Lambda$ and nonnegative integers $n^{+}$ and $n^{-}$ satisfying $n^{+}+n^{-}\ge 1$, such that for $\lambda \in \Lambda ^{-}\setminus \left \{ \mu  \right \} (resp. \ \lambda \in \Lambda ^{+}\setminus \left \{ \mu  \right \})$, system \eqref{e4'} with parameter $\lambda$ has at least $n^{-}(resp.\ n^{+})$ $\mathbb{R} $-distinct nonzero solutions $x^{\lambda}_{i}(i=1,\cdots,n^{-}(resp.\ n^{+}))$ which converge to zero in $C_{{\rm loc}}^{1}\left ( \mathbb{R},\mathbb{R}  \right ) $ as $\lambda\to \mu$.
 		\end{enumerate}
 	\end{enumerate}
 \end{thm}
 The crossing conditions and bifurcation alternatives in
 Theorem \ref{th4}(III) are illustrated schematically by
 Figure~\ref{fig:HL-thm11}, with
 $p(\lambda)=g(\lambda)f(\lambda)$ and
 $p^*=p(\mu)=(2j-1)^2\pi^2/(8\tau^2)$.
 \begin{thm}\label{thm5}
 	Suppose that $({\bf FG})$ holds with $n>1$, then we have 
 	\begin{enumerate}[{\rm (I)}]
 		\item {\rm (Necessary condition):} If $\left ( \mu ,0 \right )$ with $\mu \in \Lambda $ is a bifurcation point of system  \eqref{e3'}, and if the $n\times n$ real symmetric matrices $\nabla^{2}F(\mu ,0)$ and $\nabla^{2}G(\mu ,0)$ commute, that is, there exists a real orthogonal $n\times n$ matrix $T$ such that 
 		\begin{equation}
 			T^{-1}\nabla ^{2}G(\mu ,0)T=\mathrm{diag}(g_{1}(\mu ),\cdots ,g_{n}(\mu ))\ \mathrm{and} \ T^{-1}\nabla ^{2}F(\mu ,0)T=\mathrm{diag}(f_{1}(\mu ),\cdots ,f_{n}(\mu )),
 		\end{equation}
 		then $2\sqrt{2g_{j}(\mu  )f_{j}(\mu  )}\tau \in  \left ( 2\mathbb{N}-1  \right ) \pi $ for some $j\in \left \{ 1,\cdots ,n \right \} $.
 		\item {\rm (Sufficient condition):} For each $\lambda\in \Lambda$, suppose $\nabla ^{2}F(\lambda ,0)$ and $\nabla ^{2}G(\lambda ,0)$ commute. Let $\mu$ be an interior point of $\Lambda$, assume further that the matrix $\nabla ^{2}G(\mu ,0)$ has $n$ distinct eigenvalues $g_{j}(\mu)$, $j=1,\cdots,n$. 
 		Suppose there exists a partition $\left \{ j_{1},\cdots ,j_{k} \right \} \cup \left \{ j_{k+1},\cdots ,j_{n} \right \} $ of $\left \{ 1,\cdots ,n \right \} $ with $1\le k\le n$, such that the following conditions hold,
 		\begin{enumerate}
 			\item[{\rm (U)}] For each $l=k+1,\cdots,n$, either $2g_{j_{l}}(\mu)f_{j_{l}}(\mu)< \left ( \pi /2\tau  \right ) ^{2}$ or 
 			\begin{equation}\label{2r'}
 				\left ( 2\rho _{l}(\mu)-1 \right ) \pi <2\sqrt{2g_{j_{l}}(\mu)f_{j_{l}}(\mu)} \tau < \left ( 2\rho _{l}(\mu)+1 \right ) \pi ,
 			\end{equation}
 			for some $\rho _{l}(\mu)\in\mathbb{N} $.
 			\item[{\rm (V)}] For each $l=1,\cdots,k$, $g_{j_{l}}(\mu)f_{j_{l}}(\mu)> 0$ and $2\sqrt{2g_{j_{l}}(\mu)f_{j_{l}}(\mu)}\tau =\left ( 2\rho _{l}(\mu)-1 \right ) \pi  $ for some $\rho _{l}(\mu)\in\mathbb{N} $.
 			\item[{\rm (W)}] There exist sequence $\left ( \lambda _{q}^{\pm } \right ) _{q}\subset \Lambda $  converging to $\mu$ such that $g_{j_{l}}(\lambda_{q}^{-} )f_{j_{l}}(\lambda_{q}^{-} )<g_{j_{l}}(\mu)f_{j_{l}}(\mu)<g_{j_{l}}(\lambda_{q}^{+} )f_{j_{l}}(\lambda_{q}^{+} )$ if $ g_{j_{l}}(\mu  )>0$, and $g_{j_{l}}(\lambda_{q}^{-} )f_{j_{l}}(\lambda_{q}^{-} )>g_{j_{l}}(\mu)f_{j_{l}}(\mu)>g_{j_{l}}(\lambda_{q}^{+} )f_{j_{l}}(\lambda_{q}^{+} )$ if $ g_{j_{l}}(\mu  )<0$, for all $q\in \mathbb{N} $ and  $l=1,\cdots,k$.
 		\end{enumerate}
 		then $\left ( \mu ,0 \right )$ is a bifurcation point of system  \eqref{e3'}.
 		\item {\rm (Alternative bifurcations of Fadell-Rabinowitz type and of Rabinowitz type):} Let the assumptions of {\rm (II)} hold, with {\rm (W)} strengthened as follows,
 		\begin{enumerate}[{\rm (SW)}]
 			\item	For some $\varepsilon >0$, one of the following two conditions holds,
 			\begin{enumerate}[{\rm (SW.1)}]
 				\item $g_{j_{l}}(\lambda )f_{j_{l}}(\lambda )< g_{j_{l}}(\mu  )f_{j_{l}}(\mu  )< g_{j_{l}}(\lambda'  )f_{j_{l}}(\lambda'  )$ for all $\left ( \lambda ,\lambda ' \right )\in \left ( \mu-\varepsilon ,\mu  \right )\times \left ( \mu ,\mu +\varepsilon \right ) $ if $ g_{j_{l}}(\mu  )>0$, and $g_{j_{l}}(\lambda )f_{j_{l}}(\lambda )> g_{j_{l}}(\mu  )f_{j_{l}}(\mu  )> g_{j_{l}}(\lambda'  )f_{j_{l}}(\lambda'  )$ for all $\left ( \lambda ,\lambda ' \right )\in \left ( \mu-\varepsilon ,\mu  \right )\times \left ( \mu ,\mu +\varepsilon \right ) $ if $ g_{j_{l}}(\mu  )<0$, $l=1,\cdots,k$.
 				\item $g_{j_{l}}(\lambda )f_{j_{l}}(\lambda )> g_{j_{l}}(\mu  )f_{j_{l}}(\mu  )> g_{j_{l}}(\lambda'  )f_{j_{l}}(\lambda'  )$ for all $\left ( \lambda ,\lambda ' \right )\in \left ( \mu-\varepsilon ,\mu  \right )\times \left ( \mu ,\mu +\varepsilon \right ) $ if $ g_{j_{l}}(\mu  )>0$, and $g_{j_{l}}(\lambda )f_{j_{l}}(\lambda )< g_{j_{l}}(\mu  )f_{j_{l}}(\mu  )< g_{j_{l}}(\lambda'  )f_{j_{l}}(\lambda'  )$ for all $\left ( \lambda ,\lambda ' \right )\in \left ( \mu-\varepsilon ,\mu  \right )\times \left ( \mu ,\mu +\varepsilon \right ) $ if $ g_{j_{l}}(\mu  )<0$, $l=1,\cdots,k$.
 			\end{enumerate}
 		\end{enumerate}
 		Then for system \eqref{e3'}, there exist two possible alternatives:
 		\begin{enumerate}[{\rm (i)}]
 			\item The system \eqref{e3'} with $\lambda=\mu$ admits a sequence of nonzero $\mathbb{R} $-distinct solutions $\left \{ x_{l} \right \} _{l}^{\infty }$ such that $x_{l}\to 0$ in $C_{{\rm loc}}^{1}\left ( \mathbb{R},\mathbb{R}^{n}  \right ) $ as $l\to \infty $.
 			\item There exist left and right neighborhoods $\Lambda^{-}$ and $\Lambda^{+}$ of $\mu$ in $\Lambda$ and nonnegative integers $n^{+}$ and $n^{-}$ satisfying $n^{+}+n^{-}\ge k$, such that for $\lambda \in \Lambda ^{-}\setminus \left \{ \mu  \right \} (resp. \ \lambda \in \Lambda ^{+}\setminus \left \{ \mu  \right \})$, system \eqref{e3'} with parameter $\lambda$ has at least $n^{-}(resp.\ n^{+})$ $\mathbb{R} $-distinct nonzero solutions $x_{i}^{\lambda}(i=1,\cdots,n^{-}(resp.\ n^{+}))$ which converge to zero in $C_{{\rm loc}}^{1}\left ( \mathbb{R},\mathbb{R}^{n}  \right ) $ as $\lambda\to \mu$.
 			\item[]
 			\hspace{-\labelwidth}\hspace{-\labelsep}
 			Moreover, if $k>1$, then at least one of {\rm (i)} and the following {\rm (iii)} and {\rm (iv)} holds,
 			\item For every $\lambda \in \Lambda \setminus \left \{ \mu  \right \} $ near $\mu$, there is a nonzero solution $x^{\lambda}$ of system \eqref{e3'} with parameter $\lambda$, such that   $x^{\lambda}\to 0$ in $C_{{\rm loc}}^{1}\left ( \mathbb{R},\mathbb{R}^{n}  \right ) $ as $\lambda\to \mu$.
 			\item  For a given $\varepsilon >0$, there exists a one-sided neighborhood $\Lambda^{0}$ of $\mu$ in $\Lambda$ such that, for any $\lambda\in \Lambda ^{0}\setminus \left \{ \mu \right \} $, system \eqref{e3'} with parameter $\lambda$ has either infinitely many nonzero $\mathbb{R} $-distinct solutions $\left \{ x_{l }^{\lambda} \right \}_{l=1}^{\infty } $ such that $\left \| x_{l}^{\lambda} \mid _{\left [ 0,\tau \right ] }\right \| _{C^{1}}<\varepsilon $ for all $l\in \mathbb{N} $, or at least two nonzero $\mathbb{R} $-distinct solutions $\bar{x} ^{\lambda } $ and $\hat{x} ^{\lambda } $ such that $\left \| \bar{x}^{\lambda } \mid _{\left [ 0,\tau \right ] }\right \| _{C^{1}}<\varepsilon $, $\left \| \hat{x}^{\lambda } \mid _{\left [ 0,\tau \right ] }\right \| _{C^{1}}<\varepsilon $, and 
\begin{equation}\label{124}
	\begin{split}
		&\int_{0}^{2\tau }G(\lambda , \bar{x} ^{\lambda }(t))dt- \int_{0}^{2\tau } \left ( \nabla G\left (\lambda, \bar{x}^{\lambda } \left ( t\right ) \right ),\bar{x}^{\lambda }(t) \right )dt\\
		&+\int_{0}^{2\tau }F\left ( \lambda , \int_{0}^{\tau }\nabla G\left (\lambda, \bar{x}^{\lambda } \left ( t-s \right ) \right ) ds \right ) dt\\
		\ne &\int_{0}^{2\tau }G(\lambda , \hat{x} ^{\lambda }(t))dt- \int_{0}^{2\tau } \left ( \nabla G\left (\lambda, \hat{x}^{\lambda } \left ( t\right ) \right ),\hat{x}^{\lambda }(t) \right )dt\\
		&+\int_{0}^{2\tau }F\left ( \lambda , \int_{0}^{\tau }\nabla G\left (\lambda, \hat{x}^{\lambda } \left ( t-s \right ) \right ) ds \right ) dt.
	\end{split}
\end{equation}	
 		\end{enumerate}
 	\end{enumerate}
 \end{thm}
\begin{rem}
	Since the continuous, real symmetric $n\times n$ matrix-valued functions $\nabla ^{2}F(\lambda ,0)$ and $\nabla ^{2}G(\lambda ,0)$ commute for all $\lambda\in \Lambda	$,  Lemma 6.2 in \cite{Lu} assures that in a neighborhood $\Lambda_{0}$ of $\mu$ in $\Lambda$ there exist continuous real functions $g_{j}, f_{j}$ for $j=1,\cdots,n$, and a continuous real orthogonal $n\times n$ matrix-valued function $T$ such that for each $\lambda\in \Lambda_{0}$,
	\begin{equation}
		T(\lambda )^{-1}\nabla ^{2}G(\lambda ,0)T(\lambda )=\mathrm{diag}(g_{1}(\lambda ),\cdots ,g_{n}(\lambda ))\ \mathrm{and} \ 	T(\lambda )^{-1}\nabla ^{2}F(\lambda ,0)T(\lambda )=\mathrm{diag}(f_{1}(\lambda ),\cdots ,f_{n}(\lambda )).
	\end{equation}
\end{rem}
The bifurcation alternatives in Theorem \ref{thm5}(III) are illustrated
schematically by the solution diagrams in
Figure~\ref{fig:HL-thm12-alternatives}, with \eqref{124} replacing \eqref{c}.
For the crossing diagrams in Figure~\ref{fig:HL-thm12-crossings},
the relevant quantities are
\[
p_l(\lambda)=g_{j_l}(\lambda)f_{j_l}(\lambda),\qquad
p_l^*=p_l(\mu)
=\frac{(2\rho_l(\mu)-1)^2\pi^2}{8\tau^2}.
\]
 \begin{cor}\label{co2}
 	Let $\mathbf{F}, \mathbf{G}:\mathbb{R}^{n}\to \mathbb{R} $ be two even $C^{2}$ functions. For a nonzero real $\lambda$, consider the distributed delay differential system 
 	\begin{equation}\label{e31}
 		\begin{cases}
 			\dot{x}(t)=-\nabla\mathbf{F} \left (  \int_{0}^{\lambda } \nabla\mathbf{G}(x\left ( t-s \right ))    ds \right )-\nabla \mathbf{F} \left (\int_{ \lambda }^{2\lambda }\nabla\mathbf{G}(x\left ( t-s \right ) )   ds \right ),\\
 			x\left ( t+2\lambda  \right ) =-x\left ( t \right ),\quad  t\in \mathbb{R}, \ x(t)\in \mathbb{R}^{n}.
 		\end{cases}
 	\end{equation} 
 	Assume that matrices $\nabla ^{2}\mathbf{F} (0)$ and $\nabla ^{2}\mathbf{G} (0)$ commute, and thus there exists an orthogonal matrix $T$ of order $n$  such that 
 	\begin{equation}\label{T21}
 		T^{-1}\nabla ^{2}\mathbf{F}(0)T={\rm diag}\left ( f_{1},\cdots ,f_{n} \right )\ {\rm and}\ T^{-1}\nabla ^{2}\mathbf{G}(0) T={\rm diag}\left ( g _{1},\cdots ,g _{n} \right ).
 	\end{equation}
 	If $\left ( \mu,0 \right ) $ with $\mu \in \mathbb{R}  \setminus \left \{ 0 \right \} $ is a bifurcation point of system \eqref{e31}, then $2\left | \mu  \right | \sqrt{2g_{j}f_{j}} \in \left ( 2\mathbb{N}-1  \right )\pi $ for some $j$. Conversely, suppose $\mu \ne 0$ and there exists a partition $\left \{ j_{1},\cdots ,j_{k} \right \} \cup \left \{ j_{k+1},\cdots ,j_{n} \right \} $ of $\left \{ 1,\cdots ,n \right \} $ with $1\le k\le n$ such that the following two conditions hold,
 	\begin{enumerate}
 		\item[{\rm (X)}] For each $l=k+1,\cdots,n$, either $\left | \mu  \right | ^{2}g_{j_{l}}f_{j_{l}}<  \pi^{2} /8 $ or 
 		\begin{equation}\label{2r1}
 			\left ( 2\rho _{l}(\mu)-1 \right ) \pi <2\left | \mu  \right | \sqrt{2g_{j_{l}}f_{j_{l}}} < \left ( 2\rho _{l}(\mu)+1 \right ) \pi ,
 		\end{equation}
 		for some $\rho _{l}(\mu)\in\mathbb{N} $.
 		\item[{\rm (Y)}] For each $l=1,\cdots,k$, $g_{j_{l}}f_{j_{l}}> 0$,  and $2\left | \mu  \right | \sqrt{2g_{j_{l}}f_{j_{l}}} =\left ( 2\rho _{l}(\mu)-1 \right ) \pi  $ for some $\rho _{l}(\mu)\in\mathbb{N} $. Moreover, $g _{j_{1}},\cdots,g_{j_{k}}$ have the same sign.
 	\end{enumerate}
 	Then, for system \eqref{e31}, there exist two possible alternatives,
 	\begin{enumerate}[{\rm (i)}]
 		\item The system \eqref{e31} with $\lambda=\mu$ admits a sequence of nonzero $\mathbb{R} $-distinct solutions $\left \{ x_{l} \right \} _{l}^{\infty }$ such that $x_{l}\to 0$ in $C_{{\rm loc}}^{1}\left ( \mathbb{R},\mathbb{R}^{n}  \right ) $ as $l\to \infty $.
 		\item There exist left and right neighborhoods $\Lambda^{-}$ and $\Lambda^{+}$ of $\mu$ in $\mathbb{R} $, respectively, containing no $0$,   and nonnegative integers $n^{+}$ and $n^{-}$ satisfying $n^{+}+n^{-}\ge k$, such that for $\lambda \in \Lambda ^{-}\setminus \left \{ \mu  \right \} (resp. \ \lambda \in \Lambda ^{+}\setminus \left \{ \mu  \right \})$, system \eqref{e31} with parameter $\lambda$ has at least $n^{-}(resp.\ n^{+})$ $\mathbb{R} $-distinct nonzero solutions $x_{i}^{\lambda}(i=1,\cdots,n^{-}(resp.\ n^{+}))$ which converge to zero in $C_{{\rm loc}}^{1}\left ( \mathbb{R},\mathbb{R}^{n}  \right ) $ as $\lambda\to \mu$.
 		\item[]
 		\hspace*{-\leftmargin}%
 		\parbox{\dimexpr\linewidth+\leftmargin\relax}{%
 			Moreover, if $k>1$, then at least one of {\rm (i)} and the following {\rm (iii)} and {\rm (iv)} holds,}
 		\item For every nonzero $\lambda \in \Lambda \setminus \left \{ \mu  \right \} $ near $\mu$, there is a nonzero solution $x^{\lambda}$ of system \eqref{e31} with parameter $\lambda$, such that   $x^{\lambda}\to 0$ in $C_{{\rm loc}}^{1}\left ( \mathbb{R},\mathbb{R}^{n}  \right ) $ as $\lambda\to \mu$.
 		\item  For a given $\varepsilon >0$, there exists a one-sided neighborhood $\Lambda^{0}$ of $\mu$ in $\Lambda$ such that, for any nonzero $\lambda\in \Lambda ^{0}\setminus \left \{ \mu \right \} $, system \eqref{e31} with parameter $\lambda$ has either infinitely many nonzero $\mathbb{R} $-distinct solutions $\left \{ x_{l }^{\lambda} \right \}_{l=1}^{\infty } $ such that $\left \| x_{l}^{\lambda} \mid _{\left [ 0,\left | \lambda  \right |  \right ] }\right \| _{C^{1}}<\varepsilon $ for all $l\in \mathbb{N} $, or at least two nonzero $\mathbb{R} $-distinct solutions $\bar{x} ^{\lambda } $ and $\hat{x} ^{\lambda } $ such that $\left \| \bar{x}^{\lambda } \mid _{\left [ 0,\left | \lambda  \right |  \right ] }\right \| _{C^{1}}<\varepsilon $, $\left \| \hat{x}^{\lambda } \mid _{\left [ 0,\left | \lambda  \right |  \right ] }\right \| _{C^{1}}<\varepsilon $, and 
 		\begin{equation}\label{c12}
 		\begin{split}
 			&\int_{0}^{2\lambda  }\mathbf{G} ( \bar{x} ^{\lambda }(t))dt- \int_{0}^{2\lambda  } \left ( \nabla \mathbf{G} \left (\bar{x}^{\lambda } \left ( t\right ) \right ),\bar{x}^{\lambda }(t) \right )dt+\int_{0}^{2\lambda  }\mathbf{F} \left (  \int_{0}^{\lambda  }\nabla \mathbf{G} \left (\bar{x}^{\lambda } \left ( t-s \right ) \right ) ds \right ) dt\\
 			\ne &\int_{0}^{2\lambda  }\mathbf{G} (\hat{x} ^{\lambda }(t))dt- \int_{0}^{2\lambda } \left ( \nabla \mathbf{G} \left (\hat{x}^{\lambda } \left ( t\right ) \right ),\hat{x}^{\lambda }(t) \right )dt+\int_{0}^{2\lambda  }\mathbf{F} \left ( \int_{0}^{\lambda }\nabla \mathbf{G} \left ( \hat{x}^{\lambda } \left ( t-s \right ) \right ) ds \right ) dt.
 		\end{split}
 		\end{equation}						
 	\end{enumerate}
 \end{cor}
 The main contributions of this paper are twofold. We derive explicit formulas for the Maslov $P_{n,m}$-index and nullity (Lemmas \ref{abv} and \ref{abv2}),  which yield resonance conditions and bifurcation criteria for multiple distributed delay systems with quadratic $G$. We further exploit the anti-periodic symmetry of the two-delay case to extend these results to general even nonlinear potentials $F$ and $G$.
 
 The remainder of the paper is organized as follows. In Section \ref{sec2}, we recall the Maslov $P$-index theory and establish equivalent Hamiltonian formulations for the two classes of distributed delay systems considered in this paper. In Section \ref{sec3}, we first derive explicit formulas for the Maslov $P_{n,m}$-index and nullity in Lemmas \ref{abv} and \ref{abv2}, and then use these formulas to prove the main bifurcation results, including Theorems \ref{n1}, \ref{thm2}, \ref{th4}, and \ref{thm5} and Corollaries \ref{co1} and \ref{co2}.
\section{Variational Setting }\label{sec2}

\subsection{$P$-index}
In \cite{Liu4}, Liu extended the Maslov-type index theory for symplectic paths developed by Long in \cite{Long3}-\cite{Long} to the so-called Maslov $P$-index theory, which is adapted to the problem of a nonlinear Hamiltonian system with a $P$-boundary condition as follows,
\begin{equation}
	\left\{\begin{matrix}\label{Ps}
		\dot{x}(t)=JH'(t,x(t)), \\
		x(\tau)=Px(0),\quad \quad \ 
	\end{matrix}\right.
\end{equation}
where $J=\begin{pmatrix}
	0 &-I_{n} \\
	I_{n} &0
\end{pmatrix}$ is the standard symplectic matrix, the Hamiltonian function $H\in C^{2}\left ( \mathbb{R}\times \mathbb{R}^{2n},\mathbb{R}   \right ) $ satisfying $H(t+\tau, Px)=H(t,x)$ with $P\in {\rm Sp}(2n)$. Here,  ${\rm Sp}(2n)$ is the symplectic group defined by
\begin{equation}
	{\rm Sp}(2n)=\left \{ M\in {\rm GL}(2n,\mathbb{R} )\mid M^{T}JM=J \right \} ,
\end{equation}
where $M^{T}$ denotes the transpose of $M$. The set of symplectic paths is defined by 
\begin{equation}
	\mathcal{P}_{\tau } (2n)=\left \{ \gamma \in C\left ( \left [ 0,\tau \right ],{\rm Sp}(2n)  \right )\mid \gamma (0)=I_{2n}  \right \} ,
\end{equation}
where $I_{2n}$ denotes the $2n\times 2n$ identity matrix.

The linearized system along a solution $x(t)$ of system \eqref{Ps} is the following linear Hamiltonian system 
\begin{equation}\label{sl}
	\dot{z}(t)=JB(t)z(t),
\end{equation}
where $B(t)=H''(t,x(t))$ satisfying $P^{T}B(t+\tau)P=B(t)$. Let $\gamma (t)$ be the fundamental solution of system \eqref{sl} in the interval $\left [ 0,\tau  \right ] $, then $\gamma (t)$ is a symplectic path, that is, $\gamma \in 	\mathcal{P}_{\tau } (2n)$.  For any $\gamma \in 	\mathcal{P}_{\tau } (2n)$ and $P\in {\rm Sp}(2n)$, we associate a pair of integers $\left ( i_{\tau }^{P}(\gamma ), \nu_{\tau }^{P}(\gamma ) \right ) \in \mathbb{Z}\times \left \{ 0,1,\cdots ,2n \right \}  $, where $i_{\tau }^{P}(\gamma )$ is the Maslov $P$-index and $\nu_{\tau }^{P}(\gamma )$ is the $P$-nullity, defined as follows.
\begin{defi}\label{nu}
	For any $\gamma \in \mathcal{P}_{\tau } (2n)$, $P\in \mathrm{Sp}(2n)$, the $P$-nullity $\nu _{\tau }^{P}(\gamma )$ is defined by 
	\begin{equation}
		\nu _{\tau }^{P}(\gamma )=\dim \ker \left ( \gamma (\tau )- P \right ) .
	\end{equation}
	The Maslov $P$-index $i_{\tau }^{P}(\gamma )$ is defined by 
	\begin{equation}
		i_{\tau }^{P}(\gamma )=i_{\tau }\left ( P^{-1}\gamma \ast \xi  \right )-i_{\tau }(\xi ) ,
	\end{equation}
	where $\xi \in \mathcal{P}_{\tau } (2n)$ such that $\xi(\tau)=P^{-1}\gamma(0)=P^{-1}$ and $i_{\tau }$ is the Maslov-type index for periodic boundary condition established by Long in \cite{Long}.
\end{defi}
Via different methods, Dong in \cite{Dong1} also extended the Maslov-type index $i_{\tau}(\gamma)$ and the nullity $\nu_{\tau}(\gamma)$ for systems with periodic boundary condition to the Maslov $P$-index, denoted by $i_{\tau,P}(\gamma)$, and the $P$-nullity, denoted by $\nu_{\tau,P}(\gamma)$, for systems with $P$-boundary condition. As discussed in Proposition A.1 in \cite{Lu1}, we derive a precise relation between the index pair $\left ( i_{\tau,P}(\gamma), \nu_{\tau,P}(\gamma) \right ) $ and $\left ( i_{\tau }^{P}(\gamma ), \nu_{\tau }^{P}(\gamma ) \right ) $ as follows.
\begin{prop}\label{pin}
	For any  $\gamma \in \mathcal{P}_{\tau } (2n)$ and $P\in \mathrm{Sp}(2n)$, it holds  
	\begin{equation}
		i_{\tau,P}(\gamma)=i_{\tau }^{P}(\gamma )+k_{0} \quad {\rm and} \quad \nu_{\tau,P}(\gamma)=\nu_{\tau }^{P}(\gamma ),
	\end{equation}
	where $k_{0}\in \mathbb{Z} $ is the constant only depending on $P$.
\end{prop}
By Proposition \ref{pin}, Theorems 1.5 and 1.14 in \cite{Lu1} remain valid when formulated in terms of the index pair $\left ( i_{\tau }^{P}(\gamma ), \nu_{\tau }^{P}(\gamma ) \right ) $.

For any two paths $\gamma _{1}:\left [  0,\tau \right ] \to \mathrm{Sp}(2n)$ and $\gamma _{2}:\left [  0,\tau \right ] \to \mathrm{Sp}(2n)$ with $\gamma _{1}(\tau )=\gamma _{2}(0)$, we define their joint path by 
\begin{equation}\label{g12}
	\gamma _{2}\ast \gamma _{1}(t)=\left\{\begin{matrix}
		\gamma _{1}(2t) ,\quad \quad \ \ \ 0\le t\le \tau /2,\\
		\gamma _{2}(2t-\tau ),\quad \tau /2\le t\le \tau .
	\end{matrix}\right.
\end{equation}
Given any two $2k_{i}\times 2k_{i}$ matrices of square block form
\begin{equation}
	M_{i}=\begin{pmatrix}
		A_{i} &B_{i} \\
		C_{i}&D_{i}
	\end{pmatrix},\quad i=1,2,
\end{equation}
we define an operation $\Diamond $-sum of $M_{1}$ and $M_{2}$ to be the $2(k_{1}+k_{2})\times 2(k_{1}+k_{2})$ symplectic matrix given by 
\begin{equation}\label{Dia}
	M_{1}\Diamond M_{2}=\begin{pmatrix}
		A_{1} &0  &B_{1}  &0 \\
		0&  A_{2}&  0& B_{2}\\
		C_{1}&0  &D_{1}  &0 \\
		0&  C_{2}&  0&D_{2}
	\end{pmatrix}.
\end{equation}
We denote by $M^{\Diamond k}$ the $k$-fold $\Diamond$-sum $M\Diamond \cdots \Diamond M$. The following basic properties of the \(\Diamond\)-sum play a crucial role in the proof of the main results.
\begin{lem}\label{pDia}
	The $\Diamond $-sum  of square block matrices defined in \eqref{Dia} satisfies 
	\begin{enumerate}[{\rm (i)}]
		\item For matrices $M_{1},N_{1}\in \mathbb{R} ^{2k_{1}\times 2k_{1}}$ and $M_{2},N_{2}\in \mathbb{R} ^{2k_{2}\times 2k_{2}}$, there holds 
		\begin{equation}
			\left ( M_{1}\Diamond M_{2} \right ) \left ( N_{1}\Diamond N_{2} \right ) =\left ( M_{1}N_{1} \right ) \Diamond \left ( M_{2}N_{2} \right ).
		\end{equation}
		\item $ I_{2k_{1}}\Diamond  I_{2k_{2}}= I_{2\left ( k_{1}+k_{2} \right ) }$; the $\Diamond $-sum  $M_{1}\Diamond M_{2}$ of invertible matrices $M_{1}\in\mathbb{R} ^{2k_{1}\times 2k_{1}}$ and $M_{2}\in\mathbb{R} ^{2k_{2}\times 2k_{2}}$ is invertible, and $\left ( M_{1}\Diamond M_{2} \right )^{-1}=M_{1}^{-1}\Diamond M_{2}^{-1}$.
		\item $M_{1}\Diamond \cdots \Diamond M_{k}+N_{1}\Diamond \cdots \Diamond N_{k}=\left ( M_{1}+ N_{1} \right ) \Diamond \cdots \Diamond \left ( M_{k}+ N_{k} \right ) $.
		\item $\left ( M_{1}+\cdots +M_{k} \right ) \Diamond \left ( N_{1}+\cdots +N_{k} \right )=M_{1}\Diamond N_{1}+\cdots +M_{k}\Diamond N_{k}$ and 
		$$\exp \left ( M_{1}\Diamond \cdots \Diamond M_{k} \right ) =\exp \left ( M_{1} \right ) \Diamond \cdots \Diamond \exp \left ( M_{k} \right ).$$
	\end{enumerate}
\end{lem}
We recall below several fundamental properties of the $P$-index $i_{\tau }^{P}(\gamma )$. For detailed proofs, we refer the readers to \cite{Liu5}.
\begin{prop}(Symplectic additivity)\label{Sa}
	Suppose $P_{j}\in \mathrm{Sp}(2n_{j})$, $\gamma _{j}\in \mathcal{P}_{\tau } (2n_{j})$ for $j=0,1$. Then
	\begin{equation}
		i_{\tau }^{P}\left ( \gamma _{0}\Diamond \gamma _{1} \right ) =i_{\tau }^{P_{0}}\left ( \gamma _{0} \right ) +i_{\tau }^{P_{1}}\left ( \gamma _{1} \right ),
	\end{equation} 
	where $P=P_{0} \Diamond P_{1}$ and $ \gamma _{0}\Diamond \gamma _{1} (t)=\gamma _{0}(t)\Diamond \gamma _{1} (t)$ for all $t\in \left [ 0,\tau  \right ] $.
\end{prop}
\begin{prop}(Symplectic invariance)\label{Si}
	For any $M,P\in \mathrm{Sp}(2n)$, $\gamma\in \mathcal{P}_{\tau } (2n)$, there holds 
	\begin{equation}
		i_{\tau }^{M^{-1}PM}\left ( M^{-1}\gamma M  \right )=i_{\tau}^{P}\left ( \gamma  \right ),\quad 	\nu_{\tau }^{M^{-1}PM}\left ( M^{-1}\gamma M  \right )=\nu_{\tau}^{P}\left ( \gamma  \right ).
	\end{equation}
\end{prop}
\begin{defi}
	Given two paths $\gamma _{0},\gamma _{1}\in\mathcal{P}_{\tau } (2n)$, if there exists a map $\delta \in C\left ( \left [ 0,1 \right ],\mathcal{P}_{\tau }(2n)   \right ) $ such that $\delta (0,\cdot )=\gamma _{0}(\cdot )$, $\delta (1,\cdot )=\gamma _{1}(\cdot )$, $\delta (s,0)=I_{2n}$ and $\dim \ker \left ( \delta (s,\tau )-P \right ) $ are constant for $s\in \left [ 0,1 \right ]$, then we say they are $P$-homotopic, and write $\gamma_{0}\sim _{P}\gamma_{1}$.
\end{defi}
\begin{prop}(Homotopy invariance)\label{hi}
	Let $P\in \mathrm{Sp}(2n)$, for any two paths $\gamma _{0},\gamma _{1}\in\mathcal{P}_{\tau } (2n)$, if $\gamma_{0}\sim _{P}\gamma_{1}$ on $\left [ 0,\tau \right ]$, then 
	\begin{equation}
		i_{\tau }^{P}\left ( \gamma _{0} \right ) =i_{\tau }^{P}\left ( \gamma _{1} \right ), \quad \nu _{\tau }^{P}\left ( \gamma _{0} \right ) =\nu _{\tau }^{P}\left ( \gamma _{1} \right ).
	\end{equation}
\end{prop}
For any $2n\times 2n$ continuous symmetric matrix function $B(t)$ satisfying $P^{T}B(t+\tau)P=B(t)$, let $\gamma _{B}(t)$ be the fundamental solution of linear Hamiltonian system $\dot{z}(t)=JB(t)z(t)$.
For each $0\le t\le \tau $, define 
\begin{equation}
	\nu _{t}^{P}\left ( \gamma _{B} \mid _{\left [ 0,t \right ] }\right ) =\dim \ker \left ( \gamma _{B}(t)-P \right ) .
\end{equation}
By Proposition \ref{pin}, as discussed in Lemma 4.2 of \cite{Dong1}, we derive 
\begin{lem}\label{lpn}
	Under the above assumptions, it holds
	\begin{equation}
		i_{\tau }^{P}\left ( \gamma _{B} \right ) =\sum_{0\le t< \tau }\nu _{t}^{P}\left ( \gamma _{B} \mid _{\left [ 0,t \right ] }\right ),
	\end{equation}
	if $B$ is positive definite, and 
	\begin{equation}
		i_{\tau }^{P}\left ( \gamma _{B} \right ) =-\sum_{0<t\le  \tau }\nu _{t}^{P}\left ( \gamma _{B} \mid _{\left [ 0,t \right ] }\right ),
	\end{equation}
	if $B$ is negative definite.
\end{lem}

Liu in \cite{Liu10} proved the following new result about the Maslov index pair $\left ( i_{\tau }(\gamma ),\nu _{\tau } (\gamma )\right ) $.
\begin{thm}\label{po}
	For any 
	\begin{equation}
		\gamma =\begin{pmatrix}
			A & B\\
			-B&A
		\end{pmatrix}\in \mathcal{P}_{\tau } (2n )\cap \mathrm{O}(2n),
	\end{equation}
	suppose $\lambda _{j}(t)=\exp \left ( \sqrt{-1}\theta _{j}(t)  \right )$, $j=1,\cdots ,n$, are eigenvalues of the unitary matrix function $A(t)-\sqrt{-1}B(t) $ with $\theta _{j}\in C\left ( \left [ 0,\tau  \right ],\mathbb{R}   \right ) $. There holds 
	\begin{equation}
		i_{\tau }(\gamma )=2\sum_{j=1}^{n}E\left ( \frac{\theta _{j}\left ( \tau  \right )-\theta _{j}\left ( 0 \right )  }{2\pi }  \right )+n,
	\end{equation}
	where $E\left ( a \right )=\max \left \{ k\in \mathbb{Z} |k< a \right \}  $.
\end{thm}
\subsection{Reduction to a Hamiltonian System with a $P_{n,m}$-Boundary Condition}
We first establish an equivalent Hamiltonian formulation for system \eqref{e2}. By introducing suitable auxiliary variables associated with the distributed delay terms, the original delay differential system can be transformed into a  Hamiltonian system with a $P_{n,m}$-boundary condition. For each $\lambda \in \Lambda $, let 
\begin{equation}\label{x}
	x_{1}(t)=x(t), \cdots, x_{m}(t)=x(t-(m-1)\tau ),
\end{equation}
and
\begin{equation}\label{z}
	z_{1}\left ( t \right )= \int_{0}^{\tau }\Theta (\lambda ) x\left ( t-s \right ) ds, \cdots, z_{m}\left ( t \right )= \int_{\left ( m-1 \right )\tau  }^{m\tau }\Theta (\lambda )x\left ( t-s \right )  ds,
\end{equation}
then by a straightforward computation, it holds 
\begin{equation}\label{es1}
	\begin{pmatrix}
		\ddot{z}_{1}(t) \\
		\ddot{z}_{2}(t)\\
		\vdots \\
		\ddot{z}_{m-1}(t)\\
		\ddot{z}_{m}(t)
	\end{pmatrix}=\begin{pmatrix}
		-2\Theta(\lambda) & 0 & \cdots  &0  &0 \\
		0 & -2\Theta(\lambda)  & \cdots  &0  & 0\\
		\vdots & \vdots  & \ddots   & \vdots  &\vdots  \\
		0 &  0& \cdots  & -2\Theta(\lambda) &0 \\
		0  & 0 & \cdots  &  0&-2\Theta(\lambda)
	\end{pmatrix}_{nm\times nm}
	\begin{pmatrix}
		\nabla F(\lambda ,z_{1}(t)) \\
		\nabla F(\lambda ,z_{2}(t)) \\
		\vdots \\
		\nabla F(\lambda ,z_{m-1}(t))\\
		\nabla F(\lambda ,z_{m}(t))
	\end{pmatrix}.
\end{equation}
Let $\mathcal{A} (\lambda )$ denote the $nm\times nm$ matrix defined in \eqref{es1}, then it holds $\mathcal{A} (\lambda )^{T}=\mathcal{A} (\lambda )$ for each $\lambda \in \Lambda $. Moreover, let 
\begin{equation}\label{z1}
	z(t)=\left ( {z}_{1}(t),{z}_{2}(t),\cdots, {z}_{m-1}(t), {z}_{m}(t) \right )^{T},
\end{equation}
and 
\begin{equation}\label{H}
	H(\lambda,z)=F(\lambda,z_{1})+F(\lambda,z_{2})+\cdots +F(\lambda,z_{m-1}) +F(\lambda,z_{m}),
\end{equation}
then \eqref{es1} can be rewritten as
\begin{equation}
	\ddot{z}(t)=\mathcal{A} (\lambda ) \nabla H (\lambda,z(t)).
\end{equation}
Define $nm\times nm$ matrix
\begin{equation}\label{418}
	Q_{n,m}=\begin{pmatrix}
		0&  0& \cdots  & 0 &-I_{n} \\
		I_{n}&  0& \cdots  & 0 &0 \\
		\vdots &  \vdots & \cdots  & \vdots  &\vdots  \\
		0&  0& \cdots  & 0 & 0\\
		0&  0& \cdots  & I_{n} &0
	\end{pmatrix}_{nm\times nm}=
	\begin{pmatrix}
		0 & -I_{n}\\
		I_{n(m-1)}&0
	\end{pmatrix}_{nm\times nm},
\end{equation}
then it is easily verified that $Q_{n,m}$ satisfies 
\begin{equation}
	Q_{n,m}^{T}Q_{n,m}=I_{nm}, \quad Q_{n,m}^{2m}=I_{nm},\quad \mathcal{A} (\lambda ) Q_{n,m}=Q_{n,m}\mathcal{A} (\lambda ),\ \forall \lambda \in \Lambda .
\end{equation}
Moreover, from \eqref{x} and \eqref{z}, one can infer 
\begin{equation}
	z(t+\tau )=Q_{n,m}z(t),\quad \dot{z} (t+\tau )=Q_{n,m}\dot{z} (t).
\end{equation}

We therefore have the following lemma, for which a detailed proof can be found in Lemma 2.1 in \cite{He}.
\begin{lem}\label{l1}
	Let $z$ and $H $ be defined in \eqref{z1} and \eqref{H}, if $x(t)$ is a symmetric $2m\tau$-periodic solution of distributed delay differential system (\ref{e2}), then $z(t)$ is a $2m\tau $-periodic  solution of the following second order  Hamiltonian  system with $Q_{n,m}$-boundary condition,
	\begin{equation}\label{es}
		\begin{cases}
			\ddot{z}(t)=\mathcal{A}(\lambda ) \nabla  H (\lambda ,z(t)), \\
			\dot{z} (t+\tau )=Q_{n,m}\dot{z} (t),\\
			z(t+\tau )=Q_{n,m}z(t).
		\end{cases}
	\end{equation}
	Conversely, if $z(t)=\left ( {z}_{1}(t),\cdots,  {z}_{m}(t) \right )^{T}$, where $z_{i}(t)\in \mathbb{R}^{n} $, $i=1,\cdots,m$, is a $2m\tau $-periodic  solution of  second order  Hamiltonian  system \eqref{es}, then $x(t):=\frac{1}{2}\Theta (\lambda )^{-1} \left ( \dot{z}_{1}(t)+\cdots +\dot{z}_{m}(t) \right ) $ is a symmetric $2m\tau $-periodic solution of distributed delay differential system (\ref{e2}).
\end{lem}
It is well known that system \eqref{es} is in one-to-one correspondence with a classical first order Hamiltonian system.
In fact, for each $\lambda \in \Lambda $, denote
\begin{equation} \label{eP27}
	w:=\mathcal{A}(\lambda )^{-1}\dot{z},\ u:=\begin{pmatrix}
		z\\
		w
	\end{pmatrix},\
	J_{nm}=\begin{pmatrix}
		0 &-I_{nm} \\
		I_{nm} &0
	\end{pmatrix},
\end{equation}
then the system \eqref{es} can be transformed to the following first order Hamiltonian system
\begin{equation} \label{ef}
	\begin{cases}
		\dot{u}(t)=J_{nm}\nabla \tilde{H} (\lambda ,u(t)),\\
		u(t+\tau )=P_{n,m}u(t),
	\end{cases}
\end{equation}
with
\begin{equation} \label{Hp}
	\tilde{H}(\lambda ,u) :=H(\lambda ,z)-\frac{1}{2}\left ( \mathcal{A}(\lambda )w,w \right )  ,\
	P_{n,m}:=\begin{pmatrix}
		Q_{n,m}&0 \\
		0&Q_{n,m}
	\end{pmatrix},
\end{equation}
satisfying
\begin{equation} \label{eP30}
	\tilde{H}(\lambda ,P_{n,m}u)=\tilde{H}(\lambda,u),\ P_{n,m}^{T}J_{nm}P_{n,m}=J_{nm}\ {\rm and}\ P_{n,m}^{2m}=I_{2nm}.
\end{equation}

The following lemma is readily verified.
\begin{lem}\label{l2}
	Let $u$ and $	\tilde{H}$ be defined above, then $z(t)$ is a $2m\tau $-periodic  solution of second order  Hamiltonian  system \eqref{es} if and only if $u(t)$ is a  $2m\tau $-periodic  solution of first order Hamiltonian system \eqref{ef}. Moreover, it holds
	\begin{equation}\label{ne}
		\int_{0}^{\tau} \left [ \frac{1}{2}\left ( \mathcal{A}(\lambda )^{-1}\dot{z}(t), \dot{z}(t)  \right ) +H(\lambda ,z(t)) \right ] dt=\int_{0}^{\tau} \left [ \frac{1}{2}\left ( J_{nm}\dot{u}(t), u(t)  \right ) +\tilde{H} (\lambda ,u(t)) \right ] dt .
	\end{equation} 
\end{lem}

By Lemmas \ref{l1} and \ref{l2}, the transformations introduced above induce a bijection between the symmetric $2m\tau$-periodic solutions of \eqref{e2} and the $2m\tau$-periodic solutions of \eqref{ef}, preserving convergence to the trivial solution in $C^1_{\rm loc}$ and $\mathbb R$-distinctness. Consequently, $(\mu,0)$ is a bifurcation point of \eqref{e2} if and only if it is a bifurcation point of \eqref{ef}.

That is to say, we have successfully transfer  distributed differential delay system \eqref{e2} into Hamiltonian system \eqref{ef}  with a $P_{n,m}$-boundary condition, and there is a one-to-one correspondence of bifurcation points between the two systems. For each $\lambda \in \Lambda $, the function $\tilde{H}(\lambda ,u)$ in \eqref{Hp} is even with respect to $u$, thus $0\in \mathbb{R} ^{2nm} $ is the trivial solution of system \eqref{ef}. The linearization of the Hamiltonian system \eqref{ef} at the trivial solution $0$ is given by
\begin{equation}\label{ut2}
	\dot{u} (t)=J_{nm}\nabla ^{2}\tilde{H} (\lambda ,0) u(t), \quad \lambda\in \Lambda ,\quad u(t)\in \mathbb{R}^{2nm},
\end{equation}
where $$\nabla ^{2}\tilde{H} (\lambda ,0)=\begin{pmatrix}
		\nabla ^{2}H(\lambda,0)& 0\\
		0 &-\mathcal{A} (\lambda )
	\end{pmatrix}.$$
	 Let $\gamma_{\lambda}(t)$ be the fundamental solution of linear Hamiltonian system \eqref{ut2} with $\gamma_{\lambda}(0)=I_{2nm}$, then $$		\gamma_{\lambda}(t)=\exp\left ( tJ_{nm} \begin{pmatrix}
	\nabla ^{2}H(\lambda,0)& 0\\
	0 &-\mathcal{A} (\lambda )
	\end{pmatrix}\right )\in \mathcal{P}_{\tau } (2nm) .$$
	 For such symplectic path $\gamma _{\lambda}$, it assigns with a  Maslov $P_{n,m}$-index  $ i_{\tau}^{P_{n,m}}(\gamma _{\lambda})$ and a $P_{n,m}$-nullity $\nu_{\tau}^{P_{n,m}}(\gamma _{\lambda}) $ established in Definition \ref{nu}.
	 \subsection{Reduction to a Hamiltonian System with a $P_{n,2}$-Boundary Condition}
	 In this subsection, we consider the special case of two distributed delays with general nonlinear potentials. Unlike the multiple-delay system considered in the previous subsection, the special anti-periodic structure of the two-delay system allows us to extend the Hamiltonian reduction to the case of general nonlinear $G$. By introducing suitable auxiliary variables associated with the two distributed delay terms, we show that the nonlinear two-delay system can be equivalently transformed into a first-order Hamiltonian system with a \(P_{n,2}\)-boundary condition.

	 If $x(t)$ is a solution of system (\ref{e3'}), let
	 \begin{equation}\label{eP4}
	 	x_{1}(t)=x(t),\ x_{2}(t)=x(t-\tau ),
	 \end{equation}
	 and
	 \begin{equation}\label{eP5}
	 	z_{1}\left ( t \right )=\int_{0}^{\tau }\nabla G\left (\lambda, x\left ( t-s \right ) \right ) ds,\  z_{2}\left ( t \right )= \int_{\tau }^{2\tau }\nabla G\left (\lambda, x\left ( t-s \right ) \right )ds,
	 \end{equation}
	 then from $({\bf FG})$, it holds for $n\ge 1$ and each $\lambda\in \Lambda$
	 \begin{equation}\label{eP6}
	 \begin{pmatrix}
	 	\dot{x} _{1}(t) \\
	 	\dot{x} _{2}(t) \\
	 	\dot{z} _{1}(t)\\
	 	\dot{z} _{2}(t)
	 \end{pmatrix}=\frac{1}{\sqrt{2} } \begin{pmatrix}
	 	0&  0&  -I_{n}&-I_{n} \\
	 	0&  0&  I_{n}&-I_{n} \\
	 	I_{n}& -I_{n} & 0 & 0\\
	 	I_{n}&  I_{n}& 0 &0
	 \end{pmatrix}\begin{pmatrix}
	 	\sqrt{2}	\nabla G(\lambda,x_{1}(t))  \\
	 	\sqrt{2}	\nabla G(\lambda,x_{2}(t))  \\
	 	\sqrt{2}	\nabla F(\lambda,z_{1}(t)) \\
	 	\sqrt{2}	\nabla F(\lambda,z_{2}(t)) 
	 \end{pmatrix},
	 \end{equation}
and thus we can rewrite \eqref{eP6} as 
\begin{equation}\label{eP9}
	\dot{z } (t)=\mathcal{J}\nabla H(\lambda, z(t)), \quad z(t)\in \mathbb{R}^{4n},\quad n\ge 1,\quad  \lambda\in \Lambda,
\end{equation}
where 
	 \begin{equation}\label{eP10}
	 	z=\left ( x_{1},x_{2},z_{1},z_{2}\right )^{T}\in \mathbb{R}^{4n},\quad H(\lambda, z)=\sqrt{2}\left ( G(\lambda, x_{1})+G(\lambda, x_{2})+F(\lambda, z_{1})+F(\lambda, z_{2}) \right ),
	 \end{equation}
	 and  $\mathcal{J}$ is a $4n\times 4n$ skew-symmetric non-degenerate matrix defined by 
	 \begin{equation}\label{eP8}
	 	\mathcal{J}=\frac{1}{\sqrt{2} } \begin{pmatrix}
	 		0&  0&  -I_{n}&-I_{n} \\
	 		0&  0&  I_{n}&-I_{n} \\
	 		I_{n}& -I_{n} & 0 & 0\\
	 		I_{n}&  I_{n}& 0 &0
	 	\end{pmatrix},   
	 \end{equation}
	 satisfying $\mathcal{J}^{T}=-\mathcal{J}=\mathcal{J}^{-1}$. Moreover, given
	 \begin{equation}\label{eP12}
	 	\mathcal{M}=\begin{pmatrix}
	 		0&  -I_{n}& 0 &0 \\
	 		I_{n}& 0 & 0 &0 \\
	 		0& 0 & 0 & -I_{n}\\
	 		0& 0 & I_{n} &0
	 	\end{pmatrix} ,
	 \end{equation}
	 one can see that $\mathcal{M}$ is a $\mathcal{J}$-symplectic matrix satisfying $\mathcal{M}^{T}\mathcal{J}^{-1}\mathcal{M}=\mathcal{J}^{-1}$ with $\mathcal{M}^{4}=I_{4n}$, where
	 \begin{equation}\label{eP13}
	 	\mathcal{M}^{T}=\begin{pmatrix}
	 		0&  I_{n}& 0 &0 \\
	 		-I_{n}& 0 & 0 &0 \\
	 		0& 0 & 0 & I_{n}\\
	 		0& 0 & -I_{n} &0
	 	\end{pmatrix}=-\mathcal{M} 
	 	=\mathcal{M}^{-1},
	 \end{equation}
	 It is easily verified from (\ref{eP4}) and (\ref{eP5}) that
	 \begin{equation}\label{eP14}
	 	z(t+\tau )=\mathcal{M}z(t),
	 \end{equation}
	 which is equivalent to
	 \begin{equation}\label{eP15}
	 	z(\tau )=\mathcal{M}z(0).
	 \end{equation}
	 
	 The following lemma can be proved by a standard argument similar to that used in \cite{He2}.
	 \begin{lem}\label{l3}
	 	Let $z$ and $H $ be defined in \eqref{eP10}, if $x(t)$ is a symmetric $4\tau $-periodic solution of distributed delay differential system \eqref{e3'}, then  $z(t)$ is the $4\tau $-periodic  solution of the following Hamiltonian  system with $\mathcal{M}$-boundary condition,
	 	\begin{equation} \label{eP16}
	 		\begin{cases}
	 			\dot{z } (t)=\mathcal{J}\nabla H(\lambda, z(t)), \\
	 			z(t+\tau )=\mathcal{M}z(t).
	 		\end{cases}
	 	\end{equation}
	 	Conversely, if $z(t)=\left ( x_{1}(t),x_{2}(t),z_{1}(t),z_{2}(t)\right )^{T}\in \mathbb{R}^{4n}$ is a $4\tau $-periodic solution of Hamiltonian system \eqref{eP16}, then $x(t):=x_{1}(t)$ is a symmetric $4\tau $-periodic solution of distributed delay differential system \eqref{e3'}.
	 \end{lem}
	 Furthermore, denote 
	 \begin{equation}
	 u(t):=S^{-1}z(t),\quad 	S=\begin{pmatrix} \frac{1}{\sqrt{2} }I_{n} &\frac{1}{\sqrt{2} }I_{n} &0 &0 \\ -\frac{1}{\sqrt{2} }I_{n}& \frac{1}{\sqrt{2} }I_{n} & 0 &0 \\ 0 & 0 & I_{n} & 0\\ 0 & 0 & 0 &I_{n} \end{pmatrix},
	 \end{equation}
then it holds $SS^{T}=I_{4n}$ and $S^{T}\mathcal{J}S=J_{2n} $, where $J_{2n}$ is the standard symplectic matrix $\begin{pmatrix}
	0&-I_{2n} \\
	I_{2n}&0
\end{pmatrix}$. 
Since $S^{-1}\mathcal{M}S=\mathcal{M}$ and notice that $\mathcal{M}$ actually coincides with  the matrix $P_{n,2}=\begin{pmatrix}
	Q_{n,2}&0 \\
	0&Q_{n,2}
\end{pmatrix}$ defined in \eqref{Hp}, system  \eqref{eP16}  can be transformed into 
\begin{equation}\label{e5}
	\left\{\begin{matrix}
		\dot{u}(t)=J_{2n}\nabla \tilde{H}(\lambda ,u(t)),   \\
		u(t+\tau)=P_{n,2}u(t),\quad \
	\end{matrix}\right.
\end{equation}
with 
\begin{equation}\label{H3}
\tilde{H} (\lambda ,u):=H(\lambda ,Su),\quad \forall \left ( \lambda ,u \right ) \in \Lambda \times \mathbb{R}^{4n} ,
\end{equation}
satisfying 
\begin{equation}
	\tilde{H} (\lambda ,P_{n,2}u)=\tilde{H} (\lambda ,u),\ P_{n,2}^{T}J_{2n}P_{n,2}=J_{2n}\ {\rm and}\ P_{n,2}^{4}=I_{4n}.
\end{equation}

The following lemma is readily verified.
 \begin{lem}\label{l4} Let $u$ and $	\tilde{H}$ be defined above, then $z(t)$ is a $4\tau $-periodic  solution of  Hamiltonian  system \eqref{eP16} if and only if $u(t)$ is a  $4\tau $-periodic  solution of Hamiltonian system \eqref{e5}. Moreover, it holds
 	\begin{equation}\label{ne2}
 		\int_{0}^{\tau} \left [ \frac{1}{2}\left ( \mathcal{J} \dot{z}(t), z(t)  \right ) +H(\lambda ,z(t)) \right ] dt=\int_{0}^{\tau} \left [ \frac{1}{2}\left ( J_{2n}\dot{u}(t), u(t)  \right ) +\tilde{H} (\lambda ,u(t)) \right ] dt .
 	\end{equation} 
 	\end{lem}
 	
 	Now, we successfully transfer the original two-delay system \eqref{e3'} into a standard first order Hamiltonian system \eqref{e5} with  $P_{n,2}$-boundary condition, and there is a one-to-one correspondence of bifurcation points between the two systems. For each $\lambda \in \Lambda $, the function $\tilde{H}(\lambda ,u)$ in \eqref{H3} is even with respect to $u$, thus $0\in \mathbb{R} ^{4n} $ is the trivial solution of system \eqref{e5}. The linearization of the Hamiltonian system \eqref{e5} at the trivial solution $0$ is given by
 	\begin{equation}\label{ut21}
 		\dot{u} (t)=J_{2n}\nabla ^{2}\tilde{H} (\lambda ,0) u(t), \quad \lambda\in \Lambda ,\quad u(t)\in \mathbb{R}^{4n},
 	\end{equation}
 	where $$\nabla ^{2}\tilde{H} (\lambda ,0)=S^{T}\nabla ^{2}H(\lambda ,0)S=\begin{pmatrix}
 		I_{2}\otimes \left ( \sqrt{2}\nabla ^{2}G(\lambda ,0)  \right ) & 0\\
 		0 &I_{2}\otimes \left ( \sqrt{2}\nabla ^{2}F(\lambda ,0)  \right )
 	\end{pmatrix}.$$
 	Here, $\otimes$ denotes the Kronecker product of matrices.
 	Let $\gamma_{\lambda}(t)$ be the fundamental solution of linear Hamiltonian system \eqref{ut21} with $\gamma_{\lambda}(0)=I_{4n}$, then $$		\gamma_{\lambda}(t)=\exp\left ( tJ_{2n} \begin{pmatrix}
 		I_{2}\otimes \left ( \sqrt{2}\nabla ^{2}G(\lambda ,0)  \right ) & 0\\
 		0 &I_{2}\otimes \left ( \sqrt{2}\nabla ^{2}F(\lambda ,0)  \right )
 	\end{pmatrix}\right )\in \mathcal{P}_{\tau } (4n) .$$
 	For such symplectic path $\gamma _{\lambda}$, it assigns with a  Maslov type $P_{n,2}$-index  $ i_{\tau}^{P_{n,2}}(\gamma _{\lambda})$ and a $P_{n,2}$-nullity $\nu_{\tau}^{P_{n,2}}(\gamma _{\lambda}) $ established in Definition \ref{nu}.
 	
\section{Proofs of the Main Results}\label{sec3}
We begin with two auxiliary lemmas that will be used in the proofs of our main results.
\begin{lem}\label{abv}
	For $a,b\in \mathbb{R} $, let $\gamma _{a,b}$ be the fundamental matrix solution of 
	\begin{equation}
		\dot{u}(t)=J_{m}\begin{pmatrix}
			aI_{m}&0 \\
			0&bI_{m}
		\end{pmatrix} u(t),
	\end{equation} 
	with $\gamma _{a,b}(0)=I_{2m}$, that is, $\gamma _{a,b}(t)=\exp \left ( tJ_{m}\begin{pmatrix}
		aI_{m}&0 \\
		0&bI_{m}
	\end{pmatrix} \right ) $.  Then for $P_{1,m}:=\begin{pmatrix}
	Q_{1,m}&0 \\
	0&Q_{1,m}
	\end{pmatrix}$, where $Q_{1,m}$ is the $m\times m$ matrix given in \eqref{418} with $n=1$, it holds
	\begin{enumerate}[{\rm (I)}]
		\item  If $ab>0$, then we have 
		\begin{equation}\label{nuhg}
			\nu_{\tau} ^{P_{1,m}}(\gamma _{a,b})=\begin{cases}
				2,\quad \quad \mathrm{if} \  m\sqrt{ab} \tau \in  \left ( 2\mathbb{N}-1  \right ) \pi  , \\
				0,\quad \quad   \mathrm{otherwise}.
			\end{cases}
		\end{equation} 
		In particular, if $m\sqrt{ab} \tau<\pi$, then $\nu_{\tau} ^{P_{1,m}}(\gamma _{a,b})=0$. Furthermore
		\begin{enumerate}[{\rm (i)}]
			\item If $a>0$, then 
			\begin{equation}
				i_{\tau} ^{P_{1,m}}(\gamma _{a,b})=2E\left ( \frac{m\sqrt{ab}\tau }{2\pi }-\frac{1}{2}   \right )+2=2\rho ,
			\end{equation}
			where $	\rho \in \mathbb{Z} $ is the unique integer such that $(2\rho -1)\pi < m\sqrt{ab}\tau \le (2\rho +1)\pi $. Here, $E\left ( a \right )=\max \left \{ k\in \mathbb{Z} |k< a \right \}  $.
			\item If $a<0$, then
			\begin{equation}
				i_{\tau} ^{P_{1,m}}(\gamma _{a,b})=2E\left (- \frac{m\sqrt{ab}\tau }{2\pi }-\frac{1}{2}   \right )+2=2\bar{\rho } ,
			\end{equation}
			where $	\bar{\rho } \in \mathbb{Z} $ is the unique integer such that $(2\bar{\rho } -1)\pi < -m\sqrt{ab}\tau \le (2\bar{\rho } +1)\pi $. 
		\end{enumerate}
		\item  If $ab\le 0$, then 
		\begin{equation}
			i_{\tau} ^{P_{1,m}}(\gamma _{a,b})=0,\quad \nu_{\tau}^{P_{1,m}}(\gamma _{a,b})=0.
		\end{equation}
	\end{enumerate}
	
\end{lem}
\begin{proof}
	We next proceed by considering three cases.
	
	{\bf Case 1}( $ab>0$). In this case,  we discuss $a>0$ and $a<0$ separately. First, suppose $a>0$ and set $r=\left ( \frac{b}{a} \right )^{\frac{1}{4} }  $, then it holds $ar^{2}=br^{-2}=\sqrt{ab} $. Define the $2m\times 2m$ symmetric matrix 
	\begin{equation}\label{lj}
		L=\begin{pmatrix}
			rI_{m}&0 \\
			0 &r^{-1}I_{m}
		\end{pmatrix},
	\end{equation}
	then a simple computation shows that  $L^{T}J_{m}L=J_{m}$, i.e., $L\in \mathrm{Sp}(2m)$, and it satisfies 
	\begin{equation}\label{LJL}
		L^{-1}J_{m}\begin{pmatrix}
			aI_{m}&0 \\
			0  &bI_{m}
		\end{pmatrix}L=J_{m}L^{T}\begin{pmatrix}
			aI_{m}&0 \\
			0  &bI_{m}
		\end{pmatrix}L=J_{m}\begin{pmatrix}
			ar^{2}I_{m}&0 \\
			0  &br^{-2}I_{m}
		\end{pmatrix}=\sqrt{ab}J_{m},
	\end{equation}
	as well as 
	\begin{equation}\label{LPL}
		L^{-1}P_{1,m}L=\begin{pmatrix}
			r^{-1}I_{m}&0 \\
			0 &rI_{m}
		\end{pmatrix}\begin{pmatrix}
			Q_{1,m}&0 \\
			0&Q_{1,m}
		\end{pmatrix}\begin{pmatrix}
			rI_{m}&0 \\
			0 &r^{-1}I_{m}
		\end{pmatrix}=P_{1,m}.
	\end{equation}
	Moreover, note that the unitary matrix $Q_{1,m}$ given in \eqref{418} with $n=1$ satisfies $Q_{1,m}^{m}=-I_{m}$ and that its characteristic polynomial is $t^{m}+1$, then its eigenvalues are precisely  $\exp \left ( (2k-1)\pi i/m\right ),k=1,\cdots ,m $, and there exists $\mathcal{Q}\in \mathrm{Sp}(2m)\cap \mathrm{O}(2m) $ such that 
	\begin{equation}\label{Qm}
		\mathcal{Q}^{-1}P_{1,m}\mathcal{Q} =R(\alpha _{1})\Diamond \cdots \Diamond R(\alpha _{m}),
	\end{equation}
	where $\alpha_{k}:=(2k-1)\pi/m, k=1,\cdots ,m  $, and $R(\alpha_{k})$ is the rotating matrix defined by 
	\begin{equation}
		R(\alpha_{k})=\begin{pmatrix}
			\cos \alpha_{k} & -\sin \alpha_{k}\\
			\sin \alpha_{k} &\cos \alpha_{k}
		\end{pmatrix}.
	\end{equation}
	
	Now, let $\tilde{L}:=L\mathcal{Q}  $, then from \eqref{LJL}, \eqref{LPL} and \eqref{Qm},  $\tilde{L}\in \mathrm{Sp}(2m)$ satisfies 
	\begin{equation}\label{TLPL}
		\tilde{L} ^{-1}P_{1,m}\tilde{L}=\mathcal{Q}^{-1}L^{-1}P_{1,m}L\mathcal{Q} =\mathcal{Q}^{-1}P_{1,m}\mathcal{Q}=R(\alpha _{1})\Diamond \cdots \Diamond R(\alpha _{m}), 
	\end{equation}
	and 
	\begin{equation}\label{Lj}
		\begin{split}
			\tilde{L}^{-1}J_{m}\begin{pmatrix}
				aI_{m}&0 \\
				0  &bI_{m}
			\end{pmatrix}\tilde{L}=&\mathcal{Q}^{-1}L^{-1}J_{m}\begin{pmatrix}
				aI_{m}&0 \\
				0  &bI_{m}
			\end{pmatrix}L\mathcal{Q}\\
			=&\sqrt{ab}\mathcal{Q}^{-1}J_{m}\mathcal{Q}\\
			=&\sqrt{ab}J_{m}.
		\end{split}
	\end{equation}
	Let $\hat{\gamma }_{a,b}(t)$ be the fundamental  matrix solution of the linear Hamiltonian system
	\begin{equation}
		\dot{u} (t)=\sqrt{ab}J_{m}u(t),\quad u(t)\in \mathbb{R} ^{2m},
	\end{equation}
	with $\hat{\gamma }_{a,b}(0)=I_{2m}$, then from \eqref{Lj}, it holds
	\begin{equation}\label{HGab}
		\hat{\gamma }_{a,b}(t)=\tilde{L}^{-1}	\gamma _{a,b}(t)\tilde{L}.
	\end{equation}
	For $c\in \mathbb{R} $, set $\gamma _{c}(t)=R(ct)=\exp \left ( tcJ_{1} \right ) $,  then from Lemma \ref{pDia} (iv) we deduce 
	\begin{equation}\label{HG}
		\hat{\gamma }_{a,b}(t)=\exp \left ( t\sqrt{ab}J_{m} \right ) =\exp \left ( \left (t\sqrt{ab} J_{1} \right )^{\Diamond m}  \right )=\left ( \exp \left ( t\sqrt{ab} J_{1}  \right )  \right )^{\Diamond m}=R(\sqrt{ab}t)^{\Diamond m}=\gamma _{\sqrt{ab}}(t)^{\Diamond m}, 
	\end{equation}
	where $\Diamond m$ denotes $m$-fold $\Diamond$-sum, and therefore, \eqref{HG} together with \eqref{TLPL} and Lemma \ref{pDia} (iii) imply  
	\begin{equation}\label{HGab2}
		\hat{\gamma }_{a,b}(\tau ) -R(\alpha _{1})\Diamond \cdots \Diamond R(\alpha _{m})=\left ( \gamma _{\sqrt{ab}}(\tau )-R(\alpha _{1}) \right )\Diamond \cdots \Diamond \left ( \gamma _{\sqrt{ab}}(\tau )-R(\alpha _{m}) \right ).
	\end{equation}
	Combining \eqref{TLPL}, \eqref{HGab}, \eqref{HG} and \eqref{HGab2}, we infer 
	\begin{equation}\label{ijsv}
		\begin{split}
			\nu_{\tau} ^{P_{1,m}}(\gamma _{a,b})= 	\nu_{\tau} ^{\tilde{L} ^{-1}P_{1,m}\tilde{L}}(\tilde{L}^{-1}	\gamma _{a,b}\tilde{L})=\nu_{\tau} ^{R(\alpha _{1})\Diamond \cdots \Diamond R(\alpha _{m})}(\gamma _{\sqrt{ab}}^{\Diamond m})=\sum_{k=1}^{m} \nu_{\tau} ^{R(\alpha _{k})}\left ( \gamma _{\sqrt{ab}}\right ).
		\end{split}
	\end{equation}
	For each $k=1,\cdots , m$, from Definition \ref{nu}, it holds
	\begin{equation}
		\nu_{\tau} ^{R(\alpha _{k})}\left ( \gamma _{\sqrt{ab}}\right )=\dim \ker \left (  \gamma _{\sqrt{ab}}\left ( \tau  \right )- R(\alpha _{k}) \right )
		=\dim\ker \left (R(\sqrt{ab}\tau -\alpha _{k})-I_{2}  \right ),    
	\end{equation}
	and by a straightforward computation, we have 
	$$\det \left (R(\sqrt{ab}\tau -\alpha _{k})-I_{2}  \right )=2-2\cos \left ( \sqrt{ab}\tau -\alpha _{k} \right ) .$$
	Consequently, for each $k=1,\cdots , m$, we deduce 
	\begin{equation}\label{dnu}
		\nu_{\tau} ^{R(\alpha _{k})}\left ( \gamma _{\sqrt{ab}}\right )=\begin{cases}
			2,\quad \quad \mathrm{if} \  m\sqrt{ab}\tau =(2k+2ml-1)\pi \ \mathrm{for \ some}\  l\in  \mathbb{Z}, \\
			0,\quad \quad   \mathrm{otherwise}.
		\end{cases}
	\end{equation}
	
	Next, we claim that  at most one \(k\in\{1,\cdots,m\}\) satisfies the condition in \eqref{dnu}. Suppose that there exist \(k_1,k_2\in\{1,\cdots,m\}\) and \(l_1,l_2\in\mathbb Z\) such that
	\[
	m\sqrt{ab}\tau
	=(2k_1+2ml_1-1)\pi
	=(2k_2+2ml_2-1)\pi.
	\]
	Then
	\[
	k_1-k_2=m(l_2-l_1).
	\]
	Since \(1\le k_1,k_2\le m\), we have
	\[
	|k_1-k_2|\le m-1.
	\]
	Consequently, \(k_1-k_2\) must be zero, and hence \(k_1=k_2\). Thus, at most one \(k\in\{1,\cdots,m\}\) satisfies \eqref{dnu}, and therefore one can easily verify from \eqref{ijsv} that 
	\begin{equation}\label{nuhg}
		\nu_{\tau} ^{P_{1,m}}(\gamma _{a,b})=\begin{cases}
			2,\quad \quad \mathrm{if} \  m\sqrt{ab}\tau \in  \left ( 2\mathbb{N}-1  \right ) \pi  , \\
			0,\quad \quad   \mathrm{otherwise}.
		\end{cases}
	\end{equation}
	
	Argued as \eqref{ijsv}, from Proposition \ref{Sa} and Proposition \ref{Si}, we derive
	\begin{equation}\label{iP}
		i_{\tau} ^{P_{1,m}}(\gamma _{a,b})= 	i_{\tau}^{\tilde{L} ^{-1}P_{1,m}\tilde{L}}(\tilde{L}^{-1}	\gamma _{a,b}\tilde{L})=i_{\tau} ^{R(\alpha _{1})\Diamond \cdots \Diamond R(\alpha _{m})}(\gamma _{\sqrt{ab}}^{\Diamond m})=\sum_{k=1}^{m} i_{\tau} ^{R(\alpha _{k})}\left ( \gamma _{\sqrt{ab}}\right ).
	\end{equation}
	To compute the index $i_{\tau}^{R(\alpha _{k})}\left ( \gamma _{\sqrt{ab}}\right ) $, for each $k=1,\cdots , m$, we choose the path $\xi (t)=R(-\alpha _{k}t/\tau )$, then it holds $\xi (0)=I_{2}$ and $ \xi (\tau )=R(-\alpha _{k})$. Note that $R(-\alpha _{k})\gamma _{\sqrt{ab}}(t)=R\left ( \sqrt{ab}t-\alpha _{k} \right )$, then from \eqref{g12} we obtain 
	\begin{equation}\label{p1}
		R(-\alpha _{k})\gamma _{\sqrt{ab}}\ast \xi (t)=\left\{\begin{matrix}
			R(-2\alpha _{k}t/\tau ),\quad \quad \quad \quad \quad \  0\le t\le \tau /2, \\
			R\left ( \sqrt{ab}(2t-\tau )-\alpha _{k} \right ),\quad  \tau /2\le t\le \tau.
		\end{matrix}\right.
	\end{equation}
	For each $k=1,\cdots , m$, define the path 
	\begin{equation}\label{p2}
		\gamma _{c_{k}}(t)=R\left (c_{k}t  \right )= R\left ( \frac{\sqrt{ab}\tau -\alpha _{k}}{\tau } t \right ),\quad t\in \left [ 0,\tau  \right ] ,
	\end{equation}
	where $c_{k}:= \frac{\sqrt{ab}\tau -\alpha _{k}}{\tau }$, then clearly the two paths defined in \eqref{p1} and \eqref{p2} are homotopic in $\mathcal{P}_{\tau } (2)$ with fixed endpoints. Therefore, it follows from the homotopy invariance of the Maslov-type index $i_{\tau}$ that 
	\begin{equation}\label{icjk}
		i_{\tau }\left ( R(-\alpha _{k})\gamma _{\sqrt{ab}}\ast \xi \right ) =i_{\tau }\left ( \gamma _{c_{k}} \right ) .
	\end{equation}
	Now, take $n=1$, $\gamma (t)=\gamma _{c_{k}}(t)$ and $\theta _{1}(t)=c_{k}t$ in Theorem \ref{po}, we deduce 
	\begin{equation}\label{cjk}
		i_{\tau }\left ( \gamma _{c_{k}} \right ) =2E\left ( \frac{\sqrt{ab}\tau -\alpha _{k}}{2\pi }  \right )+1.
	\end{equation}
	Recall that $E\left ( a \right )=\max \left \{ k\in \mathbb{Z} |k< a \right \}  $ and apply Theorem \ref{po} again to the path $\xi (t)=R(-\alpha _{k}t/\tau )$, then it holds 
	\begin{equation}\label{xk}
		i_{\tau }(\xi )=2E\left ( -\frac{ \alpha _{k}}{2\pi }  \right )+1=-1,
	\end{equation}
	where the second equality follows from the fact that $-1< -\alpha _{k}/2\pi  < 0$. For each $k=1,\cdots , m$, from \eqref{icjk}, \eqref{cjk}, \eqref{xk} and Definition \ref{nu}, we thus obtain 
	\begin{equation}\label{327}
		i_{\tau}^{R(\alpha _{k})}\left ( \gamma _{\sqrt{ab}}\right )=i_{\tau }\left ( R(-\alpha _{k})\gamma _{\sqrt{ab}}\ast \xi \right )-	i_{\tau }(\xi )=2E\left ( \frac{\sqrt{ab}\tau -\alpha _{k}}{2\pi }  \right )+2.
	\end{equation}
The integer-part function $E$ satisfies the following identity
	\[
	\sum_{k=0}^{m-1} E\left(s-\frac{k}{m}\right)
	=E(ms)-m+1,
	\qquad s\in\mathbb{R},\quad m\in\mathbb{N}.
	\]
	Combining \eqref{iP} and \eqref{327} and applying this identity with
	$s=\frac{\sqrt{ab}\,\tau}{2\pi}-\frac{1}{2m}$, we obtain
	\begin{equation}
		\begin{split}
			i_{\tau} ^{P_{1,m}}(\gamma _{a,b})=&\sum_{k=1}^{m} i_{\tau} ^{R(\alpha _{k})}\left ( \gamma _{\sqrt{ab}}\right ) \\
			=&2\sum_{k=1}^{m} E\left ( \frac{\sqrt{ab}\tau -\alpha _{k}}{2\pi }  \right )+2m\\
			=&2\sum_{k=1}^{m} E\left ( \frac{\sqrt{ab}\tau }{2\pi }-\frac{2k-1}{2m}   \right )+2m\\
			=&2E\left ( \frac{m\sqrt{ab}\tau }{2\pi }-\frac{1}{2}   \right )+2\\
			=&2\rho ,
		\end{split}
	\end{equation}
	where $	\rho=E\left ( \frac{m\sqrt{ab}\tau }{2\pi }-\frac{1}{2}   \right )+1 \in \mathbb{Z} $ is the unique integer such that $(2\rho -1)\pi < m\sqrt{ab}\tau \le (2\rho +1)\pi $. 
	
	We next turn to the case $a<0$,  and again set  $r=\left ( \frac{b}{a} \right )^{\frac{1}{4} }  $. Since $a<0$, it holds $$ar^{2}=br^{-2}=-\sqrt{ab}. $$ Thus, replacing $\sqrt{ab}$ by $-\sqrt{ab}$ in the preceding argument for $a>0$, we can similarly obtain 
	\begin{equation}\label{nuhg2}
		\nu_{\tau} ^{P_{1,m}}(\gamma _{a,b})=\begin{cases}
			2,\quad \quad \mathrm{if} \  m\sqrt{ab}\tau \in  \left ( 2\mathbb{N}-1  \right ) \pi  , \\
			0,\quad \quad   \mathrm{otherwise}.
		\end{cases}
	\end{equation}
	and 
	\begin{equation}
		i_{\tau} ^{P_{1,m}}(\gamma _{a,b})=2E\left (- \frac{m\sqrt{ab}\tau }{2\pi }-\frac{1}{2}   \right )+2=2\bar{\rho } ,
	\end{equation}
	where $	\bar{\rho } \in \mathbb{Z} $ is the unique integer such that $(2	\bar{\rho } -1)\pi < -m\sqrt{ab}\tau \le (2	\bar{\rho } +1)\pi $. 
	
	{\bf Case 2}( $ab<0$). Denote 
	\begin{equation}\label{Aj2}
		A=J_{m}\begin{pmatrix}
			a  I_{m}&0 \\
			0 &b I_{m}
		\end{pmatrix} =\begin{pmatrix}
			0&- b I_{m}\\
			a I_{m} &0
		\end{pmatrix}\in \mathbb{R} ^{2m\times 2m},
	\end{equation}
	then it holds 
	\begin{equation}\label{Aj}
		A^{2}=-abI_{2m}=\tilde{c}  ^{2}I_{2m},
	\end{equation}
	where we denote $\tilde{c}:= \sqrt{-ab}  >0$. From \eqref{Aj} and the definition of  hyperbolic cosine and hyperbolic sine, we deduce
	\begin{equation}
		\begin{split}
			\gamma _{a,b}(t)=&\exp \left ( tA \right )\\
			=&\sum_{i=0}^{\infty }\frac{t^{i}A^{i}}{i!} \\
			=&\sum_{i=0}^{\infty }\frac{t^{2i}A^{2i}}{\left ( 2i \right )! } +\sum_{i=0}^{\infty }\frac{t^{2i+1}A^{2i+1}}{\left ( 2i +1\right )! }\\ =&\sum_{i=0}^{\infty }\frac{t^{2i}\tilde{c}  ^{2i}}{\left ( 2i \right )! }I_{2m} +\frac{1}{\tilde{c}} \sum_{i=0}^{\infty }\frac{t^{2i+1}\tilde{c}^{2i+1}}{\left ( 2i +1\right )! } A \\
			=&\cosh \left ( t\tilde{c} \right )I_{2m}+\frac{1}{\tilde{c}} \sinh \left ( t\tilde{c} \right ) A,
		\end{split}
	\end{equation}
	therefore we have 
	\begin{equation}
		\gamma _{a,b}(\tau )-P_{1,m}=\begin{pmatrix}
			\cosh \left ( \tau \tilde{c} \right )I_{m}-Q_{1,m} & -\frac{b}{\tilde{c}} \sinh \left ( \tau \tilde{c} \right ) I_{m}\\
			\frac{ 	a }{\tilde{c}} \sinh \left ( \tau \tilde{c} \right ) I_{m} &\cosh \left ( \tau \tilde{c} \right )I_{m}-Q_{1,m}
		\end{pmatrix},
	\end{equation}
	where $Q_{1,m}$ is the $m\times m$ matrix given in \eqref{418} with $n=1$. Note that the $m\times m$ matrices $\cosh \left ( \tau \tilde{c} \right )I_{m}-Q_{1,m}$ and $\frac{ 	a }{\tilde{c}} \sinh \left ( \tau \tilde{c}  \right ) I_{m}$ commute, then by Lemma 3 in \cite{Long} and the identities $\cosh \left ( \tau \tilde{c} \right )+ \sinh \left ( \tau \tilde{c}  \right )=\exp \left (  \tau \tilde{c} \right )$ and $\cosh \left ( \tau \tilde{c}\right )- \sinh \left ( \tau \tilde{c} \right )=\exp \left ( - \tau \tilde{c} \right )$, it holds
	\begin{equation}\label{glp}
		\begin{split}
			\det \left ( \gamma _{a,b}(\tau )-P_{1,m} \right ) =&\det \left [ \left ( 	\cosh \left ( \tau \tilde{c} \right )I_{m}-Q_{1,m}  \right )^{2} -\left (  \sinh \left ( \tau \tilde{c} \right ) I_{m} \right )^{2}  \right ] \\
			=&\det \left [ \left ( \exp \left (  \tau \tilde{c}  \right )I_{m}-Q_{1,m} \right ) \left ( \exp \left (  -\tau \tilde{c}  \right )I_{m}-Q_{1,m} \right )  \right ]\\
			=&\det \left ( \exp \left (  \tau \tilde{c} \right )I_{m}-Q_{1,m}  \right ) \det \left ( \exp \left (  -\tau \tilde{c}  \right )I_{m}-Q_{1,m} \right )\\
			=&\left (  \exp \left (  m\tau \tilde{c}  \right )+1 \right ) \left (  \exp \left ( - m\tau \tilde{c} \right )+1 \right ) \\
			=&2\cosh \left (  m\tau \tilde{c}  \right )+2 >0,
		\end{split}
	\end{equation}
	which shows that
	\begin{equation}
		\nu_{\tau} ^{P_{1,m}}(\gamma _{a,b})=0.
	\end{equation}
	For $\left ( s,t \right ) \in \left [ 0,1 \right ] \times \left [ 0,\tau  \right ] $, define the map
	\begin{equation}
		\delta \left ( s,t \right ) =	\gamma _{a,b}(st )=\begin{pmatrix}
			\cosh \left ( st \tilde{c} \right )I_{m} & -\frac{b}{\tilde{c}} \sinh \left ( st\tilde{c}  \right ) I_{m}\\
			\frac{ 	a }{\tilde{c}} \sinh \left ( st \tilde{c}  \right ) I_{m} &\cosh \left ( st\tilde{c} \right )I_{m},
		\end{pmatrix},
	\end{equation}
	then it holds $\delta_{0}(t):=\delta \left ( 0,t \right )=I_{2m}$ and $\delta_{1}(t):=\delta \left ( 1,t \right )=	\gamma _{a,b}(t )$. Moreover, by \eqref{glp}, we have, for every $s\in \left [ 0,1 \right ] $, 
	\begin{equation}
		\det \left ( \delta \left ( s,\tau  \right )-P_{1,m}  \right )=\det \left ( 	\gamma _{a,b}(s\tau  )-P_{1,m}  \right )=2\cosh \left (  ms\tau \tilde{c}  \right )+2 >0,
	\end{equation}
	and consequently
	\begin{equation}
		\dim \ker  \left ( \delta \left ( s,\tau  \right )-P_{1,m}  \right ) =0,\quad s\in\left [ 0,1 \right ],
	\end{equation}
	which implies $\delta _{0}\sim _{P_{1,m}}\delta _{1}$.  From Definition \ref{nu}, we have $i_{\tau}^{P_{1,m}}\left ( I_{2m} \right ) =0$, then it follows from Proposition \ref{hi} that 
	\begin{equation}
		i_{\tau}^{P_{1,m}}\left (\gamma _{a,b} \right )=i_{\tau}^{P_{1,m}}\left (\delta  _{1} \right )=i_{\tau}^{P_{1,m}}\left ( \delta  _{0} \right )= i_{\tau}^{P_{1,m}}\left ( I_{2m} \right ) =0.
	\end{equation}
	
	{\bf Case 3}( $ab=0$). If $a=0$ and $b\ne0$, the $2m\times 2m$ matrix $A$ in \eqref{Aj2} is given by  $	\begin{pmatrix}
		0&- b I_{m}\\
		0 &0
	\end{pmatrix}$, then we have 
	\begin{equation}
		\gamma _{a,b}(t)=\exp \left ( tA \right )=\begin{pmatrix}
			I_{m}&-btI_{m} \\
			0&I_{m}
		\end{pmatrix},
	\end{equation}
	and 
	\begin{equation}
		\gamma _{a,b}(t)-P_{1,m}=\begin{pmatrix}
			I_{m}-Q_{1,m}&-btI_{m} \\
			0&I_{m}-Q_{1,m}
		\end{pmatrix}.
	\end{equation}
	Therefore, it holds 
	\begin{equation}\label{lhg}
		\det \left (\gamma _{a,b}(t)-P_{1,m} \right )= \det \left ( I_{m}-Q_{1,m} \right )  \det \left ( I_{m}-Q_{1,m} \right )=4>0 ,
	\end{equation}
	for any $t\in \left [ 0,\tau  \right ] $. For $\left ( s,t \right ) \in \left [ 0,1 \right ] \times \left [ 0,\tau  \right ] $, define the map
	\begin{equation}
		\delta \left ( s,t \right ) =	\gamma _{a,b}(st )=\begin{pmatrix}
			I_{m}&-bstI_{m} \\
			0&I_{m}
		\end{pmatrix},
	\end{equation}
	then it satisfies 
	\begin{equation}
		\det \left (	\delta \left ( s,t \right ) -P_{1,m} \right )= \det \left ( I_{m}-Q_{1,m} \right )  \det \left ( I_{m}-Q_{1,m} \right )=4>0 ,
	\end{equation}
	for any $\left ( s,t \right ) \in \left [ 0,1 \right ] \times \left [ 0,\tau  \right ] $. Consequently, by arguing similarly as {\bf Case 2}, we obtain 
	\begin{equation}\label{inP}
		i_{\tau} ^{P_{1,m}}(\gamma _{a,b})=0,\quad	\nu_{\tau} ^{P_{1,m}}(\gamma _{a,b})=0.
	\end{equation}
	An analogous argument shows that \eqref{inP} also holds in the case $a\ne0, b=0$ and $a=b=0$, hence we omit the details.
\end{proof}
We next extend Lemma \ref{abv} to the matrix setting. The key observation is that two commuting real symmetric matrices can be simultaneously orthogonally diagonalized, reducing the problem to the scalar case.
\begin{lem}\label{abv2}
	Let $A$ and $B$ be two commuting real symmetric matrices of order $n$, then there exists an orthogonal matrix $T$ of order $n$ such that  
	\begin{equation}\label{T}
		T^{-1}AT={\rm diag}(a_{1},\cdots, a_{n}), \quad	T^{-1}BT={\rm diag}(b_{1},\cdots, b_{n}).
	\end{equation}
	Denote by $\gamma _{A,B}$  the fundamental matrix solution of 
	\begin{equation}
		\dot{u}(t)=J_{nm}\begin{pmatrix}
			I_{m}\otimes A&0 \\
			0&I_{m}\otimes B
		\end{pmatrix} u(t),
	\end{equation} 
	with $\gamma _{A,B}(0)=I_{2nm}$, here $\otimes$ denotes the Kronecker product of matrices.  Then for $P_{n,m}=\begin{pmatrix}
		Q_{n,m} &0 \\
		0& Q_{n,m}
	\end{pmatrix}$ given in \eqref{Hp}, it holds
	\begin{equation}\label{nuP}
		\nu_{\tau} ^{P_{n,m}}(\gamma _{A,B})=2\sharp \left \{ j\mid m\sqrt{a_{j}b_{j}} \tau \in  \left ( 2\mathbb{N}-1  \right ) \pi  \right \}, 
	\end{equation}
	\begin{equation}\label{iPr}
		i_{\tau} ^{P_{n,m}}(\gamma _{A,B})=\sum_{a_{j}>0,b_{j}>0}2\rho _{j}+\sum_{a_{j}<0,b_{j}<0}2\bar{\rho}  _{j},
	\end{equation}
	where for $j\in\left \{ 1,\cdots ,n \right \} $, $	\rho_{j} \in \mathbb{Z} $ is the unique integer such that $(2\rho_{j} -1)\pi <m \sqrt{a_{j}b_{j}}\tau \le (2\rho_{j} +1)\pi $, and $	\bar{\rho }_{j} \in \mathbb{Z} $ is the unique integer such that $(2\bar{\rho }_{j} -1)\pi < -m\sqrt{a_{j}b_{j}}\tau \le (2\bar{\rho }_{j} +1)\pi $. 
	
	Furthermore, if $A$ and $B$ are positive definite, then 
	\begin{equation}\label{iPAB1}
		i _{\tau }^{P_{n,m}}(\gamma _{A,B})=\sum_{0< t<\tau }2\sharp \left \{ j\mid m\sqrt{a_{j}b_{j}} t \in  \left ( 2\mathbb{N}-1  \right ) \pi  \right \},
	\end{equation}
	and if $A$ and $B$ are negative definite, then
	\begin{equation}\label{iPAB2}
		i _{\tau }^{P_{n,m}}(\gamma _{A,B})=-\sum_{0< t\le \tau }2\sharp \left \{ j\mid m\sqrt{a_{j}b_{j}} t \in  \left ( 2\mathbb{N}-1  \right ) \pi  \right \}.
	\end{equation}
\end{lem}
\begin{proof}
	Let $T$ be the orthogonal matrix given in \eqref{T}, then a direct calculation shows that  $2mn\times 2mn$  matrix $\mathcal{T}:=\begin{pmatrix}
		I_{m}\otimes T &0 \\
		0 &I_{m}\otimes T
	\end{pmatrix}$ is orthogonal symplectic and satisfies 
	\begin{equation}\label{TJT}
		\mathcal{T}^{-1}\begin{pmatrix}
			I_{m}\otimes A&0 \\
			0&I_{m}\otimes B
		\end{pmatrix}\mathcal{T}=\begin{pmatrix}
			I_{m}\otimes {\rm diag}(a_{1},\cdots, a_{n})&0 \\
			0&I_{m}\otimes {\rm diag}(b_{1},\cdots, b_{n})
		\end{pmatrix}.
	\end{equation}
	Since $Q_{n,m}=Q_{1,m}\otimes I_{n}$,  we also verify 
	\begin{equation}\label{TPT}
		\begin{split}
			\mathcal{T}^{-1}P_{n,m}\mathcal{T}=&\begin{pmatrix}
				I_{m}\otimes T^{-1} &0 \\
				0 &I_{m}\otimes T^{-1}
			\end{pmatrix}\begin{pmatrix}
				Q_{1,m}\otimes I_{n} &0 \\
				0& Q_{1,m}\otimes I_{n}
			\end{pmatrix}\begin{pmatrix}
				I_{m}\otimes T &0 \\
				0 &I_{m}\otimes T
			\end{pmatrix}\\
			=&\begin{pmatrix}
				Q_{1,m}\otimes I_{n} &0 \\
				0& Q_{1,m}\otimes I_{n}
			\end{pmatrix}=P_{n,m}.
		\end{split}
	\end{equation}
	Let $K_{n,m}$ be the $nm\times nm$  permutation matrix defined by 
	\begin{equation}\label{k}
		K_{n,m}e_{(i-1)m+j}=e_{(j-1)n+i},\quad 1\le i\le n,\ 1\le j\le m,
	\end{equation}
	where $e_{k}$ denotes the $k$-th standard basis vector of $\mathbb{R}^{nm} $. Since $K_{n,m}$ merely rearranges the standard basis vectors, it is an orthogonal matrix. Moreover, it is easily verified that 
	\begin{equation}\label{Ka}
		K_{n,m}^{-1}\left ( I_{m}\otimes {\rm diag}(a_{1},\cdots, a_{n}) \right )  K_{n,m}={\rm diag}(a_{1},\cdots, a_{n})\otimes I_{m}, 
	\end{equation}
	\begin{equation}\label{Kb}
		K_{n,m}^{-1}\left ( I_{m}\otimes {\rm diag}(b_{1},\cdots, b_{n}) \right )  K_{n,m}={\rm diag}(b_{1},\cdots, b_{n})\otimes I_{m},
	\end{equation}
	and 
	\begin{equation}\label{Kq}
		K_{n,m}^{-1}\left ( Q_{1,m}\otimes I_{n} \right )  K_{n,m}=I_{n}\otimes Q_{1,m}.
	\end{equation}
	
	Let $\mathcal{K}:=\begin{pmatrix}
		K_{n,m} & 0\\
		0 & K_{n,m}
	\end{pmatrix}$, then the $2nm\times 2nm$ matrix $\mathcal{P} :=\mathcal{T} \mathcal{K}$ is orthogonal symplectic and it follows from \eqref{TJT}, \eqref{Ka} and \eqref{Kb} that 
	\begin{equation}\label{PJP2}
		\begin{split}
			\mathcal{P}^{-1}J_{nm}\begin{pmatrix}
				I_{m}\otimes A&0 \\
				0&I_{m}\otimes B
			\end{pmatrix}\mathcal{P}=&J_{nm}\mathcal{P}^{-1}\begin{pmatrix}
				I_{m}\otimes A&0 \\
				0&I_{m}\otimes B
			\end{pmatrix}\mathcal{P}\\
			=&J_{nm}\mathcal{K}^{-1}\begin{pmatrix}
				I_{m}\otimes {\rm diag}(a_{1},\cdots, a_{n})&0 \\
				0&I_{m}\otimes {\rm diag}(b_{1},\cdots, b_{n})
			\end{pmatrix}\mathcal{K}\\
			=&J_{nm}\begin{pmatrix}
				{\rm diag}(a_{1},\cdots, a_{n})\otimes I_{m}&0 \\
				0& {\rm diag}(b_{1},\cdots, b_{n})\otimes I_{m}
			\end{pmatrix}.
		\end{split}
	\end{equation}
	Denote by $\hat{\gamma}  _{A,B}$ the fundamental matrix solution of 
	\begin{equation}
		\dot{u}(t)=J_{nm}\begin{pmatrix}
			{\rm diag}(a_{1},\cdots, a_{n})\otimes I_{m}&0 \\
			0& {\rm diag}(b_{1},\cdots, b_{n})\otimes I_{m}
		\end{pmatrix} u(t),
	\end{equation}
	with $\hat{\gamma}  _{A,B}(0)=I_{2nm}$, then from \eqref{PJP2} we deduce
	\begin{equation}\label{RAB}
		\mathcal{P}^{-1}\gamma _{A,B}(t)\mathcal{P}=\hat{\gamma}  _{A,B}(t).
	\end{equation}
	Moreover, denote $\hat{P}_{n,m}:=\mathcal{P}^{-1}P_{n,m}\mathcal{P}$, from \eqref{TPT} and \eqref{Kq}, we have 
	\begin{equation}\label{PAB}
		\hat{P}_{n,m}=\mathcal{P}^{-1}P_{n,m}\mathcal{P}=\mathcal{K}^{-1}\mathcal{T}^{-1}P_{n,m}\mathcal{T}\mathcal{K}=\mathcal{K}^{-1}P_{n,m}\mathcal{K}=\begin{pmatrix}
			I_{n}\otimes Q_{1,m}&0 \\
			0 &I_{n}\otimes Q_{1,m}
		\end{pmatrix},
	\end{equation}
	this together with Proposition \ref{Si} imply 
	\begin{equation}
		i_{\tau} ^{P_{n,m}}(\gamma _{A,B})=i_{\tau} ^{\hat{P} _{n,m}}(\hat{\gamma } _{A,B}),\quad \nu_{\tau} ^{P_{n,m}}(\gamma _{A,B})=\nu_{\tau} ^{\hat{P} _{n,m}}(\hat{\gamma } _{A,B}).
	\end{equation}
	
	For each $j=1,\cdots,n$, let $\hat{\gamma } _{A,B}^{j}$ be the fundamental matrix solution of 
	\begin{equation}
		\dot{v}(t)=J_{m}\begin{pmatrix}
			a_{j} I_{m}&0 \\
			0& b_{j}I_{m}
		\end{pmatrix} v(t),\quad v(t)\in \mathbb{R}^{2m} ,
	\end{equation}
	with $\hat{\gamma } _{A,B}^{j}(0)=I_{2m}$, that is, $$\hat{\gamma } _{A,B}^{j}(t)=\exp \left ( tJ_{m}\begin{pmatrix}
		a_{j} I_{m}&0 \\
		0& b_{j}I_{m}
	\end{pmatrix} \right ) ,$$
	 then $\hat{\gamma } _{A,B}^{j}$ coincides with $\gamma _{a_{j},b_{j}}$ defined in Lemma \ref{abv}. Observe that $J_{nm}=\left ( J_{m} \right ) ^{\Diamond n}$ and 
	\begin{equation}
		\begin{pmatrix}
			{\rm diag}(a_{1},\cdots, a_{n})\otimes I_{m}&0 \\
			0& {\rm diag}(b_{1},\cdots, b_{n})\otimes I_{m}
		\end{pmatrix}=\begin{pmatrix}
			a_{1}I_{m}&0 \\
			0 &b_{1}I_{m}
		\end{pmatrix}\Diamond \cdots \Diamond \begin{pmatrix}
			a_{n}I_{m}&0 \\
			0 &b_{n}I_{m}
		\end{pmatrix},
	\end{equation}
	by Lemma \ref{pDia}(i), we deduce 
	\begin{equation}
		J_{nm}\begin{pmatrix}
			{\rm diag}(a_{1},\cdots, a_{n})\otimes I_{m}&0 \\
			0& {\rm diag}(b_{1},\cdots, b_{n})\otimes I_{m}
		\end{pmatrix} =J_{m}\begin{pmatrix}
			a_{1}I_{m}&0 \\
			0 &b_{1}I_{m}
		\end{pmatrix}\Diamond \cdots \Diamond J_{m}\begin{pmatrix}
			a_{n}I_{m}&0 \\
			0 &b_{n}I_{m}
		\end{pmatrix}.
	\end{equation}
	By Lemma \ref{pDia}(iv), we further obtain  
	\begin{equation}
		\hat{\gamma}  _{A,B}(t)=\hat{\gamma } _{A,B}^{1}(t)\Diamond \cdots \Diamond \hat{\gamma } _{A,B}^{n}(t).
	\end{equation}
	Note that $\hat{P}_{n,m}=\begin{pmatrix}
		Q_{1,m} &0 \\
		0 &Q_{1,m}
	\end{pmatrix}^{\Diamond n}=\left ( 	P_{1,m} \right )^{\Diamond n} $, it follows from  Lemma \ref{pDia}(iii) that 
	\begin{equation}\label{hrP}
		\hat{\gamma}  _{A,B}(\tau )-\hat{P}_{n,m}=\left ( \hat{\gamma}  _{A,B}^{1}(\tau )-P_{1,m} \right )\Diamond \cdots \Diamond \left ( \hat{\gamma}  _{A,B}^{n}(\tau )-P_{1,m}  \right ) .
	\end{equation}
	\eqref{hrP}  together with Proposition \ref{Sa}, shows that  
	\begin{equation}\label{snu}
		\nu_{\tau} ^{P_{n,m}}(\gamma _{A,B})=\nu_{\tau} ^{\hat{P} _{n,m}}(\hat{\gamma } _{A,B})=\sum_{j=1}^{n} \nu_{\tau} ^{P_{1,m}}(\hat{\gamma } _{A,B}^{j})=\sum_{j=1}^{n} \nu_{\tau} ^{P_{1,m}}(\gamma _{a_{j},b_{j}}),
	\end{equation}
	and 
	\begin{equation}\label{si}
		i_{\tau} ^{P_{n,m}}(\gamma _{A,B})=i_{\tau} ^{\hat{P} _{n,m}}(\hat{\gamma } _{A,B})=\sum_{j=1}^{n} i_{\tau} ^{P_{1,m}}(\hat{\gamma } _{A,B}^{j})=\sum_{j=1}^{n} i_{\tau} ^{P_{1,m}}(\gamma _{a_{j},b_{j}}).
	\end{equation}
	Consequently, the desired equalities \eqref{nuP} and \eqref{iPr} follows from \eqref{snu}, \eqref{si} and Lemma \ref{abv}.
	
	Furthermore,  it follows from Lemma \ref{lpn} that 
	\begin{equation}
		i _{\tau }^{P_{n,m}}(\gamma _{A,B})= \sum_{0\le t<\tau }\nu_{t}^{P_{n,m}}(\gamma _{A,B}\mid _{\left [ 0,t \right ] })=\sum_{0< t<\tau }2\sharp \left \{ j\mid m\sqrt{a_{j}b_{j}} t \in  \left ( 2\mathbb{N}-1  \right ) \pi  \right \},
	\end{equation}
	if $A$ and $B$ are positive definite,  and 
	\begin{equation}
		i _{\tau }^{P_{n,m}}(\gamma _{A,B})= -\sum_{0< t\le \tau }\nu_{t}^{P_{n,m}}(\gamma _{A,B}\mid _{\left [ 0,t \right ] })=-\sum_{0< t\le \tau }2\sharp \left \{ j\mid m\sqrt{a_{j}b_{j}} t \in  \left ( 2\mathbb{N}-1  \right ) \pi  \right \},
	\end{equation}
	if $A$ and $B$ are negative definite, therefore we deduce that \eqref{iPAB1} and \eqref{iPAB2} hold. This completes the proof of Lemma \ref{abv2}.
\end{proof}

We now turn to the proofs of our main results, using the bifurcation theory presented by Lu \cite{Lu1}. Our main tool is the following theorem for autonomous Hamiltonian system with $P$-boundary condition, which we recall for the reader's convenience. 
\begin{thm}\label{thm33}
	For each $\lambda \in \Lambda $, let $u^{\lambda}$ be a constant solution of the Hamiltonian $P$-boundary value problem
	\begin{equation}\label{s33}
		\begin{cases}
			\dot{u} (t)=J\nabla H(\lambda , u(t)),\\
			u(t+\tau)=Pu(t),\quad t\in \mathbb{R}, 
		\end{cases}
	\end{equation}
	where the function $H(\lambda ,\cdot )$ is of class $C^{2}$, the gradient $\nabla H(\lambda ,u)$  with respect to $u$, as well as the corresponding Hessian matrix $\nabla^{2} H(\lambda ,u)$ depend continuously on $(\lambda ,u)\in \Lambda \times \mathbb{R} ^{2n}$.
	Let $\gamma_{\lambda}$ be the fundamental matrix solution of 
	\begin{equation}
		\dot{u} (t)=J\nabla^{2}H(\lambda, u^{\lambda})u(t),
	\end{equation}
	then $\gamma_{\lambda}(t)=\exp(tJ\nabla^{2}H(\lambda, u^{\lambda}))$, and we have 
		\begin{enumerate}[{\rm (I)}]
		\item {\rm (Necessary condition):} If $\left ( \mu ,u^{\mu} \right )$ with $\mu \in \Lambda $ is a bifurcation point of system  \eqref{s33} with respect to the branch $\left \{ (\lambda ,u^{\lambda })\mid \lambda \in \Lambda  \right \} $, then $\nu_{\tau}^{P}(\gamma_{\mu})\ne 0$. 
		\item {\rm (Sufficient condition):} Let $\Lambda$ be first countable. Suppose that for some $\mu \in \Lambda$, there exist two sequences $\left (\lambda _{k}^{-}\right ) $ and $\left (\lambda _{k}^{+}\right ) $ in $\Lambda$ converging to $\mu$ such that for each $k\in \mathbb{N} $, 
		\begin{equation}
			\left [ i_{\tau}^{P}(\gamma _{\lambda _{k}^{-}}),i_{\tau}^{P}(\gamma _{\lambda _{k}^{-}})+\nu_{\tau}^{P}(\gamma _{\lambda _{k}^{-}}) \right ]\cap \left [ i_{\tau}^{P}(\gamma _{\lambda _{k}^{+}}),i_{\tau}^{P}(\gamma _{\lambda _{k}^{+}})+\nu_{\tau}^{P}(\gamma _{\lambda _{k}^{+}}) \right ]=\emptyset  , 
		\end{equation} 
		and either $\nu_{\tau}^{P}(\gamma _{\lambda _{k}^{-}})=0$ or $\nu_{\tau}^{P}(\gamma _{\lambda _{k}^{+}})=0$. Then $(\mu,u^{\mu})$  is a bifurcation point of system  \eqref{s33} with respect to the branch $\left \{ (\lambda ,u^{\lambda })\mid \lambda \in \Lambda  \right \} $.
		\item {\rm (Alternative bifurcations of Fadell-Rabinowitz type and of Rabinowitz type):} For some interior point $\mu$ of $\Lambda$, assume that the following two conditions are satisfied,
		\begin{enumerate}[{\rm (a')}]
			\item	The following problem 
			\begin{equation}
					\begin{cases}
					\dot{u} (t)=J\nabla^{2} H(\mu, u^{\mu})u(t),\\
					u(t+\tau)=Pu(t),\quad t\in \mathbb{R}, 
				\end{cases}
			\end{equation}
			has no nonzero constant solutions, i.e., $\ker (P-I_{2n})\cap \ker (\nabla ^{2}H(\mu, u^{\mu}))=\left \{ 0 \right \} $.
			\item $\dim \ker (\gamma _{\mu}(\tau)-P)\ne 0$, $\dim \ker (\gamma _{\lambda}(\tau)-P)= 0$ for each $\lambda  \in \Lambda \setminus \left \{ \mu  \right \} $ near $\mu$, and $i_{\tau}^{P}(\gamma_{\lambda})$ takes, respectively, values $i_{\tau}^{P}(\gamma_{\mu})$ and $i_{\tau}^{P}(\gamma_{\mu})+\nu_{\tau}^{P}(\gamma_{\mu})$ as $\lambda \in \Lambda$ varies in two deleted half neighborhoods of $\mu$. 
		\end{enumerate}
		Then for system \eqref{s33}, there exist two possible alternatives:
		\begin{enumerate}[{\rm (i)}]
			\item The system \eqref{s33} with $\lambda=\mu$ admits a sequence of  $\mathbb{R} $-distinct solutions, $\bar{u}_{k}\notin \mathbb{R}\cdot u^{\mu},k=1,2,\cdots$, which converges to $u^{\mu}$ on any compact interval $I\subset \mathbb{R} $ in $C^{1}$-topology. 
			\item There exist left and right neighborhoods $\Lambda^{-}$ and $\Lambda^{+}$ of $\mu$ in $\Lambda$ and  integers $n^{+},n^{-}\ge 0$  satisfying $n^{+}+n^{-}\ge \nu_{\tau}^{P}(\gamma _{\mu})/2$, such that for $\lambda \in \Lambda ^{-}\setminus \left \{ \mu  \right \} (resp. \ \lambda \in \Lambda ^{+}\setminus \left \{ \mu  \right \})$, system \eqref{s33} with parameter $\lambda$ has at least $n^{-}(resp.\ n^{+})$ $\mathbb{R} $-distinct  solutions $u_{i}^{\lambda}\notin \mathbb{R}\cdot u^{\lambda} (i=1,\cdots,n^{-}(resp.\ n^{+}))$ which converge to  $u^{\lambda}$ on any compact interval $I\subset \mathbb{R} $ in $C^{1}$-topology as $\lambda\to \mu$.
			\item[]
			\hspace{-\labelwidth}\hspace{-\labelsep}
			Moreover, if $\nu_{\tau}^{P}(\gamma _{\mu})\ge 3$, then at least one of {\rm (i)} and the following {\rm (iii)} and {\rm (iv)} holds,
			\item For every $\lambda \in \Lambda \setminus \left \{ \mu  \right \} $ near $\mu$, there is a solution $\bar{u}^{\lambda}\notin \mathbb{R}\cdot u^{\lambda}$ of system \eqref{s33} with parameter $\lambda$, such that   $\bar{u}^{\lambda}$ converges to $u^{\lambda}$ on any compact interval $I\subset \mathbb{R} $ in $C^{1}$-topology as $\lambda\to \mu$.
			\item  For a given $\varepsilon >0$, there exists a one-sided neighborhood $\Lambda^{0}$ of $\mu$ in $\Lambda$ such that, for any $\lambda\in \Lambda ^{0}\setminus \left \{ \mu \right \} $, system \eqref{s33} with parameter $\lambda$ has either infinitely many  $\mathbb{R} $-distinct solutions $ \bar{u} _{k }^{\lambda} \notin \mathbb{R} \cdot u^{\lambda }$ such that $\left \| \bar{u} _{k}^{\lambda} \mid _{\left [ 0,\tau \right ] }-u^{\lambda }\mid _{\left [ 0,\tau \right ] }\right \| _{C^{1}}<\varepsilon,k=1,2,\cdots $, or at least two  $\mathbb{R} $-distinct solutions $\bar{u} ^{\lambda }\notin \mathbb{R} \cdot u^{\lambda } $ and $\hat{u} ^{\lambda }\notin \mathbb{R} \cdot u^{\lambda } $ such that $\left \| \bar{u} ^{\lambda} \mid _{\left [ 0,\tau \right ] }-u^{\lambda }\mid _{\left [ 0,\tau \right ] }\right \| _{C^{1}}<\varepsilon $, $\left \| \hat{u} ^{\lambda} \mid _{\left [ 0,\tau \right ] }-u^{\lambda }\mid _{\left [ 0,\tau \right ] }\right \| _{C^{1}}<\varepsilon $, and that
			\begin{equation}\label{cc}
			\int_{0}^{\tau } \left [ \frac{1}{2}\left ( J\dot{\bar{u} }^{\lambda }(t),\bar{u} ^{\lambda }(t)  \right )+H(\lambda ,\bar{u}^{\lambda } (t))   \right ] dt \ne \int_{0}^{\tau } \left [ \frac{1}{2}\left ( J\dot{\hat{u} }^{\lambda }(t),\hat{u} ^{\lambda }(t)  \right )+H(\lambda ,\hat{u}^{\lambda } (t))   \right ] dt. 
			\end{equation}						
		\end{enumerate}
	\end{enumerate}
\end{thm}

The first part of Theorem \ref{thm33}(III) is an alternative bifurcation result of Fadell-Rabinowitz's type, and the second part is that of Rabinowitz's type. We now proceed to the proofs.
\begin{proof}[\bf Proof of Theorem \ref{n1}] We divide the proof into the following three steps.

{\bf Step 1}(Proof of {\rm (I)}). Suppose that $\left ( \mu ,0 \right )$ with $\mu \in \Lambda $ is a bifurcation point of system \eqref{fx}, then as argued in Lemma \ref{l2}, $\left ( \mu ,0 \right )$ is also a bifurcation point of the first order Hamiltonian system 
	\begin{equation} \label{uef}
		\begin{cases}
			\dot{u}(t)=J_{m}\nabla  \tilde{H} (\lambda ,u(t)),\\
			u(t+\tau )=P_{1,m}u(t), \quad u(t)\in \mathbb{R}^{2m} ,
		\end{cases}
	\end{equation}
	where $	P_{1,m}:=\begin{pmatrix}
		Q_{1,m}&0 \\
		0&Q_{1,m}
	\end{pmatrix}$ with the $m\times m$ matrix $Q_{1,m}$ given in \eqref{418} with $n=1$.  
 For such $\mu \in \Lambda $, let $\gamma _{\mu }(t)$ be the fundamental solution of the following Hamiltonian system  linearized at the  trivial solution $0\in \mathbb{R}^{2m} $,
\begin{equation}
	\dot{u} (t)=J_{m}\nabla ^{2}\tilde{H} (\mu ,0) u(t), \quad \mu \in \Lambda ,\quad u(t)\in \mathbb{R}^{2m},
\end{equation}
with $\gamma _{\mu }(0)=I_{2m}$ and it follows from Theorem \ref{thm33}(I) that $\nu_{\tau}^{P_{1,m}}(\gamma _{\mu })\ne 0$. Specifically, the Hessian $\nabla ^{2}\tilde{H} (\mu ,0)$ has the block-diagonal form 
\begin{equation}
\nabla ^{2}\tilde{H} (\mu ,0)=\begin{pmatrix}
		f(\mu )I_{m}& 0\\
		0 &2\theta  (\mu )I_{m}
	\end{pmatrix},
\end{equation}
thus $\gamma _{\mu }$ coincides with $\gamma _{	f(\mu ),2\theta  (\mu )}$, and assertion {\rm (I)}  follows immediately from Lemma \ref{abv}.

{\bf Step 2}(Proof of {\rm (III)}). To prove {\rm (III)}, we next verify that  the first part of Theorem \ref{thm33}(III) can be applied to Hamiltonian system \eqref{uef}. Note that $P_{1,m}^{m}=-I_{2m}$ implies $\ker (P_{1,m}-I_{2m})=\left \{ 0 \right \} $, and thus we infer the condition (a') of Theorem \ref{thm33}(III) is satisfied. Now, it suffices to prove that  the condition (b') of Theorem \ref{thm33}(III) holds. From condition (a) in {\rm (II)}, there exists $j\in \mathbb{N} $ such that $m\sqrt{2f(\mu  )\theta (\mu  )}\tau=\left ( 2j-1 \right ) \pi $, then Lemma \ref{abv} gives 
\begin{equation}
\nu_{\tau}^{P_{1,m}}(\gamma _{\mu })=2, \quad i_{\tau}^{P_{1,m}}(\gamma _{\mu })=\left\{\begin{matrix}
	2(j-1), \quad {\rm if}\ f(\mu )>0,\\
	-2j,\quad \quad \ \  {\rm if} \ f(\mu )<0.
\end{matrix}\right.
\end{equation}
To compute the $P_{1,m}$-nullity $\nu_{\tau}^{P_{1,m}}(\gamma _{\lambda })$ and $P_{1,m}$-index $i_{\tau}^{P_{1,m}}(\gamma _{\lambda })$ for $\lambda \in \Lambda \setminus \left \{ \mu  \right \} $ near $\mu$, where $\gamma_{\lambda }$ agrees with $\gamma _{	f(\lambda ),2\theta  (\lambda )}$, we next consider $f(\mu )>0$ and $f(\mu )<0$ separately.

First, suppose that  $f(\mu )>0$, then since $f(\mu  )\theta (\mu  )>0$, we have $\theta (\mu  )>0$ as well.  Since both $f(\lambda)$ and $\theta (\lambda )$ are continuous, by further shrinking $\varepsilon >0$ if necessary, we may assume that $f(\lambda)>0$ and $\theta (\lambda )>0$ for all $\lambda \in \left ( \mu -\varepsilon ,\mu +\varepsilon  \right ) $, and that 
\begin{equation}\label{csb1}
	(2j-3)\pi <m\sqrt{2f(\lambda   )\theta (\lambda  )}\tau< \left ( 2j-1 \right ) \pi=\left ( 2\left ( j-1 \right )+1  \right ) \pi <m\sqrt{2f(\lambda'   )\theta (\lambda ' )}\tau<(2j+1)\pi ,
\end{equation}
for all $\left ( \lambda ,\lambda ' \right ) \in \left ( \mu -\varepsilon ,\mu  \right ) \times \left ( \mu ,\mu +\varepsilon  \right ) $  if the condition (sb.1) holds, and 
\begin{equation}\label{csb2}
(2j-3)\pi <m\sqrt{2f(\lambda'   )\theta (\lambda ' )}\tau< \left ( 2j-1 \right ) \pi=\left ( 2\left ( j-1 \right )+1  \right ) \pi <m\sqrt{2f(\lambda   )\theta (\lambda  )}\tau<(2j+1)\pi,
\end{equation}
for all $\left ( \lambda ,\lambda ' \right ) \in \left ( \mu -\varepsilon ,\mu  \right ) \times \left ( \mu ,\mu +\varepsilon  \right ) $  if the condition (sb.2) holds. Clearly, in both cases, we have 
\begin{equation}
	\nu_{\tau}^{P_{1,m}}(\gamma _{\lambda })=0,\quad 0<\left | \lambda -\mu  \right | <\varepsilon. 
\end{equation}
Moreover, if (sb.1) holds, it follows from \eqref{csb1} and Lemma \ref{abv} that  
\begin{equation}\label{1}
	i_{\tau}^{P_{1,m}}(\gamma _{\lambda })=\left\{\begin{matrix}
		2j-2=i_{\tau}^{P_{1,m}}(\gamma _{\mu  }),\quad \quad  \quad  \quad \  \ \forall \lambda \in \left ( \mu -\varepsilon ,\mu  \right ),  \\
		2j=i_{\tau}^{P_{1,m}}(\gamma _{\mu  })+\nu_{\tau}^{P_{1,m}}(\gamma _{\mu  }),\quad \forall \lambda \in \left ( \mu  ,\mu +\varepsilon \right ) .
	\end{matrix}\right.
\end{equation}
Similarly, if (sb.2) holds, \eqref{csb2} together with Lemma \ref{abv} yields
 \begin{equation}\label{2}
 	i_{\tau}^{P_{1,m}}(\gamma _{\lambda })=\left\{\begin{matrix}
 		2j=i_{\tau}^{P_{1,m}}(\gamma _{\mu  })+\nu_{\tau}^{P_{1,m}}(\gamma _{\mu  }),\quad \forall \lambda \in \left ( \mu -\varepsilon ,\mu  \right ),  \\
 		2j-2=i_{\tau}^{P_{1,m}}(\gamma _{\mu  }),\quad \quad  \quad  \quad \  \ \forall \lambda \in \left ( \mu  ,\mu +\varepsilon \right ) .
 	\end{matrix}\right.
 \end{equation}
 
 Now, suppose that $f(\mu )<0$, then $\theta (\mu  )<0$ as well. By further shrinking $\varepsilon >0$ if necessary, we may assume that $f(\lambda)<0$ and $\theta (\lambda )<0$ for all $\lambda \in \left ( \mu -\varepsilon ,\mu +\varepsilon  \right ) $, and that 
 \begin{equation}\label{ccsb1}
 	-(2j+1)\pi <-m\sqrt{2f(\lambda'   )\theta (\lambda ' )}\tau< \left ( 2\left ( -j \right )+1  \right ) \pi <-m\sqrt{2f(\lambda   )\theta (\lambda  )}\tau<\left ( 2(-j+1)+1 \right ) \pi ,
 \end{equation}
 for all $\left ( \lambda ,\lambda ' \right ) \in \left ( \mu -\varepsilon ,\mu  \right ) \times \left ( \mu ,\mu +\varepsilon  \right ) $  if the condition (sb.1) holds, and 
  \begin{equation}\label{ccsb2}
 	-(2j+1)\pi <-m\sqrt{2f(\lambda   )\theta (\lambda  )}\tau< \left ( 2\left ( -j \right )+1  \right ) \pi <-m\sqrt{2f(\lambda '  )\theta (\lambda ' )}\tau<\left ( 2(-j+1)+1 \right ) \pi ,
 \end{equation}
 for all $\left ( \lambda ,\lambda ' \right ) \in \left ( \mu -\varepsilon ,\mu  \right ) \times \left ( \mu ,\mu +\varepsilon  \right ) $  if the condition (sb.2) holds. Clearly, in both cases, we also obtain
 \begin{equation}
 	\nu_{\tau}^{P_{1,m}}(\gamma _{\lambda })=0,\quad 0<\left | \lambda -\mu  \right | <\varepsilon. 
 \end{equation}
 As argued above, if  (sb.1) holds, \eqref{ccsb1} and Lemma \ref{abv} imply
 \begin{equation}\label{3}
 	i_{\tau}^{P_{1,m}}(\gamma _{\lambda })=\left\{\begin{matrix}
 		-2j+2=i_{\tau}^{P_{1,m}}(\gamma _{\mu  })+\nu_{\tau}^{P_{1,m}}(\gamma _{\mu  }),\quad \forall \lambda \in \left ( \mu -\varepsilon ,\mu  \right ),  \\
 		-2j=i_{\tau}^{P_{1,m}}(\gamma _{\mu  }),\quad \quad  \quad  \quad \quad  \quad \quad   \quad \ \forall \lambda \in \left ( \mu  ,\mu +\varepsilon \right ) .
 	\end{matrix}\right.
 \end{equation}
 If  (sb.2) holds, then from \eqref{ccsb2} and  Lemma \ref{abv}, we deduce 
 \begin{equation}\label{4}
 		i_{\tau}^{P_{1,m}}(\gamma _{\lambda })=\left\{\begin{matrix}
 		-2j=i_{\tau}^{P_{1,m}}(\gamma _{\mu  }),\quad \quad  \quad  \quad \quad  \quad \quad  \quad \forall \lambda \in \left ( \mu -\varepsilon ,\mu  \right ),  \\
 		-2j+2=i_{\tau}^{P_{1,m}}(\gamma _{\mu  })+\nu_{\tau}^{P_{1,m}}(\gamma _{\mu  }),\quad  \forall \lambda \in \left ( \mu  ,\mu +\varepsilon \right ) .
 	\end{matrix}\right.
 \end{equation}
By \eqref{1}, \eqref{2}, \eqref{3} and \eqref{4}, condition (b') of Theorem \ref{thm33}(III) is satisfied. Applying the first part of that theorem  to Hamiltonian system \eqref{uef}, we obtain the following alternatives.
\begin{enumerate}[{\rm (i')}]
	\item Hamiltonian system \eqref{uef} with $\lambda=\mu$ admits a sequence of nonzero $\mathbb{R} $-distinct solutions $\left \{ u_{l} \right \} _{l}^{\infty }$ such that $u_{l}\to 0$ in $C_{{\rm loc}}^{1}\left ( \mathbb{R},\mathbb{R}^{2m}  \right ) $ as $l\to \infty $.
	\item There exist left and right neighborhoods $\Lambda^{-}$ and $\Lambda^{+}$ of $\mu$ in $\Lambda$ and nonnegative integers $n^{+}$ and $n^{-}$ satisfying $n^{+}+n^{-}\ge \nu_{\tau}^{P_{1,m}}(\gamma _{\mu })/2=1$, such that for $\lambda \in \Lambda ^{-}\setminus \left \{ \mu  \right \} ({\rm resp. } \ \lambda \in \Lambda ^{+}\setminus \left \{ \mu  \right \})$, system \eqref{uef}  with parameter $\lambda$ has at least $n^{-}({\rm resp. } \ n^{+})$ $\mathbb{R} $-distinct nonzero solutions $u^{\lambda}_{i}(i=1,\cdots,n^{-}({\rm resp. } \ n^{+}))$ which converge to zero in $C_{{\rm loc}}^{1}\left ( \mathbb{R},\mathbb{R}^{2m}  \right ) $ as $\lambda\to \mu$.
\end{enumerate}
By the discussion following Lemma \ref{l2}, these two alternatives are equivalent to assertions (i) and (ii) of Theorem \ref{n1}(III), respectively.

{\bf Step 3}(Proof of {\rm (II)}). Let $\varepsilon $ be given above, we may assume $\left ( \lambda _{q}^{\pm } \right )_{q}\subset \left ( \mu -\varepsilon ,\mu +\varepsilon  \right )\setminus \left \{ \mu  \right \} $, then it holds either 
\begin{equation}\label{cccsb1}
	\left ( \lambda _{q}^{-} \right )_{q}\subset \left ( \mu -\varepsilon ,\mu  \right )\ {\rm and}\ \left ( \lambda _{q}^{+} \right )_{q}\subset \left ( \mu  ,\mu+\varepsilon  \right ),
\end{equation}
or 
\begin{equation}\label{cccsb2}
\left ( \lambda _{q}^{+} \right )_{q}\subset \left ( \mu -\varepsilon ,\mu  \right )\ {\rm and}\ \left ( \lambda _{q}^{-} \right )_{q}\subset \left ( \mu  ,\mu+\varepsilon  \right ).
\end{equation}
Therefore, as argued in the proof of (III), from condition (b) in (II),  we have for any $q\in \mathbb{N} $,
\begin{equation}
		\nu_{\tau}^{P_{1,m}}(\gamma _{\lambda_{q}^{+} })=\nu_{\tau}^{P_{1,m}}(\gamma _{\lambda_{q}^{-} })=0.
\end{equation}
Moreover, suppose that $f(\mu )>0$, then we deduce for any $q\in \mathbb{N} $,
\begin{equation}
	i_{\tau}^{P_{1,m}}(\gamma _{\lambda })=\left\{\begin{matrix}
		2j-2,  \quad  \ \lambda =\lambda _{q}^{-},  \\
		2j,\quad \quad  \quad \lambda =\lambda _{q}^{+}.
	\end{matrix}\right.
\end{equation}
if \eqref{cccsb1} holds, and 
\begin{equation}
	i_{\tau}^{P_{1,m}}(\gamma _{\lambda })=\left\{\begin{matrix}
2j,\quad \quad  \quad \lambda =\lambda _{q}^{+}	,  \\
			2j-2,  \quad  \ \lambda =\lambda _{q}^{-}.
	\end{matrix}\right.
\end{equation}
if \eqref{cccsb2} holds. Consequently, for any $q\in \mathbb{N}$, we obtain
\begin{equation}\label{IIp}
	\left [ i_{\tau}^{P_{1,m}}(\gamma _{\lambda_{q}^{-} }),i_{\tau}^{P_{1,m}}(\gamma _{\lambda_{q}^{-} })+\nu_{\tau}^{P_{1,m}}(\gamma _{\lambda_{q}^{-} }) \right ]\cap \left [ i_{\tau}^{P_{1,m}}(\gamma _{\lambda_{q}^{+} }),i_{\tau}^{P_{1,m}}(\gamma _{\lambda_{q}^{+} })+\nu_{\tau}^{P_{1,m}}(\gamma _{\lambda_{q}^{+} }) \right ] =\emptyset .
\end{equation}
A similar argument shows that \eqref{IIp} also holds when $f(\mu )<0$, thus it follows from Theorem \ref{thm33}(II) that $(\mu,0)$ is a bifurcation point of system \eqref{uef}, which implies that $(\mu,0)$ is a a bifurcation point of system \eqref{fx}.
\end{proof}

\begin{proof}[\bf Proof of Theorem \ref{thm2}] The proof is divided into three steps.

	{\bf Step 1}(Proof of {\rm (I)}). If $\left ( \mu ,0 \right )$ with $\mu \in \Lambda $ is a bifurcation point of system  \eqref{e2}, then as shown in the discussion following Lemma \ref{l2}, $\left ( \mu ,0 \right )$ is also a bifurcation point of Hamiltonian system 
\begin{equation}\label{ut}
	\begin{cases}
		\dot{u}(t)=J_{nm}\nabla \tilde{H} (\lambda ,u(t)),\\
		u(t+\tau )=P_{n,m}u(t),\quad u(t)\in \mathbb{R}^{2nm} .
	\end{cases}
\end{equation}
Let $\gamma _{\mu }(t)$ be the fundamental solution of linear Hamiltonian system 
\begin{equation}
	\dot{u} (t)=J_{nm}\nabla ^{2}\tilde{H} (\mu ,0) u(t), \quad \mu \in \Lambda ,\quad u(t)\in \mathbb{R}^{2nm},
\end{equation}
with $\gamma _{\mu }(0)=I_{2nm}$, it then follows from Theorem \ref{thm33}(I) that $\nu_{\tau}^{P_{n,m}}(\gamma _{\mu })\ne 0$. Specifically, the Hessian $\nabla ^{2}\tilde{H} (\mu ,0)$ has the following block-diagonal form 
\begin{equation}
\nabla ^{2}\tilde{H} (\mu  ,0)=\begin{pmatrix}
	\nabla ^{2}H(\mu,0)& 0\\
	0 &-\mathcal{A} (\mu )
\end{pmatrix}=\begin{pmatrix}
	I_{m}\otimes \nabla ^{2}F(\mu ,0) & 0\\
	0 &I_{m}\otimes 2\Theta (\mu  )
\end{pmatrix},
\end{equation}
here $\otimes$ denotes the Kronecker product, thus $\gamma _{\mu }$ coincides with $\gamma _{	\nabla ^{2}F(\mu ,0),2\Theta (\mu  )}$ and  Lemma \ref{abv2} yields
\begin{equation}\label{nt}
\nu _{\tau}^{P_{n,m}}(\gamma _{\mu })=2\sharp \left \{ j\mid m\sqrt{2f_{j}(\mu )\theta _{j}(\mu )} \tau \in  \left ( 2\mathbb{N}-1  \right ) \pi  \right \}, 
\end{equation}
which implies that (I) holds.

{\bf Step 2}(Proof of {\rm (III)}). By Lemma \ref{l1} and Lemma \ref{l2}, we only need to verify that Theorem \ref{thm33}(III) can be applied to first order Hamiltonian system \eqref{ut}. Since $P_{n,m}^{m}=-I_{2nm}$, it holds $\ker \left ( P_{n,m}-I_{2nm} \right )=\left \{ 0 \right \}  $. Therefore,  condition (a') of Theorem \ref{thm33}(III) is satisfied. It remains to verify that condition (b') of Theorem \ref{thm33}(III) holds.

Let $\gamma _{f_{j}(\mu),2\theta _{j}(\mu)}$ be the path defined in Lemma \ref{abv}, and let $\gamma _{\mu }$ be the path introduced above. By condition (B) in {\rm (II)}, equation \eqref{nt} yields
\begin{equation}\label{nP}
	\nu _{\tau}^{P_{n,m}}(\gamma _{\mu })=2k.
\end{equation}
We next compute the index $i _{\tau}^{P_{n,m}}(\gamma _{\mu })$. On one hand, for each $l\in \left \{ 1,\cdots ,k \right \} $, we have
\begin{equation}
		i _{\tau}^{P_{1,m}}(\gamma _{f_{j_{l}}(\mu),2\theta _{j_{l}}(\mu)})=\left\{\begin{matrix}
		2\left ( \rho _{l}(\mu)-1 \right ), \quad {\rm if}\ f_{j_{l}}(\mu)>0, \\
		-2\rho _{l}(\mu),\quad \quad \quad {\rm if}\ f_{j_{l}}(\mu)<0.
	\end{matrix}\right.
\end{equation}
On the other hand, let $\Gamma \subset \left \{ k+1,\cdots ,n \right \} $ denote the set of indices $l$ for which the inequality \eqref{2r} holds. We decompose $\Gamma$ as
\begin{equation}
\Gamma =\Gamma ^{+}\cup \Gamma ^{-},
\end{equation}
 where 
 \begin{equation}
 	\Gamma ^{+}:=\left \{ l\in\Gamma \mid f_{j_{l}}(\mu)>0 \right \},\quad \Gamma ^{-}:=\left \{ l\in\Gamma \mid f_{j_{l}}(\mu)<0 \right \}.
 \end{equation}
 Then it follows from condition (A) and Lemma \ref{abv} that for each $l\in\left \{ k+1,\cdots ,n \right \}$,
 \begin{equation}\label{iP1}
 	i _{\tau}^{P_{1,m}}(\gamma _{f_{j_{l}}(\mu),2\theta _{j_{l}}(\mu)})=\left\{\begin{matrix}
 		0,\quad \quad \quad \  {\rm if}\ l\notin  \Gamma , \\
 		2\rho _{l}(\mu), \quad \ \   {\rm if}\ l\in \Gamma ^{+},\\
 		-2\rho _{l}(\mu),\quad {\rm if}\ l\in \Gamma ^{-}.
 	\end{matrix}\right.
 \end{equation}
 Consequently, by Lemma \ref{abv2}, we deduce
 \begin{equation}\label{iPru}
 	\begin{split}
 		i _{\tau}^{P_{n,m}}(\gamma _{\mu })=&\sum_{j=1}^{n} i _{\tau}^{P_{1,m}}(\gamma _{f_{j}(\mu),2\theta _{j}(\mu)})\\
 		=&\sum_{l\le k,f_{j_{l}}(\mu)>0}2\left ( \rho _{l}(\mu)-1 \right )+\sum_{l\le k,f_{j_{l}}(\mu)<0}	\left ( -2\rho _{l}(\mu)  \right )+\sum_{l\in\Gamma ^{+}}2\rho _{l}(\mu)+\sum_{l\in\Gamma ^{-}}\left ( -2\rho _{l}(\mu) \right ).  
 	\end{split}
 \end{equation}

 Let $\gamma _{\lambda }:=\gamma _{	\nabla ^{2}F(\lambda  ,0),2\Theta (\lambda  )}$ be as in Lemma \ref{abv2},  we next analyze the behavior of $i _{\tau}^{P_{n,m}}(\gamma _{\lambda })$ and $\nu _{\tau}^{P_{n,m}}(\gamma _{\lambda })$ for $\lambda$ in a neighborhood of $\mu$. Since $f_{j}(\lambda )\theta _{j}(\lambda ), j=1,\cdots ,n$, depends continuously on $\lambda$, we may,  after further shrinking $\varepsilon >0$ if necessary, rewrite conditions  (A) and (SC) as follows,
 	\begin{enumerate}
 	\item[{\rm (A')}] For each $l=k+1,\cdots,n$, the following hold,
 	\begin{equation}
 		2f_{j_{l}}(\lambda)\theta _{j_{l}}(\lambda)< \left ( \pi /m\tau  \right ) ^{2},\quad {\rm for\ all}\ \lambda \in \left ( \mu -\varepsilon ,\mu +\varepsilon  \right ),\  {\rm if}\ l\notin \Gamma,  
 	\end{equation}
 	\begin{equation}
 		\left ( 2\rho _{l}(\mu)-1 \right ) \pi <m\sqrt{2f_{j_{l}}(\lambda)\theta _{j_{l}}(\lambda)} \tau < \left ( 2\rho _{l}(\mu)+1 \right ) \pi , \quad {\rm for\ all}\ \lambda \in \left ( \mu -\varepsilon ,\mu +\varepsilon  \right ),\  {\rm if}\ l\in \Gamma.
 	\end{equation}
\item[{\rm (SC')}] For some $\varepsilon>0$, one of the following two conditions holds.
 		\begin{enumerate}
 			\item[{\rm (SC'.1)}] For any $l=1,\cdots,k$, if $ f_{j_{l}}(\mu  )>0$, then for any $\left ( \lambda ,\lambda ' \right )\in \left ( \mu-\varepsilon ,\mu  \right )\times \left ( \mu ,\mu +\varepsilon \right ) $, we have $ f_{j_{l}}(\lambda  )>0$ and $ f_{j_{l}}(\lambda'  )>0$, and
 			\begin{equation}
 				\left ( 2\rho _{l}(\mu)-3 \right )\pi  <m\sqrt{2f_{j_{l}}(\lambda )\theta _{j_{l}}(\lambda )}\tau  < \left ( 2\rho _{l}(\mu)-1 \right ) \pi < m\sqrt{2f_{j_{l}}(\lambda'  )\theta_{j_{l}} (\lambda' )}\tau<\left ( 2 \rho _{l}(\mu)+1   \right )\pi;  
 			\end{equation} 
 			and if $ f_{j_{l}}(\mu  )<0$, then for any $\left ( \lambda ,\lambda ' \right )\in \left ( \mu-\varepsilon ,\mu  \right )\times \left ( \mu ,\mu +\varepsilon \right ) $, we have $ f_{j_{l}}(\lambda  )<0$ and $ f_{j_{l}}(\lambda'  )<0$, and 
 			\begin{equation}
 				\left ( 2 \rho _{l}(\mu)+1   \right )\pi>m\sqrt{2f_{j_{l}}(\lambda )\theta _{j_{l}}(\lambda )}\tau  > \left ( 2\rho _{l}(\mu)-1 \right ) \pi > m\sqrt{2f_{j_{l}}(\lambda'  )\theta_{j_{l}} (\lambda' )}\tau>\left ( 2\rho _{l}(\mu)-3 \right )\pi.
 			\end{equation}
 			\item[{\rm (SC'.2)}] For any $l=1,\cdots,k$, if $ f_{j_{l}}(\mu  )>0$, then for any $\left ( \lambda ,\lambda ' \right )\in \left ( \mu-\varepsilon ,\mu  \right )\times \left ( \mu ,\mu +\varepsilon \right ) $, we have $ f_{j_{l}}(\lambda  )>0$ and $ f_{j_{l}}(\lambda'  )>0$, and
 			\begin{equation}
 				\left ( 2 \rho _{l}(\mu)+1   \right )\pi>m\sqrt{2f_{j_{l}}(\lambda )\theta _{j_{l}}(\lambda )}\tau  > \left ( 2\rho _{l}(\mu)-1 \right ) \pi > m\sqrt{2f_{j_{l}}(\lambda'  )\theta_{j_{l}} (\lambda' )}\tau>\left ( 2\rho _{l}(\mu)-3 \right )\pi;
 			\end{equation}
 			and if $ f_{j_{l}}(\mu  )<0$, then for any $\left ( \lambda ,\lambda ' \right )\in \left ( \mu-\varepsilon ,\mu  \right )\times \left ( \mu ,\mu +\varepsilon \right ) $, we have $ f_{j_{l}}(\lambda  )<0$ and $ f_{j_{l}}(\lambda'  )<0$, and 
 			\begin{equation}
 				\left ( 2\rho _{l}(\mu)-3 \right )\pi  <m\sqrt{2f_{j_{l}}(\lambda )\theta _{j_{l}}(\lambda )}\tau  < \left ( 2\rho _{l}(\mu)-1 \right ) \pi < m\sqrt{2f_{j_{l}}(\lambda'  )\theta_{j_{l}} (\lambda' )}\tau<\left ( 2 \rho _{l}(\mu)+1   \right )\pi.  
 			\end{equation} 
 		\end{enumerate}
 		\end{enumerate}
 	Combined condition (A') and (SC'), one can easily check 
 	\begin{equation}\label{n0}
 		\nu _{\tau}^{P_{n,m}}(\gamma _{\lambda })=0,\quad 	\forall \lambda \in \left ( \mu-\varepsilon ,\mu +\varepsilon  \right )\setminus \left \{ \mu  \right \},
 	\end{equation}
 	and for any $\left ( \lambda ,l \right )\in \left ( \mu -\varepsilon ,\mu +\varepsilon  \right ) \times \left \{ k+1,\cdots ,n \right \}  $,
 	\begin{equation}\label{iP2}
 		i _{\tau}^{P_{1,m}}(\gamma _{f_{j_{l}}(\lambda),2\theta _{j_{l}}(\lambda)})=\left\{\begin{matrix}
 			0,\quad \quad \quad \  {\rm if}\ l\notin  \Gamma , \\
 			2\rho _{l}(\mu), \quad \ \   {\rm if}\ l\in \Gamma ^{+},\\
 			-2\rho _{l}(\mu),\quad {\rm if}\ l\in \Gamma ^{-}.
 		\end{matrix}\right.
 	\end{equation}
 	
 	Suppose now (SC'.1) holds. For any $\left ( \lambda ,l \right )\in \left ( \mu -\varepsilon ,\mu   \right ) \times \left \{ 1,\cdots ,k \right \}  $, by  Lemma \ref{abv},  we can deduce from (SC'.1) that 
 	\begin{equation}\label{iP3}
 		i _{\tau}^{P_{1,m}}(\gamma _{f_{j_{l}}(\lambda),2\theta _{j_{l}}(\lambda)})=\left\{\begin{matrix}
 			2\left ( \rho _{l}(\mu)-1 \right ), \quad   {\rm if}\ f_{j_{l}}(\mu)>0, \\
 			-2\rho _{l}(\mu) ,\quad \quad \quad {\rm if}\ f_{j_{l}}(\mu)<0.
 		\end{matrix}\right.
 	\end{equation}
 	Therefore, for any $\lambda\in \left ( \mu -\varepsilon ,\mu   \right )$, it follows from \eqref{iPru}, \eqref{iP2}, \eqref{iP3} and Lemma \ref{abv2} that
 	\begin{equation}\label{ip}
 		\begin{split}
 			i _{\tau}^{P_{n,m}}(\gamma _{\lambda })=&\sum_{j=1}^{n} i _{\tau}^{P_{1,m}}(\gamma _{f_{j}(\lambda ),2\theta _{j}(\lambda )})\\
 			=&\sum_{l\le k,f_{j_{l}}(\mu)>0}2\left ( \rho _{l}(\mu)-1 \right )+\sum_{l\le k,f_{j_{l}}(\mu)<0}(	-2\rho _{l}(\mu))+\sum_{l\in\Gamma ^{+}}2\rho _{l}(\mu)+\sum_{l\in\Gamma ^{-}}\left ( -2\rho _{l}(\mu) \right )\\
 			=&i _{\tau}^{P_{n,m}}(\gamma _{\mu }).
 		\end{split}
 	\end{equation}
 	We next turn to the case  $\left ( \lambda' ,l \right )\in \left ( \mu ,\mu +\varepsilon    \right ) \times \left \{ 1,\cdots ,k \right \}  $. By  Lemma \ref{abv},  condition (SC'.1) yields that
 		\begin{equation}\label{iP4}
 			i _{\tau}^{P_{1,m}}(\gamma _{f_{j_{l}}(\lambda'),2\theta _{j_{l}}(\lambda')})=\left\{\begin{matrix}
 				2\rho _{l}(\mu), \quad \quad \quad \quad \quad \ \   {\rm if}\ f_{j_{l}}(\mu)>0, \\
 				2\left ( 1-\rho _{l}(\mu) \right )  ,\quad \quad \quad {\rm if}\ f_{j_{l}}(\mu)<0.
 			\end{matrix}\right.
 		\end{equation}
 		Hence, for every $\lambda' \in \left ( \mu ,\mu +\varepsilon    \right )$, from \eqref{nP}, \eqref{iPru}, \eqref{iP2}, \eqref{iP4} and Lemma \ref{abv2}, we deduce 
 	\begin{equation}
 		\begin{split}\label{ip1}
 			i _{\tau}^{P_{n,m}}(\gamma _{\lambda' })=&\sum_{j=1}^{n} i _{\tau}^{P_{1,m}}(\gamma _{f_{j}(\lambda '),2\theta _{j}(\lambda' )})\\
 			=&\sum_{l\le k,f_{j_{l}}(\mu)>0}2\rho _{l}(\mu)+\sum_{l\le k,f_{j_{l}}(\mu)<0}2(	1-\rho _{l}(\mu))+\sum_{l\in\Gamma ^{+}}2\rho _{l}(\mu)+\sum_{l\in\Gamma ^{-}}\left ( -2\rho _{l}(\mu) \right )\\
 			=&\sum_{l\le k,f_{j_{l}}(\mu)>0}2(\rho _{l}(\mu)-1)+\sum_{l\le k,f_{j_{l}}(\mu)<0}(-2\rho _{l}(\mu))+\sum_{l\in\Gamma ^{+}}2\rho _{l}(\mu)+\sum_{l\in\Gamma ^{-}}\left ( -2\rho _{l}(\mu) \right )+2k\\
 			=&i _{\tau}^{P_{n,m}}(\gamma _{\mu })+\nu _{\tau}^{P_{n,m}}(\gamma _{\mu }).
 		\end{split}
 	\end{equation}
 	
 	Similarly, if condition (SC'.2) holds, the same argument yields
 	\begin{equation}\label{ip2}
 		i _{\tau}^{P_{n,m}}(\gamma _{\lambda })=\left\{\begin{matrix}
 			i _{\tau}^{P_{n,m}}(\gamma _{\mu })+\nu _{\tau}^{P_{n,m}}(\gamma _{\mu }),\quad  \forall \lambda \in \left ( \mu -\varepsilon ,\mu   \right ),\\
 			i _{\tau}^{P_{n,m}}(\gamma _{\mu }), \quad \quad \quad \quad\quad\quad \ \ \forall \lambda \in \left ( \mu ,\mu +\varepsilon    \right ).
 		\end{matrix}\right.
 	\end{equation}
 	Then it follows from \eqref{n0}, \eqref{ip}, \eqref{ip1} and \eqref{ip2} that condition (b') of Theorem \ref{thm33}(III) is satisfied. Therefore, applying the first part of Theorem \ref{thm33}(III) to Hamiltonian system \eqref{ut}, we obtain the following alternatives:
 		\begin{enumerate}[{\rm (i')}]
 		\item The system \eqref{ut} with $\lambda=\mu$ admits a sequence of nonzero $\mathbb{R} $-distinct solutions $\left \{ u_{l} \right \} _{l}^{\infty }$ such that $u_{l}\to 0$ in $C_{{\rm loc}}^{1}\left ( \mathbb{R},\mathbb{R}^{2nm}  \right ) $ as $l\to \infty $.
 		\item There exist left and right neighborhoods $\Lambda^{-}$ and $\Lambda^{+}$ of $\mu$ in $\Lambda$ and nonnegative integers $n^{+}$ and $n^{-}$ satisfying $n^{+}+n^{-}\ge \nu _{\tau}^{P_{n,m}}(\gamma _{\mu })/2= k$, such that for $\lambda \in \Lambda ^{-}\setminus \left \{ \mu  \right \} ({\rm resp. } \ \lambda \in \Lambda ^{+}\setminus \left \{ \mu  \right \})$, system \eqref{ut} with parameter $\lambda$ has at least $n^{-}({\rm resp. }\ n^{+})$ $\mathbb{R} $-distinct nonzero solutions $u_{i}^{\lambda}(i=1,\cdots,n^{-}({\rm resp. }\ n^{+}))$ which converge to zero in $C_{{\rm loc}}^{1}\left ( \mathbb{R},\mathbb{R}^{2nm}  \right ) $ as $\lambda\to \mu$.
	   \item[]
\hspace*{-\leftmargin}%
\parbox{\dimexpr\linewidth+\leftmargin\relax}{%
 	Moreover, if $k>1$, then $\nu _{\tau}^{P_{n,m}}(\gamma _{\mu })=2k>3$, and hence we may apply the second part of Theorem \ref{thm33}(III) to system \eqref{ut} and conclude that at least one of (i') and the following (iii') and (iv') holds,}
 		\item For every $\lambda \in \Lambda \setminus \left \{ \mu  \right \} $ near $\mu$, there is a nonzero solution $u^{\lambda}$ of system \eqref{ut} with parameter $\lambda$, such that   $u^{\lambda}\to 0$ in $C_{{\rm loc}}^{1}\left ( \mathbb{R},\mathbb{R}^{2nm}  \right ) $ as $\lambda\to \mu$.
 	\item  For a given $\varepsilon >0$, there exists a one-sided neighborhood $\Lambda^{0}$ of $\mu$ in $\Lambda$ such that, for any $\lambda\in \Lambda ^{0}\setminus \left \{ \mu \right \} $, system \eqref{ut} with parameter $\lambda$ has either infinitely many nonzero $\mathbb{R} $-distinct solutions $\left \{ u_{l }^{\lambda} \right \}_{l=1}^{\infty } $ such that $\left \| u_{l}^{\lambda} \mid _{\left [ 0,\tau \right ] }\right \| _{C^{1}}<\varepsilon $ for all $l\in \mathbb{N} $, or at least two nonzero $\mathbb{R} $-distinct solutions $\bar{u} ^{\lambda } $ and $\hat{u} ^{\lambda } $ such that $\left \| \bar{u}^{\lambda } \mid _{\left [ 0,\tau \right ] }\right \| _{C^{1}}<\varepsilon $, $\left \| \hat{u}^{\lambda } \mid _{\left [ 0,\tau \right ] }\right \| _{C^{1}}<\varepsilon $, and 
 	\begin{equation}\label{ne1}
 	\int_{0}^{\tau } \left [ \frac{1}{2}\left ( J_{nm}\dot{\bar{u} }^{\lambda }(t),\bar{u} ^{\lambda }(t)  \right )+\tilde{H}(\lambda ,\bar{u}^{\lambda } (t))   \right ] dt \ne \int_{0}^{\tau } \left [ \frac{1}{2}\left ( J_{nm}\dot{\hat{u} }^{\lambda }(t),\hat{u} ^{\lambda }(t)  \right )+\tilde{H}(\lambda ,\hat{u}^{\lambda } (t))   \right ] dt ,
 	\end{equation}
 	where $\bar{u}^{\lambda }(t) =\left ( \bar{z}^{\lambda }(t), \bar{	w}^{\lambda } (t)\right ) ^{T}$, $\bar{	w}^{\lambda } (t)=\mathcal{A}(\lambda )^{-1}\dot{\bar{z} }^{\lambda }(t)$ and $\hat{u}^{\lambda }(t) =\left ( \hat{z}^{\lambda }(t), \hat{	w}^{\lambda } (t)\right ) ^{T}$, $\hat{	w}^{\lambda } (t)=\mathcal{A}(\lambda )^{-1}\dot{\hat{z} }^{\lambda }(t)$.						
 	\end{enumerate}
 	By the correspondence induced by Lemmas \ref{l1} and \ref{l2},  the preceding bifurcation alternatives for system \eqref{ut} yield the corresponding alternatives for system  \eqref{e2}. In addition, from \eqref{ne}, inequality \eqref{ne1} implies 
 	\begin{equation}\label{int0}
 		\int_{0}^{\tau} \left [ \frac{1}{2}\left ( \mathcal{A}(\lambda )^{-1}\dot{\bar{z} }^{\lambda } (t), \dot{\bar{z} }^{\lambda } (t)  \right ) +H(\lambda ,\bar{z}^{\lambda } (t)) \right ] dt\ne \int_{0}^{\tau} \left [ \frac{1}{2}\left ( \mathcal{A}(\lambda )^{-1}\dot{\hat{z} }^{\lambda } (t), \dot{\hat{z} }^{\lambda } (t)  \right ) +H(\lambda ,\hat{z}^{\lambda } (t)) \right ] dt,
 	\end{equation}
 	where $\bar{z}^{\lambda } (t)=\left ( \bar{z} _{1}^{\lambda }(t),\cdots , \bar{z} _{m}^{\lambda }(t)\right )^{T}, \bar{z} _{j}^{\lambda }\left ( t \right )= \int_{\left ( j-1 \right )\tau  }^{j\tau }\Theta (\lambda )\bar{x}^{\lambda } \left ( t-s \right )  ds $ and $\hat{z}^{\lambda } (t)=\left ( \hat{z} _{1}^{\lambda }(t),\cdots , \hat{z} _{m}^{\lambda }(t)\right )^{T}, \hat{z} _{j}^{\lambda }\left ( t \right )= \int_{\left ( j-1 \right )\tau  }^{j\tau }\Theta (\lambda )\hat{x}^{\lambda } \left ( t-s \right )  ds $. Since $\dot{\bar{z}}  _{j}^{\lambda }\left ( t \right ) =\Theta (\lambda ) \left ( \bar{x}^{\lambda } \left ( t-(j-1)\tau  \right )- \bar{x}^{\lambda } \left ( t-j\tau  \right )  \right ) $ and $\bar{x}^{\lambda } \left ( t-m\tau  \right )=-\bar{x}^{\lambda } (t)$, we deduce 
 	\begin{equation}\label{fra}
 		\begin{split}
 			\frac{1}{2}\int_{0}^{\tau}\left ( \mathcal{A}(\lambda )^{-1}\dot{\bar{z} }^{\lambda } (t), \dot{\bar{z} }^{\lambda } (t)  \right )  dt&=-\frac{1}{4} \sum_{j=1}^{m}\int_{0}^{\tau }\left ( \Theta (\lambda )^{-1}\dot{\bar{z} }_{j}^{\lambda }(t), \dot{\bar{z} }_{j}^{\lambda }(t) \right )dt\\
 			&=-\frac{1}{4}\int_{0}^{m\tau } \left ( \Theta (\lambda )\left ( \bar{x}^{\lambda } (t)-\bar{x}^{\lambda } (t-\tau ) \right ), \bar{x}^{\lambda } (t)-\bar{x}^{\lambda } (t-\tau ) \right )dt. 
 		\end{split}
 	\end{equation}
 	Moreover, since the function $F(\lambda ,\cdot )$ is even, we have
 	\begin{equation}\label{int}
 		\int_{0}^{\tau } H(\lambda ,\bar{z}^{\lambda } (t))dt=\sum_{j=1}^{m}\int_{0}^{\tau }F(\lambda ,\bar{z} _{j}^{\lambda }(t) )dt=\int_{0}^{m\tau }F\left ( \lambda ,\int_{0}^{\tau }\Theta (\lambda ) \bar{x} ^{\lambda }\left ( t-s \right ) ds \right ) dt .
 	\end{equation}
 	Combining \eqref{int0}, \eqref{fra} and \eqref{int}, we obtain  \eqref{c},  which completes the proof of Theorem \ref{thm2}(III).
 	
 	{\bf Step 3}(Proof of {\rm (II)}). Suppose that condition (C) holds. Choose $\varepsilon>0$ as in {\bf Step 2}, and take the sequences $\lambda_q^\pm$ sufficiently close to $\mu$ so that they lie in $(\mu-\varepsilon,\mu+\varepsilon)\setminus{\mu}$. By passing to a subsequence if necessary, we may assume that  either 
 	\begin{equation}\label{cccsb11}
 		\left ( \lambda _{q}^{-} \right )_{q}\subset \left ( \mu -\varepsilon ,\mu  \right )\ {\rm and}\ \left ( \lambda _{q}^{+} \right )_{q}\subset \left ( \mu  ,\mu+\varepsilon  \right ),
 	\end{equation}
 	or 
 	\begin{equation}\label{cccsb22}
 		\left ( \lambda _{q}^{+} \right )_{q}\subset \left ( \mu -\varepsilon ,\mu  \right )\ {\rm and}\ \left ( \lambda _{q}^{-} \right )_{q}\subset \left ( \mu  ,\mu+\varepsilon  \right ).
 	\end{equation}
From condition (C) in (II),  we have for any $q\in \mathbb{N} $,
 	\begin{equation}
 		\nu_{\tau}^{P_{n,m}}(\gamma _{\lambda_{q}^{+} })=\nu_{\tau}^{P_{n,m}}(\gamma _{\lambda_{q}^{-} })=0.
 	\end{equation}
 	Moreover,   as argued in the proof of (III),   we deduce for any $q\in \mathbb{N} $,
 	\begin{equation}
 		i_{\tau}^{P_{n,m}}(\gamma _{\lambda })=\left\{\begin{matrix}
 				i _{\tau}^{P_{n,m}}(\gamma _{\mu }),  \quad \quad \quad \quad \quad  \quad \quad  \quad\lambda =\lambda _{q}^{-},  \\
 				i _{\tau}^{P_{n,m}}(\gamma _{\mu })+\nu _{\tau}^{P_{n,m}}(\gamma _{\mu }),\quad \quad  \ \lambda =\lambda _{q}^{+}.
 		\end{matrix}\right.
 	\end{equation}
 	if \eqref{cccsb11} holds, and 
 	\begin{equation}
 		i_{\tau}^{P_{n,m}}(\gamma _{\lambda })=\left\{\begin{matrix}
 			i _{\tau}^{P_{n,m}}(\gamma _{\mu })+\nu _{\tau}^{P_{n,m}}(\gamma _{\mu }),\quad \quad  \quad \lambda =\lambda _{q}^{+}	,  \\
 		i _{\tau}^{P_{n,m}}(\gamma _{\mu }),  \quad \quad  \quad  \quad  \quad  \quad \quad  \quad  \ \lambda =\lambda _{q}^{-}.
 		\end{matrix}\right.
 	\end{equation}
 	if \eqref{cccsb22} holds. In both cases, since $\nu _{\tau}^{P_{n,m}}(\gamma _{\mu })=2k$ with $k\ge 1$, for any $q\in \mathbb{N}$, we deduce
 	\begin{equation}\label{IIp1}
 		\left [ i_{\tau}^{P_{n,m}}(\gamma _{\lambda_{q}^{-} }),i_{\tau}^{P_{n,m}}(\gamma _{\lambda_{q}^{-} })+\nu_{\tau}^{P_{n,m}}(\gamma _{\lambda_{q}^{-} }) \right ]\cap \left [ i_{\tau}^{P_{n,m}}(\gamma _{\lambda_{q}^{+} }),i_{\tau}^{P_{n,m}}(\gamma _{\lambda_{q}^{+} })+\nu_{\tau}^{P_{n,m}}(\gamma _{\lambda_{q}^{+} }) \right ] =\emptyset .
 	\end{equation}
Consequently, it follows from Theorem \ref{thm33}(II) that $(\mu,0)$ is a bifurcation point of system \eqref{ut}, which implies that $(\mu,0)$ is a a bifurcation point of system \eqref{e2}.
\end{proof}

\begin{proof}[\bf Proof of Corollary \ref{co1}]
For each $\lambda\ne0$, the time-rescaling $y(t)=x(\lambda t)$ transforms \eqref{e3} into 
\begin{equation}\label{e4}
	\begin{cases}
		\dot{y}(t)=-\lambda \nabla \mathbf{F} \left (  \int_{0}^{1} \lambda \mathbf{\Theta }  y\left ( t-s \right )   ds \right )- \cdots-\lambda \nabla \mathbf{F} \left (\int_{m-1  }^{m} \lambda \mathbf{\Theta }  y\left ( t-s \right )   ds \right ),\\
		y\left ( t+m \right ) =-y\left ( t \right ),\quad  t\in \mathbb{R}, \ y(t)\in \mathbb{R}^{n}.
	\end{cases}
\end{equation}
This rescaling induces a bijection
between the solution sets of the two systems, preserving nontriviality, $\mathbb R$-distinctness, and $C^1_{\rm loc}$-convergence to the trivial solution. In particular, for $\mu\ne0$, the two systems have the same bifurcation behavior at $(\mu,0)$. It therefore remains only to verify that Theorem \ref{thm2} is applicable
to \eqref{e4}.

Define
\[
F(\lambda,x):=\lambda \mathbf{F}(x),
\qquad
\Theta(\lambda):=\lambda\mathbf{\Theta}.
\]
Then assumption $({\bf FG})$ is satisfied. Moreover, the orthogonal matrix $T$ in \eqref{T2} satisfies
\begin{equation}
	T^{-1}\nabla ^{2}F(\lambda,0)T={\rm diag}\left ( \lambda f_{1},\cdots ,\lambda f_{n} \right )\ {\rm and}\ T^{-1}\Theta(\lambda) T={\rm diag}\left ( \lambda \theta _{1},\cdots ,\lambda \theta _{n} \right ).
\end{equation}
If $(\mu,0)$ is a bifurcation point of system \eqref{e3}, then,
by the correspondence established above, it is also a bifurcation
point of system \eqref{e2} with
\[
F(\lambda,x)=\lambda\mathbf F(x),\qquad
\Theta(\lambda)=\lambda\mathbf\Theta,
\qquad \tau=1.
\]
Therefore, Theorem \ref{thm2}(I) yields
\[
m|\mu|\sqrt{2f_j\theta_j}
\in (2\mathbb N-1)\pi,
\]
for some $j\in\{1,\cdots,n\}$. Clearly, condition (D)  coincides with condition (A) in Theorem \ref{thm2}(II), while condition (E) coincides with condition (B) therein. It therefore remains only to verify that either  condition (SC.1) or condition (SC.2) is satisfied for all $l=1,\cdots,k$.

Since the numbers $f_{j_1},\cdots,f_{j_k}$ have the same sign, there are two cases to consider. Suppose first that
\[
f_{j_l}>0,\qquad l=1,\cdots,k.
\]
Then, for each $l=1,\cdots,k$, if $\mu f_{j_l}>0$, we have
\begin{equation}
	\lambda^2 f_{j_l}\theta_{j_l}
	<
	\mu^2 f_{j_l}\theta_{j_l}
	<
	{\lambda'}^2 f_{j_l}\theta_{j_l},
	\qquad
	0<\lambda<\mu<\lambda',
\end{equation}
whereas, if $\mu f_{j_l}<0$, we have
\begin{equation}
	\lambda^2 f_{j_l}\theta_{j_l}
	>
	\mu^2 f_{j_l}\theta_{j_l}
	>
	{\lambda'}^2 f_{j_l}\theta_{j_l},
	\qquad
	\lambda<\mu<\lambda'<0.
\end{equation}
Thus condition (SC.1) is satisfied for all $l=1,\ldots,k$.

Suppose next that
\[
f_{j_l}<0,\qquad l=1,\cdots,k.
\]
Then, for each $l=1,\cdots,k$, if $\mu f_{j_l}>0$, we have
\begin{equation}
	\lambda^2 f_{j_l}\theta_{j_l}
	>
	\mu^2 f_{j_l}\theta_{j_l}
	>
	{\lambda'}^2 f_{j_l}\theta_{j_l},
	\qquad
	\lambda<\mu<\lambda'<0,
\end{equation}
whereas, if $\mu f_{j_l}<0$, we have
\begin{equation}
	\lambda^2 f_{j_l}\theta_{j_l}
	<
	\mu^2 f_{j_l}\theta_{j_l}
	<
	{\lambda'}^2 f_{j_l}\theta_{j_l},
	\qquad
	0<\lambda<\mu<\lambda'.
\end{equation}
Hence condition (SC.2) is satisfied for all $l=1,\cdots,k$.

Therefore, Theorem \ref{thm2}(III) is applicable to system
\eqref{e4}. By the correspondence established above, the
bifurcation alternatives obtained for \eqref{e4} carry over to
system \eqref{e3}. This completes the proof of Corollary \ref{co1}.
\end{proof}

\begin{proof}[\bf Proof of Theorem \ref{th4}] As established in Lemmas \ref{l3} and  \ref{l4}, for each $\lambda\in \Lambda$, $(\lambda,0)$ is a bifurcation point of system \eqref{e4'} if and only if $(\lambda,0)$ is a bifurcation point of the associated Hamiltonian system 
	\begin{equation}\label{3126}
		\left\{\begin{matrix}
			\dot{u}(t)=J_{2}\nabla \tilde{H}(\lambda ,u(t)),   \\
			u(t+\tau)=P_{1,2}u(t),\quad \
		\end{matrix}\right.
	\end{equation}
where the $4\times 4$ matrix $P_{1,2}$ is given by 	
\begin{equation}\label{eP12}
P_{1,2}=\begin{pmatrix}
		0&  -1& 0 &0 \\
		1& 0 & 0 &0 \\
		0& 0 & 0 & -1\\
		0& 0 & 1 &0
	\end{pmatrix} .
\end{equation}
Let $\gamma _{\lambda}(t)$ be the fundamental solution of the following Hamiltonian system  linearized at the  trivial solution $0\in \mathbb{R}^{4} $
\begin{equation}
	\dot{u} (t)=J_{2}\nabla ^{2}\tilde{H} (\lambda ,0) u(t), \quad \lambda\in \Lambda ,\quad u(t)\in \mathbb{R}^{4},
\end{equation}
with $\gamma _{\lambda }(0)=I_{4}$. Here, the Hessian $\nabla ^{2}\tilde{H} (\lambda ,0)$ has the block-diagonal form 
\begin{equation}
	\nabla ^{2}\tilde{H} (\lambda ,0)=\begin{pmatrix}
		\sqrt{2}g(\lambda)I_{2} & 0\\
		0 &\sqrt{2}f(\lambda)I_{2}
	\end{pmatrix},
\end{equation}
and thus $\gamma _{\lambda}(t)$ coincides with $\gamma _{	\sqrt{2}g(\lambda),\sqrt{2}f(\lambda)}$ defined in Lemma \ref{abv}. We next prove the three assertions separately.

	{\bf Step 1}(Proof of {\rm (I)}). 
Suppose that \((\mu,0)\) is a bifurcation point of system \eqref{e4'}. By the necessary condition in Theorem \ref{thm33}(I), the \(P_{1,2}\)-nullity of the linearized Hamiltonian system at \(\lambda=\mu\) is nonzero, that is,
\[
\nu_\tau^{P_{1,2}}(\gamma_\mu)\ne0.
\]
Applying Lemma \ref{abv} with \(m=2\), we obtain
\[
2\sqrt{2g(\mu)f(\mu)}\,\tau
\in(2\mathbb N-1)\pi.
\]
 Hence the necessary condition (I) follows.
 
 	{\bf Step 2}(Proof of {\rm (III)}). Suppose that the assumptions of (II) hold and that condition (d) is replaced by (sd). By condition (c), there exists \(j\in\mathbb N\) such that
 \begin{equation}\label{131}
 	 	2\sqrt{2g(\mu)f(\mu)}\,\tau=(2j-1)\pi.
 \end{equation}
Since \(g(\mu)f(\mu)>0\), both \(g(\mu)\) and \(f(\mu)\) are nonzero and have the same sign. By continuity, after shrinking \(\varepsilon>0\) if necessary, we may assume that \(g(\lambda)f(\lambda)>0\) and that \(g(\lambda)\) has the same sign as \(g(\mu)\) for every \(\lambda\in(\mu-\varepsilon,\mu+\varepsilon)\).

We first consider the case \(g(\mu)>0\). By Lemma \ref{abv} and \eqref{131}, we have
\begin{equation}
	\nu_\tau^{P_{1,2}}(\gamma_\mu)=2,
	\qquad
	i_\tau^{P_{1,2}}(\gamma_\mu)=2(j-1).
\end{equation}
Moreover, by further reducing \(\varepsilon\) and arguing similarly as \eqref{csb1} and \eqref{csb2},  we deduce
\begin{equation}\label{c132}
	i_\tau^{P_{1,2}}(\gamma_\lambda)=
	\begin{cases}
		2(j-1)=i_\tau^{P_{1,2}}(\gamma_\mu),&\lambda\in(\mu-\varepsilon,\mu),\\
		2j=i_\tau^{P_{1,2}}(\gamma_\mu)+\nu_\tau^{P_{1,2}}(\gamma_\mu),&\lambda\in(\mu,\mu+\varepsilon),
	\end{cases}
\end{equation}
if (sd.1) holds, 
\begin{equation}\label{c133}
	i_\tau^{P_{1,2}}(\gamma_\lambda)=
	\begin{cases}
		2j=i_\tau^{P_{1,2}}(\gamma_\mu)+\nu_\tau^{P_{1,2}}(\gamma_\mu),&\lambda\in(\mu-\varepsilon,\mu),\\
		2(j-1)=i_\tau^{P_{1,2}}(\gamma_\mu),&\lambda\in(\mu,\mu+\varepsilon),
	\end{cases}
\end{equation}
if (sd.2) holds, and in both cases, we have 
\begin{equation}\label{c134}
	\nu_\tau^{P_{1,2}}(\gamma_\lambda)=0,
	\qquad
	\text{for }0<|\lambda-\mu|<\varepsilon.
\end{equation}

Next, suppose that \(g(\mu)<0\). In this case, Lemma \ref{abv} yields
\begin{equation}\label{c135}
	\nu_\tau^{P_{1,2}}(\gamma_\mu)=2,
	\qquad
	i_\tau^{P_{1,2}}(\gamma_\mu)=-2j.
\end{equation}
Using the same argument, we obtain
\begin{equation}\label{c136}
	i_\tau^{P_{1,2}}(\gamma_\lambda)=
	\begin{cases}
		-2j+2=i_\tau^{P_{1,2}}(\gamma_\mu)+\nu_\tau^{P_{1,2}}(\gamma_\mu),&\lambda\in(\mu-\varepsilon,\mu),\\
		-2j=i_\tau^{P_{1,2}}(\gamma_\mu),&\lambda\in(\mu,\mu+\varepsilon),
	\end{cases}
\end{equation}
if (sd.1) holds, and
\begin{equation}\label{c137}
	i_\tau^{P_{1,2}}(\gamma_\lambda)=
	\begin{cases}
		-2j=i_\tau^{P_{1,2}}(\gamma_\mu),&\lambda\in(\mu-\varepsilon,\mu),\\
		-2j+2=i_\tau^{P_{1,2}}(\gamma_\mu)+\nu_\tau^{P_{1,2}}(\gamma_\mu),&\lambda\in(\mu,\mu+\varepsilon),
	\end{cases}
\end{equation}
if (sd.2) holds. Combining \eqref{c132}-\eqref{c134} with \eqref{c136} and \eqref{c137}, we conclude that assumption (b') of Theorem \ref{thm33}(III) is satisfied. Moreover, since $P_{1,2}^{2}=-I_4$, we have $\ker(P_{1,2}-I_4)={0}$, which verifies assumption (a') of the same theorem. 

Applying Theorem \ref{thm33}(III) to the Hamiltonian system \eqref{3126}, we obtain the alternatives stated in (III). By the correspondence established in Lemmas \ref{l3} and \ref{l4}, these alternatives carry over to system \eqref{e4'}. This completes the proof of (III).

	{\bf Step 3}(Proof of {\rm (II)}). The desired conclusion follows from a slight modification of the proof of Theorem \ref{n1}(II).
	\end{proof}

\begin{proof}[\bf Proof of Theorem \ref{thm5}] By Lemmas \ref{l3} and  \ref{l4}, system  \eqref{e3'} is equivalent, in the sense of preserving bifurcation, to the Hamiltonian system 
	\begin{equation}\label{3144}
		\left\{\begin{matrix}
			\dot{u}(t)=J_{2n}\nabla \tilde{H}(\lambda ,u(t)),   \\
			u(t+\tau)=P_{n,2}u(t),\quad \
		\end{matrix}\right.
	\end{equation}
	with  the matrix $P_{n,2}=\begin{pmatrix}
		Q_{n,2}&0 \\
		0&Q_{n,2}
	\end{pmatrix}$ defined in \eqref{Hp}.  Let \(\gamma_\lambda\) denote the fundamental solution of the linearized Hamiltonian system 
		\begin{equation}\label{ut211}
		\dot{u} (t)=J_{2n}\nabla ^{2}\tilde{H} (\lambda ,0) u(t), \quad \lambda\in \Lambda ,\quad u(t)\in \mathbb{R}^{4n},
	\end{equation}
	where Hessian matrix has the block-diagonal form
	$$\nabla ^{2}\tilde{H} (\lambda ,0)=\begin{pmatrix}
		I_{2}\otimes \left ( \sqrt{2}\nabla ^{2}G(\lambda ,0)  \right ) & 0\\
		0 &I_{2}\otimes \left ( \sqrt{2}\nabla ^{2}F(\lambda ,0)  \right )
	\end{pmatrix},$$
	and thus \(\gamma_\lambda\) coincides with $\gamma _{ \sqrt{2}\nabla ^{2}G(\lambda ,0),\sqrt{2}\nabla ^{2}F(\lambda ,0)}$ defined in Lemma \ref{abv2}. Next, we prove the three assertions separately.
	
		{\bf Step 1}(Proof of {\rm (I)}). Suppose that \((\mu,0)\) is a bifurcation point of system \eqref{e3'}. By the necessary bifurcation condition in Theorem \ref{thm33}(I), we have \(\nu_\tau^{P_{n,2}}(\gamma_\mu)\ne0\). Applying Lemma \ref{abv2} with \(m=2\), \(a_j=\sqrt{2}\,g_j(\mu)\), and \(b_j=\sqrt{2}\,f_j(\mu)\), we obtain
		\begin{equation}
		\nu_\tau^{P_{n,2}}(\gamma_\mu)
		=2\sharp \left \{ j\mid 2\sqrt{2g_j(\mu)f_j(\mu)} \tau \in  \left ( 2\mathbb{N}-1  \right ) \pi  \right \},
		\end{equation}
which  yields the necessary condition in (I).
	
	{\bf Step 2}(Proof of {\rm (III)}). By Lemmas \ref{l3} and  \ref{l4}, it suffices to verify the hypotheses of Theorem \ref{thm33}(III) for the Hamiltonian system \eqref{3144}. Since \(P_{n,2}^2=-I_{4n}\), we have \(\ker(P_{n,2}-I_{4n})=\{0\}\), so assumption (a') is satisfied.
	
	By condition (V) in (II), Lemma \ref{abv2} gives
	\[
	\nu_\tau^{P_{n,2}}(\gamma_\mu)=2k.
	\]
	For each \(l=1,\cdots,k\), we have 
	\[
	i_\tau^{P_{1,2}}
	\bigl(\gamma_{\sqrt{2}g_{j_l}(\mu),\,\sqrt{2}f_{j_l}(\mu)}\bigr)
	=
	\begin{cases}
		2(\rho_l(\mu)-1),&g_{j_l}(\mu)>0,\\
		-2\rho_l(\mu),&g_{j_l}(\mu)<0.
	\end{cases}
	\]
Let \(\Gamma\subset\{k+1,\cdots,n\}\) consist of those indices for which the  inequality \eqref{2r'} holds, and write \(\Gamma=\Gamma^+\cup\Gamma^-\), where
	\[
	\Gamma^\pm
	=\{l\in\Gamma\mid \pm g_{j_l}(\mu)>0\}.
	\]
	By conditions (U) and (V), together with Lemmas \ref{abv} and \ref{abv2}, the Maslov \(P_{n,2}\)-index at \(\mu\) is
\begin{equation}
	 		i _{\tau}^{P_{n,2}}(\gamma _{\mu })=\sum_{l\le k,g_{j_{l}}(\mu)>0}2\left ( \rho _{l}(\mu)-1 \right )+\sum_{l\le k,g_{j_{l}}(\mu)<0}	\left ( -2\rho _{l}(\mu)  \right )+\sum_{l\in\Gamma ^{+}}2\rho _{l}(\mu)+\sum_{l\in\Gamma ^{-}}\left ( -2\rho _{l}(\mu) \right ).
\end{equation}

Following the same argument as in the proof of Theorem \ref{thm2}(III) after \eqref{iPru}, we use the index formulas in Lemma \ref{abv} to determine the behavior of \(i_\tau^{P_{n,2}}(\gamma_\lambda)\) on either side of \(\mu\). In view of condition (SW), we obtain
 \begin{equation}
	i _{\tau}^{P_{n,2}}(\gamma _{\lambda })=\left\{\begin{matrix}
		i _{\tau}^{P_{n,2}}(\gamma _{\mu }), \quad \quad \quad \quad\quad\quad \ \forall \lambda \in \left ( \mu -\varepsilon ,\mu   \right ),\\
		i _{\tau}^{P_{n,2}}(\gamma _{\mu })+\nu _{\tau}^{P_{n,2}}(\gamma _{\mu }),\quad  \forall \lambda \in \left ( \mu ,\mu +\varepsilon    \right ),
	\end{matrix}\right.
\end{equation}
if (SW.1) holds, the two relations are reversed if (SW.2) holds. Thus, assumption (b') of Theorem \ref{thm33}(III) is satisfied.

	Applying Theorem \ref{thm33}(III) to \eqref{3144}, we obtain the alternatives in (III) for the Hamiltonian system. By  Lemmas \ref{l3} and  \ref{l4}, the corresponding alternatives hold for system \eqref{e3'}. Furthermore, when \(k>1\), we have \(\nu_\tau^{P_{n,2}}(\gamma_\mu)=2k>3\), so the second part of Theorem \ref{thm33}(III) apply. In particular, we obtain the inequality 
	\begin{equation}
			\int_{0}^{\tau } \left [ \frac{1}{2}\left ( J_{2n} \dot{\bar{u} }^{\lambda }(t),\bar{u} ^{\lambda }(t)  \right )+\tilde{H} (\lambda ,\bar{u}^{\lambda } (t))   \right ] dt \\
			\ne \int_{0}^{\tau } \left [ \frac{1}{2}\left ( J_{2n}\dot{\hat{u} }^{\lambda }(t),\hat{u} ^{\lambda }(t)  \right )+\tilde{H} (\lambda ,\hat{u}^{\lambda } (t))   \right ] dt ,
	\end{equation}	
	where $\bar{u} ^{\lambda } $ and $\hat{u} ^{\lambda } $ are two nonzero $\mathbb{R} $-distinct solutions of Hamiltonian system \eqref{3144} with $ \bar{u} ^{\lambda }(t)=S^{-1}\bar{z}^{\lambda } (t)$ and $ \hat{u} ^{\lambda }(t)=S^{-1}\hat{z}^{\lambda } (t)$. Then it follows from Lemma \ref{l4} that 
		\begin{equation}\label{c2}
				\int_{0}^{\tau } \left [ \frac{1}{2}\left ( \mathcal{J} \dot{\bar{z} }^{\lambda }(t),\bar{z} ^{\lambda }(t)  \right )+H(\lambda ,\bar{z}^{\lambda } (t))   \right ] dt \\
				\ne \int_{0}^{\tau } \left [ \frac{1}{2}\left ( \mathcal{J} \dot{\hat{z} }^{\lambda }(t),\hat{z} ^{\lambda }(t)  \right )+H(\lambda ,\hat{z}^{\lambda } (t))   \right ] dt ,   
		\end{equation}	
		where $\bar{z}^{\lambda }(t)=\left ( \bar{x} _{1}^{\lambda }(t),\bar{x} _{2}^{\lambda }(t),\bar{z} _{1}^{\lambda }(t),\bar{z} _{2}^{\lambda }(t)\right )^{T}$, $\bar{x} _{i}^{\lambda }(t)=\bar{x}^{\lambda }(t-(i-1)\tau)$, $\bar{z} _{i}^{\lambda }(t)= \int_{(i-1)\tau }^{i\tau }\nabla G\left (\lambda, \bar{x}^{\lambda } \left ( t-s \right ) \right )ds$, and $\hat{z}^{\lambda }(t)=\left ( \hat{x} _{1}^{\lambda }(t),\hat{x} _{2}^{\lambda }(t),\hat{z} _{1}^{\lambda }(t),\hat{z} _{2}^{\lambda }(t)\right )^{T}$, $\hat{x} _{i}^{\lambda }(t)=\hat{x}^{\lambda }(t-(i-1)\tau)$, $\hat{z} _{i}^{\lambda }(t)= \int_{(i-1)\tau }^{i\tau }\nabla G\left (\lambda, \hat{x}^{\lambda } \left ( t-s \right ) \right )ds$, $i=1,2$. By a direct computation, we have 
		\begin{equation}\label{126}
			\begin{split}
				\frac{1}{2}\int_{0}^{\tau }\left ( \mathcal{J} \dot{\bar{z} }^{\lambda }(t),\bar{z} ^{\lambda }(t)  \right )dt=&\frac{1}{2\sqrt{2} }\left ( \int_{0}^{\tau } \left ( \dot{\bar{z} }^{\lambda }_{1}(t)-\dot{\bar{z} }^{\lambda }_{2}(t),\bar{x} _{2}^{\lambda }(t) \right )dt -\int_{0}^{\tau } \left ( \dot{\bar{z} }^{\lambda }_{1}(t)+\dot{\bar{z} }^{\lambda }_{2}(t),\bar{x} _{1}^{\lambda }(t) \right )dt  \right )\\
				&+ \frac{1}{2\sqrt{2} }\left ( \int_{0}^{\tau } \left ( \dot{\bar{x} }^{\lambda }_{1}(t)-\dot{\bar{x} }^{\lambda }_{2}(t),\bar{z} _{1}^{\lambda }(t) \right )dt +\int_{0}^{\tau } \left ( \dot{\bar{x} }^{\lambda }_{1}(t)+\dot{\bar{x} }^{\lambda }_{2}(t),\bar{z} _{2}^{\lambda }(t) \right )dt  \right ).
			\end{split}
		\end{equation}
		Using integration by parts, we obtain from the  symmetry $	 \bar{x} ^{\lambda }\left ( t+2\tau  \right ) =- \bar{x} ^{\lambda }\left ( t \right )=- \bar{x}_{1} ^{\lambda }\left ( t \right )$ that 
		\begin{equation}\label{127}
			\begin{split}
				&\int_{0}^{\tau } \left ( \dot{\bar{x} }^{\lambda }_{1}(t)-\dot{\bar{x} }^{\lambda }_{2}(t),\bar{z} _{1}^{\lambda }(t) \right )dt +\int_{0}^{\tau } \left ( \dot{\bar{x} }^{\lambda }_{1}(t)+\dot{\bar{x} }^{\lambda }_{2}(t),\bar{z} _{2}^{\lambda }(t) \right )dt\\
				=& \int_{0}^{\tau } \left ( \dot{\bar{z} }^{\lambda }_{1}(t)-\dot{\bar{z} }^{\lambda }_{2}(t),\bar{x} _{2}^{\lambda }(t) \right )dt -\int_{0}^{\tau } \left ( \dot{\bar{z} }^{\lambda }_{1}(t)+\dot{\bar{z} }^{\lambda }_{2}(t),\bar{x} _{1}^{\lambda }(t) \right )dt\\
				=&-2\left ( \int_{0}^{\tau } \left ( \nabla G\left (\lambda, \bar{x}^{\lambda }_{2} \left ( t\right ) \right ),\bar{x} _{2}^{\lambda }(t) \right )dt -\int_{0}^{\tau } \left ( \nabla G\left (\lambda, \bar{x}^{\lambda }_{1} \left ( t\right ) \right ),\bar{x} _{1}^{\lambda }(t) \right )dt  \right )\\
				=& -2\int_{0}^{2\tau } \left ( \nabla G\left (\lambda, \bar{x}^{\lambda } \left ( t\right ) \right ),\bar{x}^{\lambda }(t) \right )dt,
			\end{split}
		\end{equation}
		where the second equality holds since
		$$  \dot{\bar{z} }^{\lambda }_{1}(t)-\dot{\bar{z} }^{\lambda }_{2}(t)=-2\nabla G\left (\lambda, \bar{x}^{\lambda }_{2} \left ( t\right ) \right )\ {\rm and}\  \dot{\bar{z} }^{\lambda }_{1}(t)+\dot{\bar{z} }^{\lambda }_{2}(t)=2\nabla G\left (\lambda, \bar{x}^{\lambda }_{1} \left ( t\right ) \right ).$$
		\eqref{127} together with \eqref{126} yields 
		\begin{equation}\label{128}
			\frac{1}{2}\int_{0}^{\tau }\left ( \mathcal{J} \dot{\bar{z} }^{\lambda }(t),\bar{z} ^{\lambda }(t)  \right )dt=-\sqrt{2} \int_{0}^{2\tau } \left ( \nabla G\left (\lambda, \bar{x}^{\lambda } \left ( t\right ) \right ),\bar{x}^{\lambda }(t) \right )dt.
		\end{equation}
		Moreover, since the functions $F(\lambda ,\cdot )$ and $G(\lambda ,\cdot )$ are even, we have 
		\begin{equation}\label{129}
			\begin{split}
				\int_{0}^{\tau } H(\lambda ,\bar{z}^{\lambda } (t))dt=&\sqrt{2} \left ( \int_{0}^{\tau }G(\lambda , \bar{x} _{1}^{\lambda }(t))dt+\int_{0}^{\tau } G(\lambda , \bar{x} _{2}^{\lambda }(t))dt+\int_{0}^{\tau } F(\lambda , \bar{z} _{1}^{\lambda }(t))dt+\int_{0}^{\tau } F(\lambda , \bar{z} _{2}^{\lambda }(t))dt \right )\\
				=& \sqrt{2} \int_{0}^{2\tau }G(\lambda , \bar{x} ^{\lambda }(t))dt+ \sqrt{2} \int_{0}^{2\tau }F(\lambda , \bar{z} _{1}^{\lambda }(t))dt\\
				=& \sqrt{2} \int_{0}^{2\tau }G(\lambda , \bar{x} ^{\lambda }(t))dt+ \sqrt{2} \int_{0}^{2\tau }F\left ( \lambda , \int_{0}^{\tau }\nabla G\left (\lambda, \bar{x}^{\lambda } \left ( t-s \right ) \right ) ds \right ) dt.
			\end{split}
		\end{equation}
		Combined with \eqref{128} and \eqref{129}, it follows that the inequality \eqref{124} holds, which completes the proof of Theorem \ref{thm5}(III).

{\bf Step 3}(Proof of {\rm (II)}). The desired conclusion follows by a slight modification of the proof of Theorem \ref{thm2}(II). 
	\end{proof}

\begin{proof}[\bf Proof of Corollary \ref{co2}] For each $\lambda\ne0$, the time rescaling $y(t)=x(\lambda t)$ transforms \eqref{e31} into 
	\begin{equation}\label{e41}
		\begin{cases}
			\dot{y}(t)=-\lambda \nabla \mathbf{F} \left (  \int_{0}^{1} \lambda \nabla \mathbf{G}  \left ( y\left ( t-s \right ) \right )   ds \right )-\lambda \nabla \mathbf{F} \left (\int_{1  }^{2} \lambda \nabla \mathbf{G}  \left ( y\left ( t-s \right ) \right )   ds \right ),\\
			y\left ( t+2 \right ) =-y\left ( t \right ),\quad  t\in \mathbb{R}, \ y(t)\in \mathbb{R}^{n},
		\end{cases}
	\end{equation}
which is precisely system \eqref{e3'} with \(\tau=1\) and potentials $F(\lambda ,\cdot )=\lambda \mathbf{F}(\cdot )$ and $G(\lambda ,\cdot )=\lambda \mathbf{G}(\cdot)$. The time rescaling establishes a bijection between the corresponding solution sets, preserving nontriviality, $\mathbb R$-distinctness, and convergence to the trivial solution in \(C_{\mathrm{loc}}^1\). Hence, for every \(\mu\neq0\), the bifurcation behavior of system \eqref{e31} at \((\mu,0)\) is equivalent to that of system \eqref{e41}. 
	
	From \eqref{T21}, the orthogonal matrix $T$ satisfies 
\begin{equation}
	T^{-1}\nabla ^{2}F(\lambda ,0)T={\rm diag}\left ( \lambda f_{1},\cdots ,\lambda f_{n} \right )\ {\rm and}\ T^{-1}\nabla ^{2}G(\lambda ,0 ) T={\rm diag}\left (\lambda g _{1},\cdots ,\lambda g _{n} \right ),
\end{equation}	
	and thus, in the notation of Theorem \ref{thm5}, the corresponding eigenvalue functions are
	\[
	f_j(\lambda)=\lambda f_j,\qquad
	g_j(\lambda)=\lambda g_j,
	\quad j=1,\ldots,n,
	\]
	where \(f_j\) and \(g_j\) on the right-hand side denote the constant eigenvalues in  \eqref{T21}. In particular,
	\begin{equation}\label{3157}
	g_j(\lambda)f_j(\lambda)=\lambda^2g_jf_j.
	\end{equation}

	Suppose first that \((\mu,0)\) is a bifurcation point of system \eqref{e31}. By the preceding equivalence, it is also a bifurcation point of system \eqref{e41}. The necessary condition in Theorem \ref{thm5}(I), applied with \(\tau=1\), therefore gives
	\[
	2\sqrt{2g_j(\mu)f_j(\mu)}
	=2|\mu|\sqrt{2g_jf_j}
	\in(2\mathbb N-1)\pi,
	\]
	for some \(j\in\{1,\cdots,n\}\). This proves the necessary condition.
	
	Conversely, suppose that \(\mu\neq0\) and that conditions (X) and (Y) hold. By \eqref{3157}, condition (X) is exactly condition (U) of Theorem \ref{thm5}(II) for \(\tau=1\), while condition (Y) is exactly condition (V). It remains to verify the strengthened crossing condition (SW) in Theorem \ref{thm5}(III).
	
	Since \(g_{j_1},\cdots,g_{j_k}\) have the same sign, we distinguish two cases.
	
	{\bf Case 1:}\ \(g_{j_l}>0\) for \(l=1,\cdots,k\).  Since \(g_{j_l}f_{j_l}>0\), we also have \(f_{j_l}>0\) for every \(l\). If \(\mu g_{j_l}>0\), then for any \(\lambda<\mu<\lambda'\) sufficiently close to \(\mu\),
	\[
	\lambda^2g_{j_l}f_{j_l}
	<
	\mu^2g_{j_l}f_{j_l}
	<
	\lambda'^2g_{j_l}f_{j_l}.
	\]
	On the other hand, if \(\mu g_{j_l}<0\), then for any \(\lambda<\mu<\lambda'\) sufficiently close to \(\mu\),
	\[
	\lambda^2g_{j_l}f_{j_l}
	>
	\mu^2g_{j_l}f_{j_l}
	>
	\lambda'^2g_{j_l}f_{j_l},
	\]
 hence (SW.1) is  satisfied.
	
		{\bf Case 2:}\  \(g_{j_l}<0\) for \(l=1,\cdots,k\).  If \(\mu g_{j_l}>0\), then for any \(\lambda<\mu<\lambda'\) sufficiently close to \(\mu\),
	\[
	\lambda^2g_{j_l}f_{j_l}
>
	\mu^2g_{j_l}f_{j_l}
	>
	\lambda'^2g_{j_l}f_{j_l}.
	\]
	On the other hand, if \(\mu g_{j_l}<0\), then for any \(\lambda<\mu<\lambda'\) sufficiently close to \(\mu\), 
	\[
	\lambda^2g_{j_l}f_{j_l}
	<
	\mu^2g_{j_l}f_{j_l}
	<
	\lambda'^2g_{j_l}f_{j_l},
	\]
hence (SW.2) is  satisfied.
	
	Therefore, Theorem \ref{thm5}(III) is applicable to system \eqref{e41}. By the correspondence established above, the
	bifurcation alternatives obtained for \eqref{e41} carry over to
	system \eqref{e31} . This completes the proof of Corollary \ref{co2}.
	\end{proof}

\end{document}